\documentclass[]{elsarticle}

\usepackage{lineno,hyperref}

\usepackage{graphicx}
\usepackage{bm}
\usepackage{multirow}
\usepackage{amsmath}
\usepackage{amsfonts}
\usepackage{amssymb}
\usepackage{amsthm}
\usepackage{secdot}

\usepackage{multicol}
\usepackage{booktabs}

\journal{Applied Numerical Mathematics}

\theoremstyle{plain}

\theoremstyle{definition}

\theoremstyle{remark}

\numberwithin{equation}{section}
\numberwithin{theorem}{section}
\numberwithin{remark}{section}

\newcommand{\porder}{p}
\newcommand{\qorder}{\mathcal{P}}

\newcommand{\PP}{{P}}
\newcommand{\vv}{{v}}
\newcommand{\ver}{{{\rm vert}}}

\begin{document}

\begin{frontmatter}

\title{Foundations of Space-Time Finite Element Methods: Polytopes, Interpolation, and Integration}

\author{Cory V. Frontin}
\address{Department of Aeronautics and Astronautics, Massachusetts Institute of Technology, Cambridge, Massachusetts 02139}

\author{Gage S. Walters} 
\address{Department of Mechanical Engineering, The Pennsylvania State University, University Park, Pennsylvania 16802}

\author{Freddie D. Witherden} 
\address{Department of Ocean Engineering, Texas A\&M University, College Station, Texas 77843}

\author{Carl W. Lee} 
\address{Department of Mathematics, University of Kentucky, Lexington, Kentucky 40506}

\author{David M. Williams \corref{mycorrespondingauthor}} 
\address{Department of Mechanical Engineering, The Pennsylvania State University, University Park, Pennsylvania 16802}

\cortext[mycorrespondingauthor]{Corresponding author}
\ead{david.m.williams@psu.edu}

\author{David L. Darmofal}
\address{Department of Aeronautics and Astronautics, Massachusetts Institute of Technology, Cambridge, Massachusetts 02139}

\begin{abstract}
The main purpose of this article is to facilitate the implementation of space-time finite element methods in four-dimensional space. 
In order to develop a finite element method in this setting, it is necessary to create a numerical foundation, or equivalently a numerical infrastructure. This foundation should include a collection of suitable elements (usually hypercubes, simplices, or closely related polytopes), numerical interpolation procedures (usually orthonormal polynomial bases), and numerical integration procedures (usually quadrature rules). It is well known that each of these areas has yet to be fully explored, and in the present article, we attempt to directly address this issue. We begin by developing a concrete, sequential procedure for constructing generic four-dimensional elements (4-polytopes). Thereafter, we review the key numerical properties of several canonical elements: the tesseract, tetrahedral prism, and pentatope. Here, we provide explicit expressions for orthonormal polynomial bases on these elements. Next, we construct fully symmetric quadrature rules with positive weights that are capable of exactly integrating high-degree polynomials, e.g.~up to degree 17 on the tesseract. Finally, the quadrature rules are successfully tested using a set of canonical numerical experiments on polynomial and transcendental functions.
\end{abstract}

\begin{keyword}
space-time; finite element methods; quadrature; tesseract; tetrahedral prism; pentatope 
\MSC[2010] 52B11, 65D05, 65D32, 74S05, 76M10 
\end{keyword}

\end{frontmatter}

\section{Introduction} \label{sec;introduction}

Space-time finite element methods are a promising approach for producing highly accurate and efficient solutions to partial differential equations (PDEs). More specifically, these methods can adaptively redistribute their degrees of freedom in both space and time in order to optimally resolve critical features of the solution.  Despite their potential, the application of these methods to real-world engineering problems remains an open area of research. While one may simulate space-time problems in 1$d+t$ (two-dimensional space) or 2$d+t$ (three-dimensional space) by leveraging well-established finite element techniques, many of these techniques do not carry over to 3$d+t$ (four-dimensional space). In some sense, four-dimensional space should be approached as a new frontier, and we must recognize the unique geometric characteristics which arise in this context. 

There are several specific challenges for space-time finite element methods in the areas of meshing, interpolation, and numerical integration (i.e.~quadrature/cubature). In what follows, we will briefly review these challenges, along with a short discussion of the associated literature. Our focus will be limited to state-of-the-art challenges associated with polynomial-based finite element methods. Note: for a more fundamental discussion of the methods themselves, including their basic construction, one may consult the textbooks of Ern and Guermond~\cite{Ern13}, and Johnson~\cite{johnson2012numerical}, (see Chapters 6 and 8, respectively). In addition, for a discussion of non-uniform rational B-spline (NURB)-based methods, please consult the recent work on isogeometric analysis (IGA)~\cite{langer2016space,takizawa2016ram,takizawa2017turbocharger,otoguro2017space,takizawa2017heart,takizawa2018stabilization,otoguro2019turbocharger,langer2019adaptive,terahara2020ventricle}.


Our first challenge is to develop unstructured space-time meshes in order to simulate problems with complex boundaries. In this setting, one may extend the traditional notion of an unstructured mesh by defining `fully-unstructured' meshes as grids which are unstructured in both space and time, and `partially-unstructured' meshes as grids which are unstructured in space only. The latter meshes are usually generated by extruding unstructured spatial meshes in the temporal direction. Evidently, the partially-unstructured meshes are less optimal than the fully-unstructured meshes, as all elements are required to take the same time-step. Nevertheless, fully-unstructured meshes are more difficult to generate. The problem of mesh generation may become significantly easier if hybrid meshes are employed, however then it becomes necessary to select a suitable conforming set of elements (4-polytopes) for populating the mesh. In three-dimensional space, one may select combinations of hexahedron, triangular prism, square pyramid, and tetrahedron elements. However, in four-dimensional space there are significantly more elements to consider (as we will demonstrate). In order to avoid these complications, previous researchers have used non-hybrid meshes composed entirely by one of the following standard elements: tesseracts, tetrahedral prisms, or pentatopes. For example, partially-unstructured meshes of tesseract elements were used by van der Vegt and coworkers~\cite{van2002spaceOne,van2002spaceTwo,klaij2006space} along with discontinuous Galerkin (DG) methods to solve laminar compressible flow problems. In addition, Diosady and Murman~\cite{diosady2015higher,diosady2017tensor,diosady2018linear,diosady2019scalable} have used a similar approach to solve turbulent compressible flow problems. In these studies, the primary objective has been to exploit the tensor-product structure of the tesseract and maximize computational efficiency. A contrasting approach is exemplified by Tezduyar and coworkers~\cite{tezduyar2006space,tezduyar2007modelling,tezduyar2010space,tezduyar2019space}, who used partially-unstructured meshes of tetrahedral prism elements along with Galerkin least-squares (GLS) in space/DG in time methods to solve fluid-structure interaction (FSI) problems. In these studies, the primary objective was to leverage the flexibility of unstructured tetrahedral grids for treating complicated geometries. Yang et al.~\cite{yang2020unstructured} have used a similar approach along with a hybrid finite-difference/finite-element method to solve the \emph{fractional} Bloch-Torrey equations for biomedical flows. (Note: for detailed discussions of fractional PDEs, including the Fredholm, Schr{\"o}dinger, and shallow water equations, one may consult the writings of Al-Smadi et al.~\cite{al2017numerical,al2019computational,al2020attractive,al2020numerical}, and Karniadakis and coworkers~\cite{karniadakis2015special,kharazmi2020fractional}). In an alternative fashion, Behr~\cite{behr2008simplex} has shown that tetrahedral prisms can be subdivided in order to generate pentatope elements. The resulting meshes have been used by Behr and coworkers~\cite{karyofylli2018simplex,karyofylli2019simplex,von2019simplex} in conjunction with GLS in space/DG in time methods to solve two-phase incompressible flow problems and single-phase compressible flow problems. In addition, Wang~\cite{wang2015discontinuous} has used these meshes in conjunction with space-time DG methods to solve compressible, viscous flow problems with rotating bodies. Furthermore, Lehrenfeld~\cite{lehrenfeld2015nitsche} has used these meshes in conjunction with enriched space-time DG methods (XFEM-DG methods) to solve  convection-diffusion problems with moving bodies. These efforts (above) have all focused on partially-unstructured meshes. In order to move towards fully-unstructured meshes, Mont~\cite{mont2012adaptive} used a tent-pitching algorithm~\cite{ungor2000tent} to generate unstructured pentatopic meshes that locally satisfy a causality constraint. In addition, Neum{\"u}ller and Steinbach~\cite{neumuller2011refinement,neumuller2013space} generated isotropic refinement strategies for fully-unstructured meshes of pentatope elements. Similar work has been subsequently performed by Foteinos et al.~\cite{foteinos20154d} and Belda-Ferr{\'\i}n et al.~\cite{belda2018local}. Most recently, Caplan and coworkers~\cite{caplan2017anisotropic,caplan2019extension,caplan2019four,caplan2020four} have developed highly-anisotropic, fully-unstructured meshes of pentatope elements. In principle, these meshes possess the greatest flexibility of any meshes generated to-date for discretizing space-time problems. Despite this flexibility, there are still opportunities to further improve the generality of space-time meshes by considering a broader class of elements beyond the classical tesseract, tetrahedral prism, and pentatope elements; and furthermore by constructing hybrid meshes (as discussed earlier).

The next challenge is to develop a family of modal basis functions and/or nodal basis functions. For example, one may construct a modal basis on the reference tesseract $[-1,1]^4$ by taking tensor products of one-dimensional Legendre polynomials, and a nodal basis can be constructed in a similar fashion by using Lagrange polynomials. Unfortunately, this simplicity-of-construction does not carry over to other element types. For example, the modal basis on the tetrahedral prism is usually obtained from a tensor product of the one-dimensional Legendre basis and the three-dimensional PKDO basis (originally derived in two dimensions by Proriol~\cite{proriol1957famille}, Koornwinder~\cite{koornwinder1975two}, Dubiner~\cite{dubiner1991spectral}, and Owens~\cite{owens1998spectral}). Once this modal basis is obtained, a nodal basis can be numerically computed using a generalized Vandermonde matrix, following the procedure on p.~410 of~\cite{Hesthaven07} and p.~124 of~\cite{karniadakis2013spectral}. However, the resulting nodal basis functions do not have explicit expressions, and the conditioning of the Vandermonde matrix can be sensitive to the polynomial order and the placement of  nodal points. In practice, the latter issue is not especially serious for low-to-moderate polynomial orders ($\porder \leq 6$)~\cite{Hesthaven07}. Note that it is possible to completely avoid this issue by explicitly constructing nodal basis functions using the Lagrange-mapping procedure of Hughes~(\cite{hughes2012finite} pp.~164-174). However, this procedure is not general, and must be tailored for each polynomial order. Therefore, the Vandermonde procedure of~\cite{Hesthaven07} is usually preferred. A similar set of issues arise for pentatope basis functions. In this case, the modal basis is chosen to be the four-dimensional PKDO basis (defined in~\cite{warburton2003constants}), and the nodal basis can be obtained using the aforementioned Vandermonde procedure. 
In summary, the construction of modal and nodal basis functions for tesseract, tetrahedral prism, and pentatope elements is well-established, with some caution required for nodal basis functions at high polynomial orders. Nevertheless, there does not appear to be a well-known resource which aggregates these basis function definitions in one place. In addition, there does not appear to be any research on the construction of modal or nodal basis functions for non-standard element types.

The last challenge is to numerically integrate polynomial functions and transcendental functions on each element. One can easily construct quadrature rules for the standard tesseract element by taking tensor products of one-dimensional Gauss-Legendre rules. The resulting rules are `good' in the following sense: i) they exactly integrate polynomials of degree $\leq \porder$, where $\porder = 2M-1$ and $M$ is the number of quadrature points in one dimension; ii) all quadrature points reside strictly within the tesseract; iii) all quadrature points are arranged symmetrically; and iv) all quadrature weights are positive. Unfortunately, while tensor-product rules are `good', they are not optimal as they possess an excessive number of quadrature points ($N_p = M^4$ points). Furthermore, they cannot be mapped to tetrahedral prism or pentatope elements without violating symmetry. As a result, it is important to generate good, non-tensor-product rules which are optimized for each element type. In what follows, we briefly review the best known quadrature rules of this kind on the tesseract and pentatope. Note: to the authors' knowledge, no such rules have been constructed for the tetrahedral prism. On the tesseract, Stroud~\cite{stroud1971approximate}, and S{\o}revik and Espelid~\cite{sorevik1987fully} constructed rules which exactly integrate polynomials of degrees $\porder = 3$ and 5, with 8 and 24 points, respectively. 
Thereafter, Majorana et al.~\cite{majorana1982shortened}, and S{\o}revik and Espelid~\cite{sorevik1989fully} constructed a rule which exactly integrates a polynomial of degree $\porder = 7$ with 57 points. Finally, S{\o}revik~\cite{sorevik1988reliable} discovered a rule which exactly integrates a polynomial of degree $\porder = 9$ with 160 points. In a similar fashion, many researchers have constructed quadrature rules for the pentatope (see the recent review in~\cite{williams2020family}). In particular, Hammer and Stroud~\cite{hammer1956numerical,stroud1966some} constructed rules which exactly integrate polynomials of degrees $\porder = 2$ and 3, with 5 and 15 points, respectively (although the latter rule has points on the boundary). Thereafter, Gusev et al.~\cite{gusev2018symbolic} have constructed rules which exactly integrate polynomials of degrees $4,5,6,7,$ and 8 with 20, 30, 56, 76, and 110 points, respectively. In summary, there is a critical absence of good quadrature rules capable of exactly integrating polynomials with $\porder \gtrsim 9$ on tesseract and pentatope elements. Furthermore, only tensor-product based rules are known to exist on the tetrahedral prism. 

The purpose of this article is to address some of the key shortcomings discussed above. In particular, we will: i) identify new types of space-time elements, and ii) review/expand the interpolation and integration procedures associated with tesseract, tetrahedral prism, and pentatope elements. The format of this article is as follows. In section 2, we describe a systematic procedure for generating sequences of space-time elements. In section 3, we define the standard tesseract, tetrahedral prism, and pentatope reference elements, and review their essential properties, including their symmetry groups and orthonormal basis functions. In section 4, we describe a computational strategy for constructing quadrature rules on these standard elements. Thereafter, in sections 5 and 6, we summarize the resulting quadrature rules, and demonstrate their effectiveness on a suite of numerical experiments. Finally, in section 7, we summarize the main conclusions of the article.

\section{Identification of Element Sequences} \label{sec;elements}

In this section, our objective is to construct sequences of space-time elements. Evidently, there are an infinite number of possible sequences, and within each sequence there are an unlimited number of polytopes which can be considered. However, for the sake of convenience, we narrow our focus to finite sequences of convex polytopes which can be formed by degenerating the vertices of the tesseract. These polytopes are most amenable to practical applications, as in principle, a polynomial basis can be formed on each element by constructing a suitable mapping from the tesseract to the associated subpolytope. In what follows, we will generate sequences of elements in this fashion, and introduce formal mathematical notation and terminology to aid in this process.  

The unit $d$-cube ($d$-dimensional cube) is the $d$-polytope  $I^d$, the $d$-fold prism of the unit interval $I=[0,1]$.  
Its vertices are all points of the form 
$(x_1,\ldots,x_d)$ such that each $x_i$ equals either 0 or 1.
A 0/1-polytope $\PP$ is the convex hull of a nonempty subset $V$ of the vertices of the $d$-cube; in other words, the vertex set of $\PP$, $\ver(\PP)$, will be denoted by $V$.
If $V$ consists of the points $(0,\ldots,0)$, $(1,0,\ldots,0)$, $(0,1,0,\ldots,0)$, \ldots, $(0,0,\ldots,0,1)$, we will call $\PP$ the standard $d$-simplex.
See~\cite{ziegler,aichholzer2000extremal,baumeister2009permutation,cihangir20150} for a description of some of the properties of 0/1-polytopes.

In the remainder of this section, we describe the construction of two sequences of $d$-elements ($d$-dimensional elements). These sequences are noteworthy, as for $d =2$ and $d =3$ they contain standard elements which are frequently used to construct 2D and 3D meshes. Naturally, here we are interested in the generalization of these sequences to 4D, and possibly higher dimensions.

First, let us consider the construction of the following sequence
\[
\PP_0,\PP_1,\ldots,\PP_{2^d-2d},\ldots,\PP_{2^d-d-1}
\]
together with a sequence of 0/1-points
\[
\vv_0,\ldots,\vv_{2^d-d-2}
\]
with the following properties:
\begin{itemize}
\item The element $\PP_i$ is a $d$-dimensional 0/1-polytope, $i=0,\ldots,2^d-d-1$.
\item The point $\vv_i$ is a vertex of $\PP_i$, $i=0,\ldots,2^d-d-2$.
\item The element $\PP_{i+1}$ is the convex hull of $\ver(\PP_i)\setminus\vv_i$,  $i=0,\ldots,2^d-d-2$.  Thus $\PP_{i+1}$ is the result of degenerating the vertex $\vv_i$ of $\PP_i$.
\item The first element of the sequence $\PP_0$ is the unit $d$-cube.
\item The intermediate element of the sequence $\PP_{2^d-2d}$ is a prism formed by extruding the $(d-1)$-simplex into one dimension higher. In other words, the element $\PP_{2^d-2d}$ is a prism over the standard $(d-1)$-simplex.
\item The last element of the sequence $\PP_{2^d-d-1}$ is the standard $d$-simplex.
\end{itemize}

Note in particular that when $d=4$, then the $3$-simplex is a $3$-dimensional tetrahedron, and  $\PP_8$ is the tetrahedral prism obtained by extruding the tetrahedron in the direction of the fourth  coordinate.  

The construction of the sequence is recursive on $d=1,2,\ldots$.
Naturally, when $d=1$, $\PP_0=\PP_{2^d-2d}=\PP_{2^d-d-1}=[0,1]$ and the vertex sequence is empty. When $d=2$, $\vv_0=(1,1)$, $\PP_0=\PP_{2^d-2d}$ (the unit square), and $\PP_1=\PP_{2^d-d-1}$ (the standard $2$-simplex). For $d\geq3$, the vertex sequence $\vv'_i$ in dimension $d$ is constructed from the vertex sequence $\vv_0,\ldots,\vv_{2^{d-1}-(d-1)-2}$ in dimension $d-1$ by using the following steps.
\begin{enumerate}
\item For each vertex $\vv_i$ create two new vertices $\vv'_{2i}=(\vv_i,1)$ and $\vv'_{2i+1}=(\vv_i,0)$ by appending to $\vv_i$ first a 1 and second a 0, for $i=0,\ldots,2^{d-1}-(d-1)-2$.
The result is $\vv'_0,\ldots,\vv'_{2^{d}-2d-1}$.
\begin{sloppypar}
\item Complete the sequence by setting $\vv'_{2^{d}-2d}=(1,0,0,0,\ldots,0,1)$,
$\vv'_{2^{d}-2d+1}=(0,1,0,0,\ldots,0,1)$,
$\vv'_{2^{d}-2d+2}=(0,0,1,0,\ldots,0,1)$,
\ldots,$\vv'_{2^{d}-d-2}=(0,0,0,0,\ldots,1,1)$.
\end{sloppypar}
\end{enumerate}
If we regard the $d$-cube as a prism over the $(d-1)$-cube, the above sequence of degenerations alternates by removing a vertex of the `upper' $(d-1)$-cube, followed by removing the corresponding vertex of the `lower' $(d-1)$-cube, thus replicating the $(d-1)$-dimensional degeneration sequence on both the top and the bottom.  At the end of phase~1 above, we will have a prism over the standard $(d-1)$-simplex.  
(In fact, $\PP'_{2i}$ will be a prism over $\PP_i$, $i=0,\ldots,2^{d-1}-d-1$\@.)
From that point on, the sequence in phase~2  strips away vertices on the top until a single point remains, resulting in a standard $d$-simplex.

Table~\ref{vertexsequences} shows the vertex sequences for dimensions 2, 3, and 4.  (For convenience we use an abbreviated form of the vertex coordinates\@.)
\begin{table}[h!]
\centering
\[
\begin{array}{ccc}
d=2&d=3&d=4\\
\hline
11&111&1111\\
&110&1110\\
&101&1101\\
&011&1100\\
&&1011\\
&&1010\\
&&0111\\
&&0110\\
&&1001\\
&&0101\\
&&0011
\end{array}
\]
\caption{Sequences of vertices in low dimensions}
\label{vertexsequences}
\end{table}

\begin{figure}[h!]
\centering
\begin{tabular}{ccc}
\includegraphics[scale=.45]{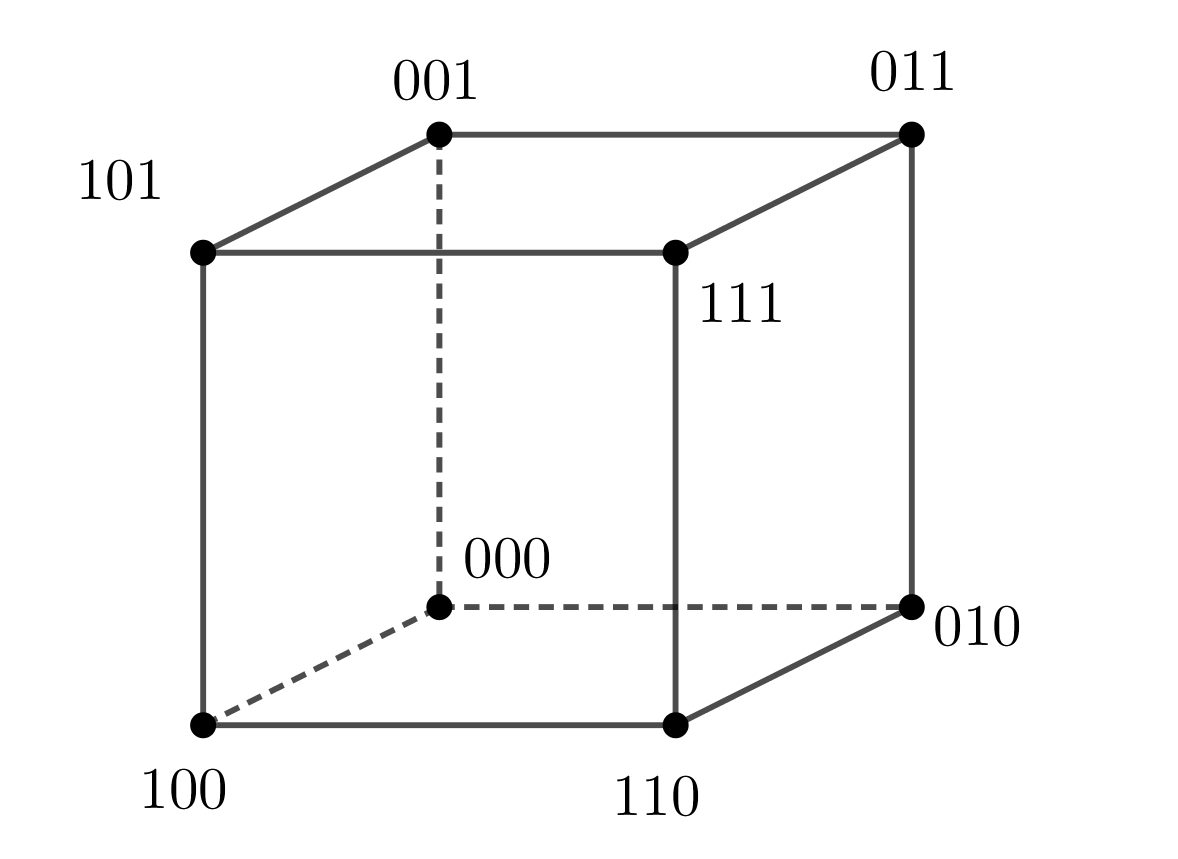}&
\includegraphics[scale=.45]{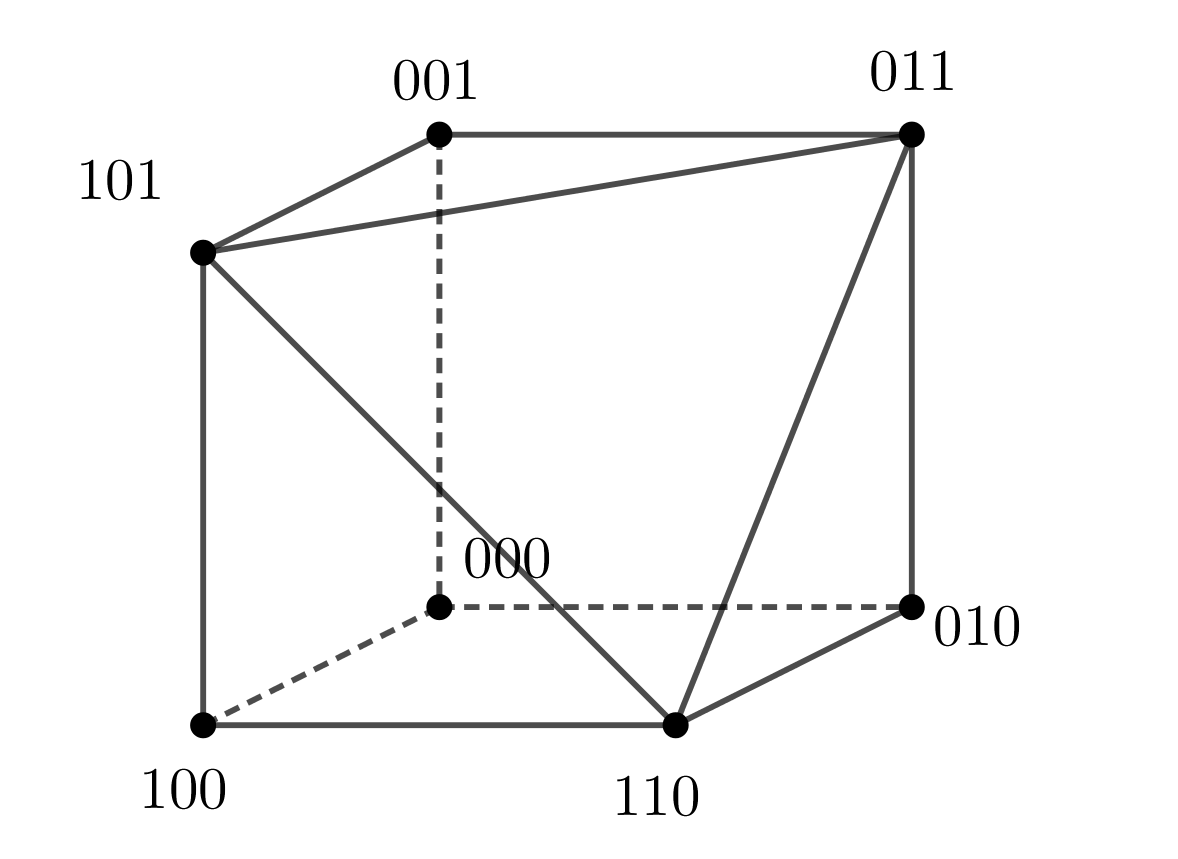}&
\includegraphics[scale=.45]{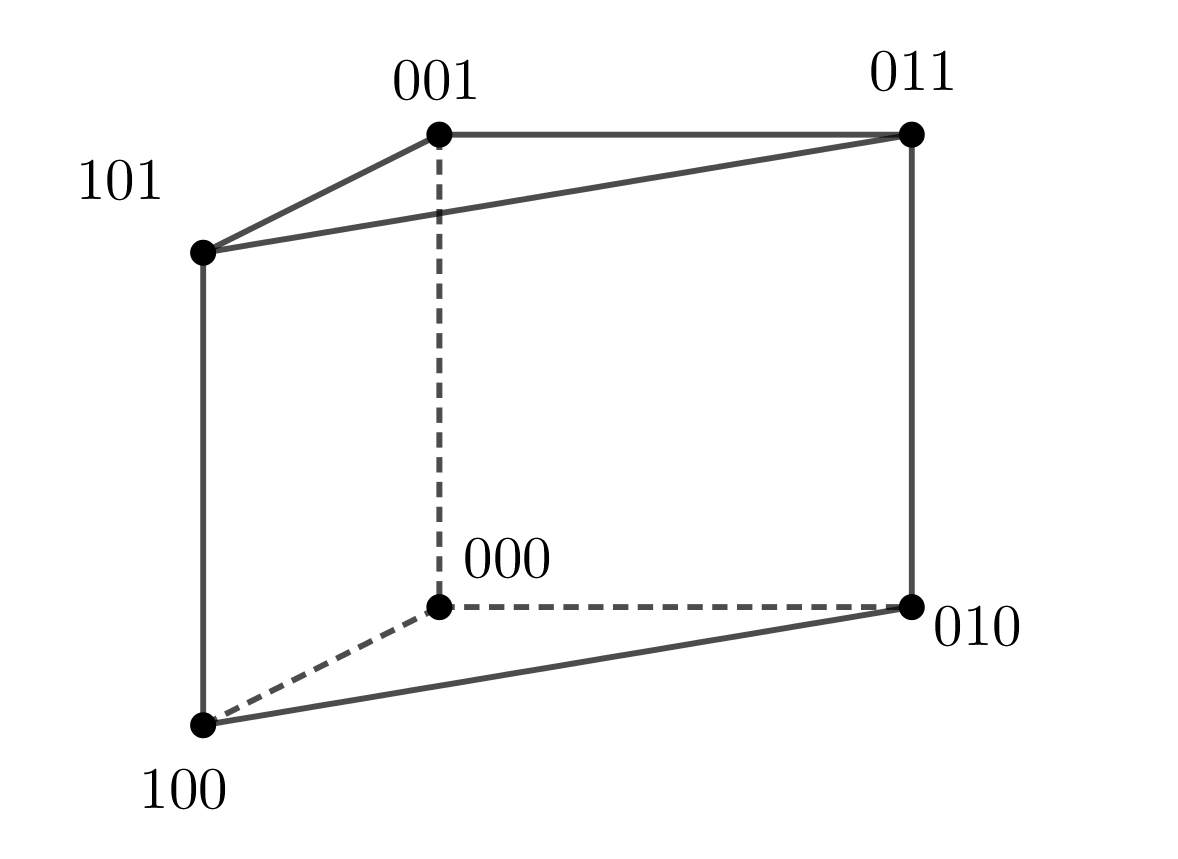}\\
$\PP_0$&$\PP_1$&$\PP_2$\\
\includegraphics[scale=.45]{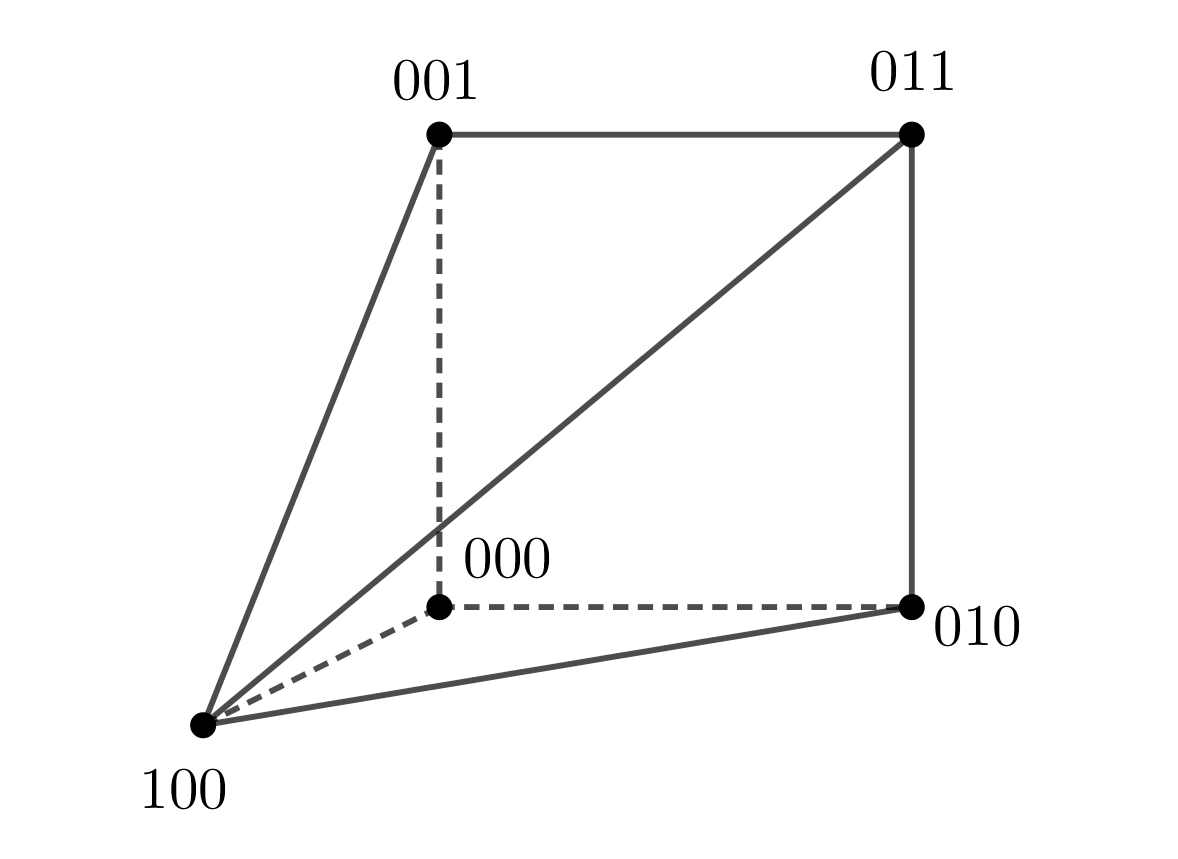}&
\includegraphics[scale=.45]{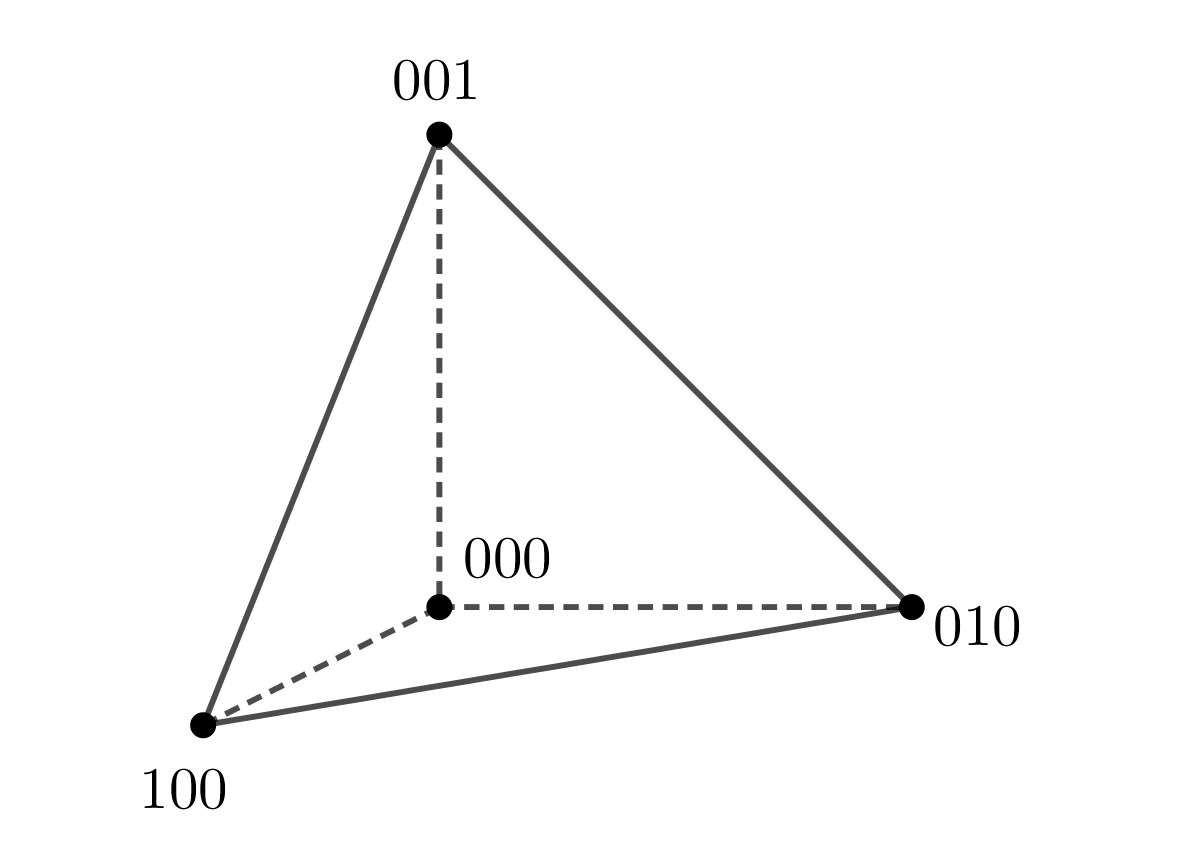}&\\
$\PP_3$&$\PP_4$&
\end{tabular}
\caption{Sequence of elements in dimension 3}
\label{3elements}
\end{figure}

For $d=3$ the resulting elements $P_i$ are displayed in Figure~\ref{3elements}.   

\begin{figure}[h!]
\centering
\begin{tabular}{ccc}
\includegraphics[scale=.3]{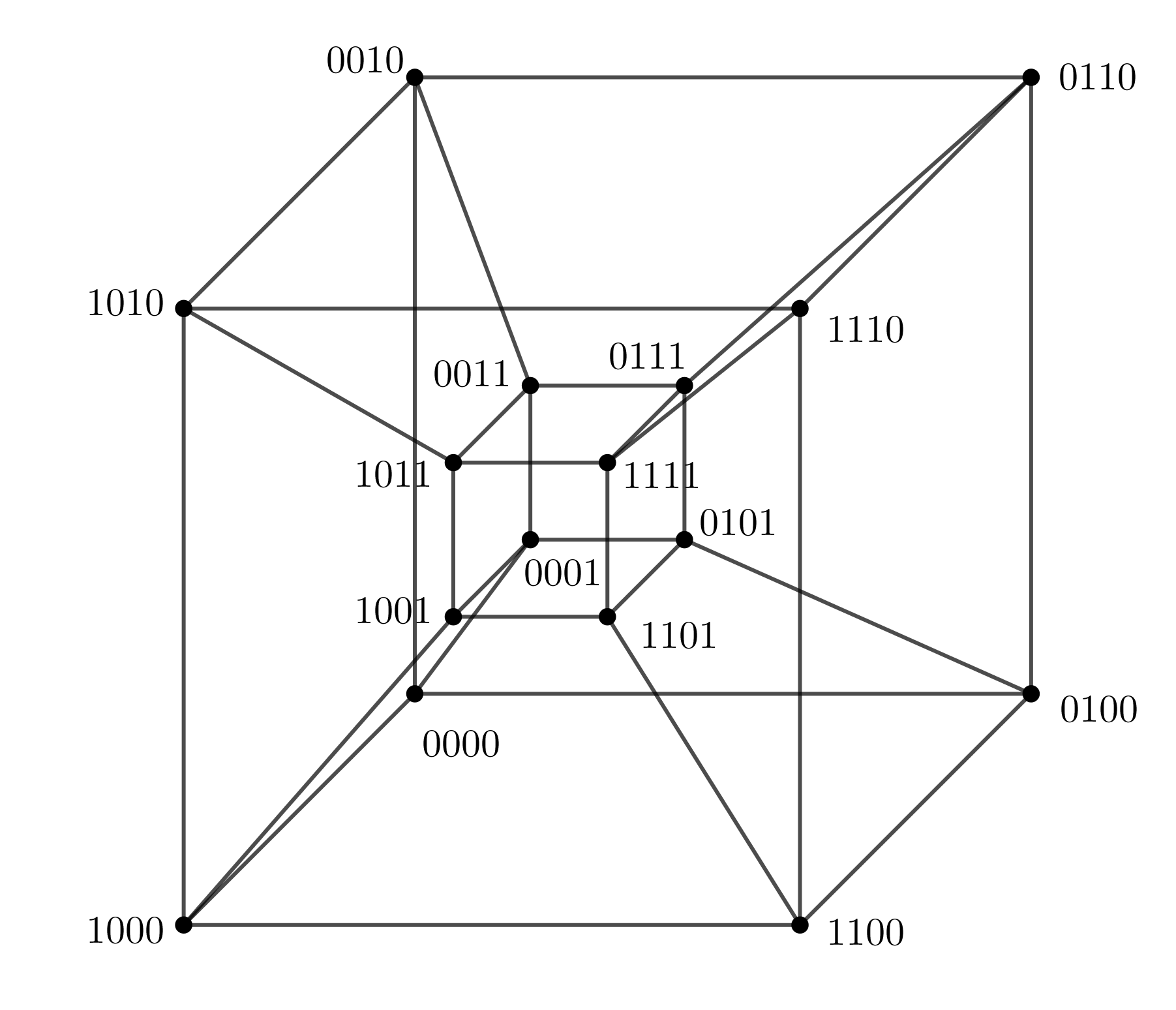}&
\includegraphics[scale=.3]{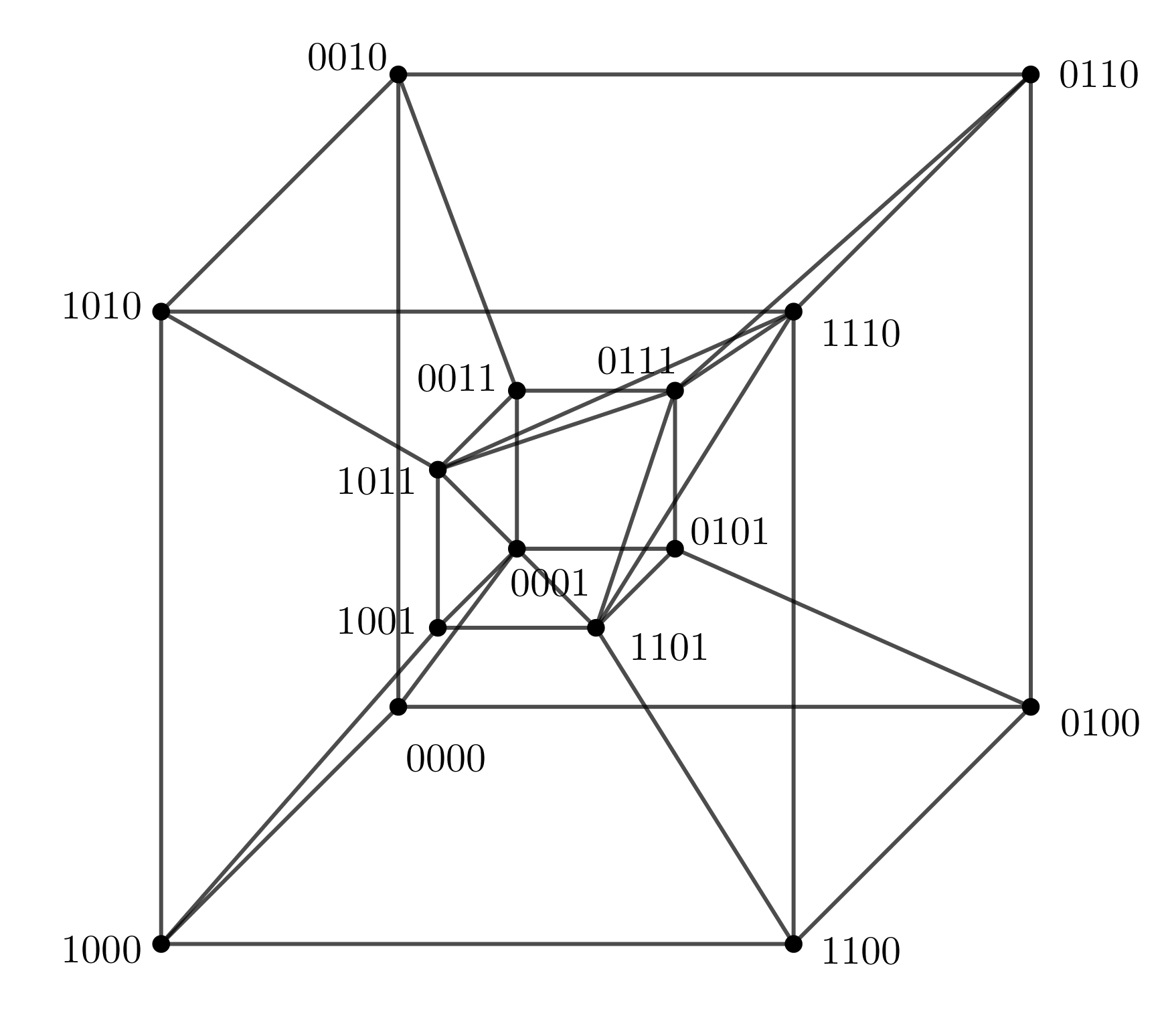}&
\includegraphics[scale=.3]{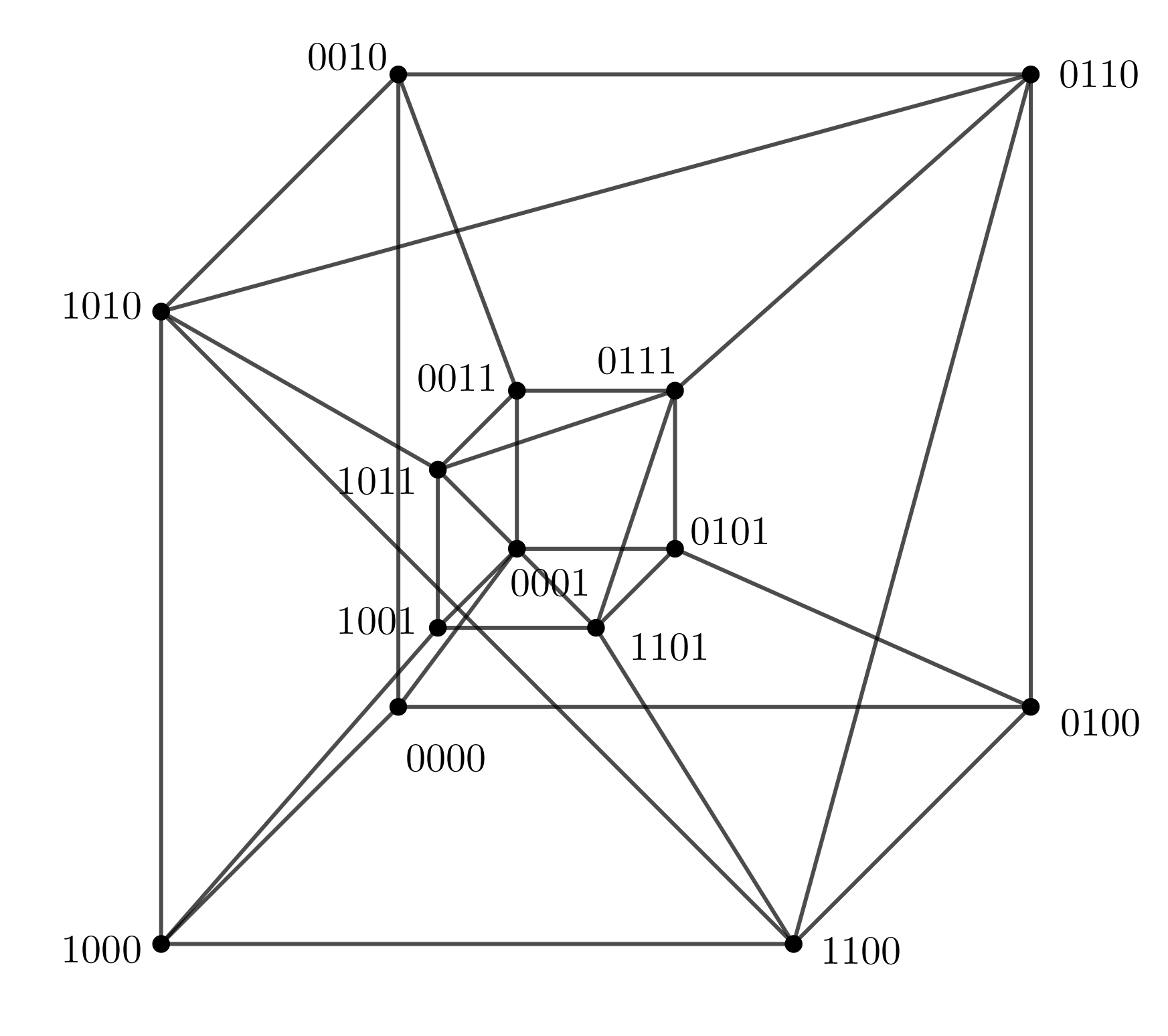}\\
$\PP'_0$&$\PP'_1$&$\PP'_2$\\
\includegraphics[scale=.3]{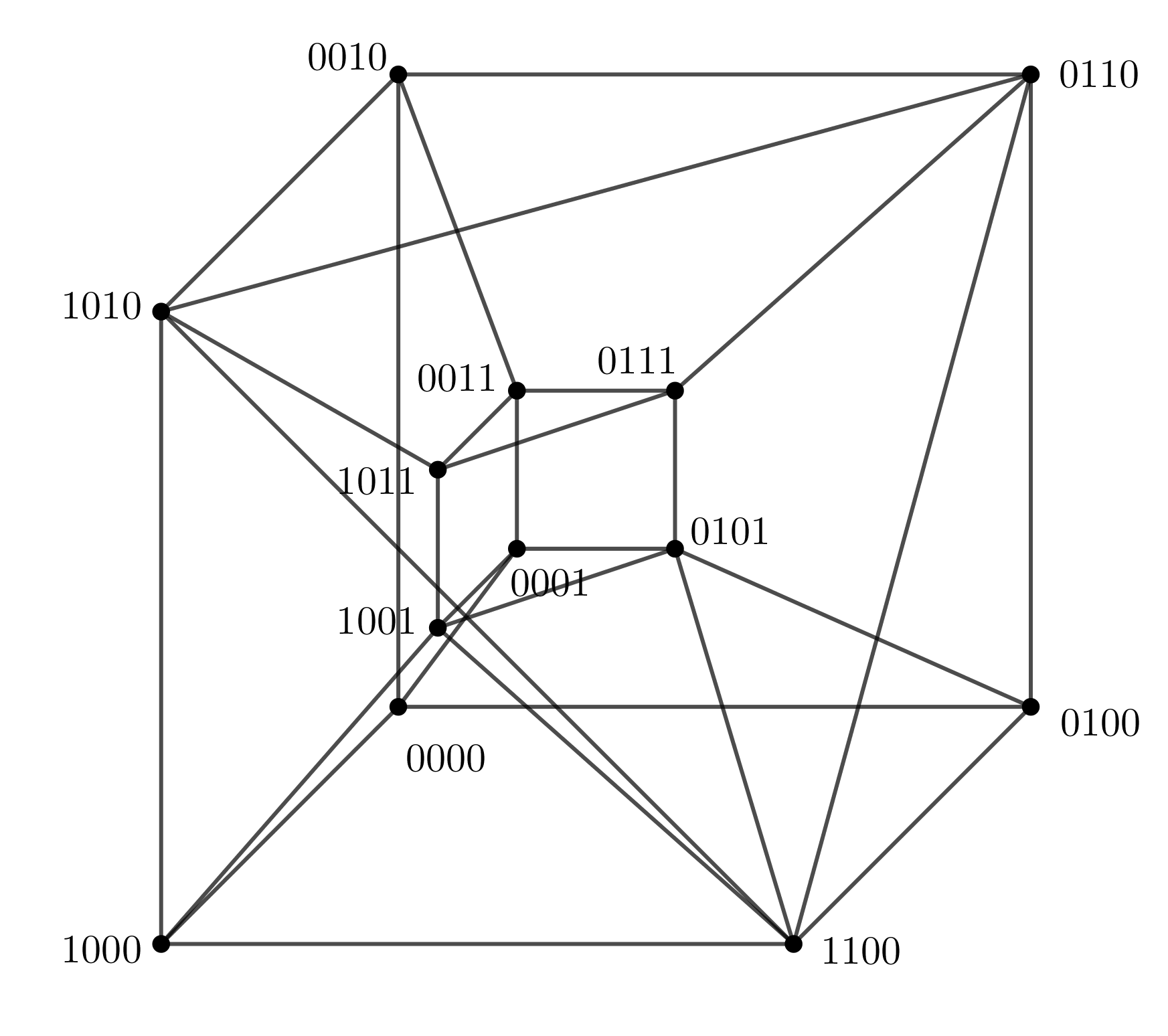}&
\includegraphics[scale=.3]{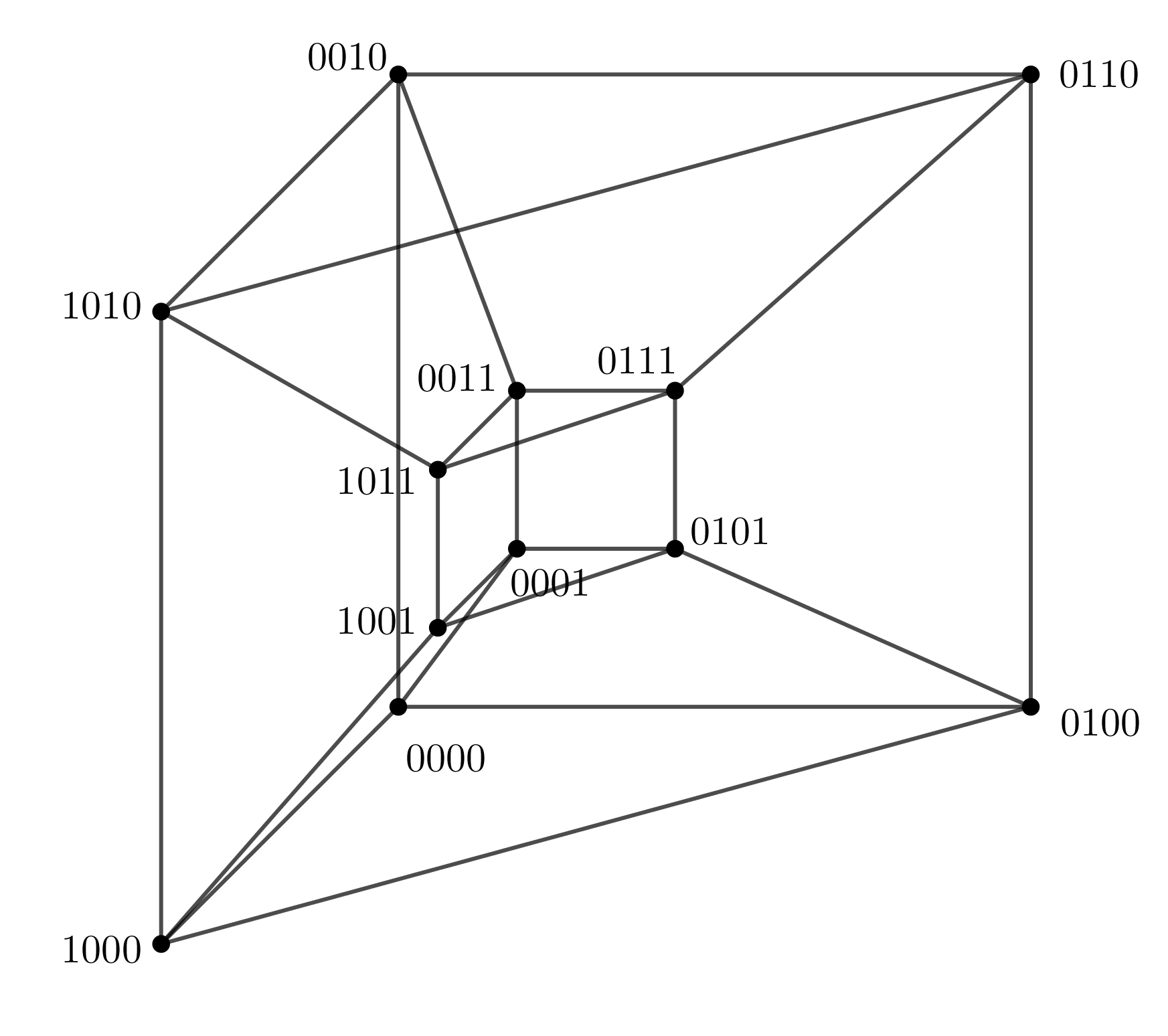}&
\includegraphics[scale=.3]{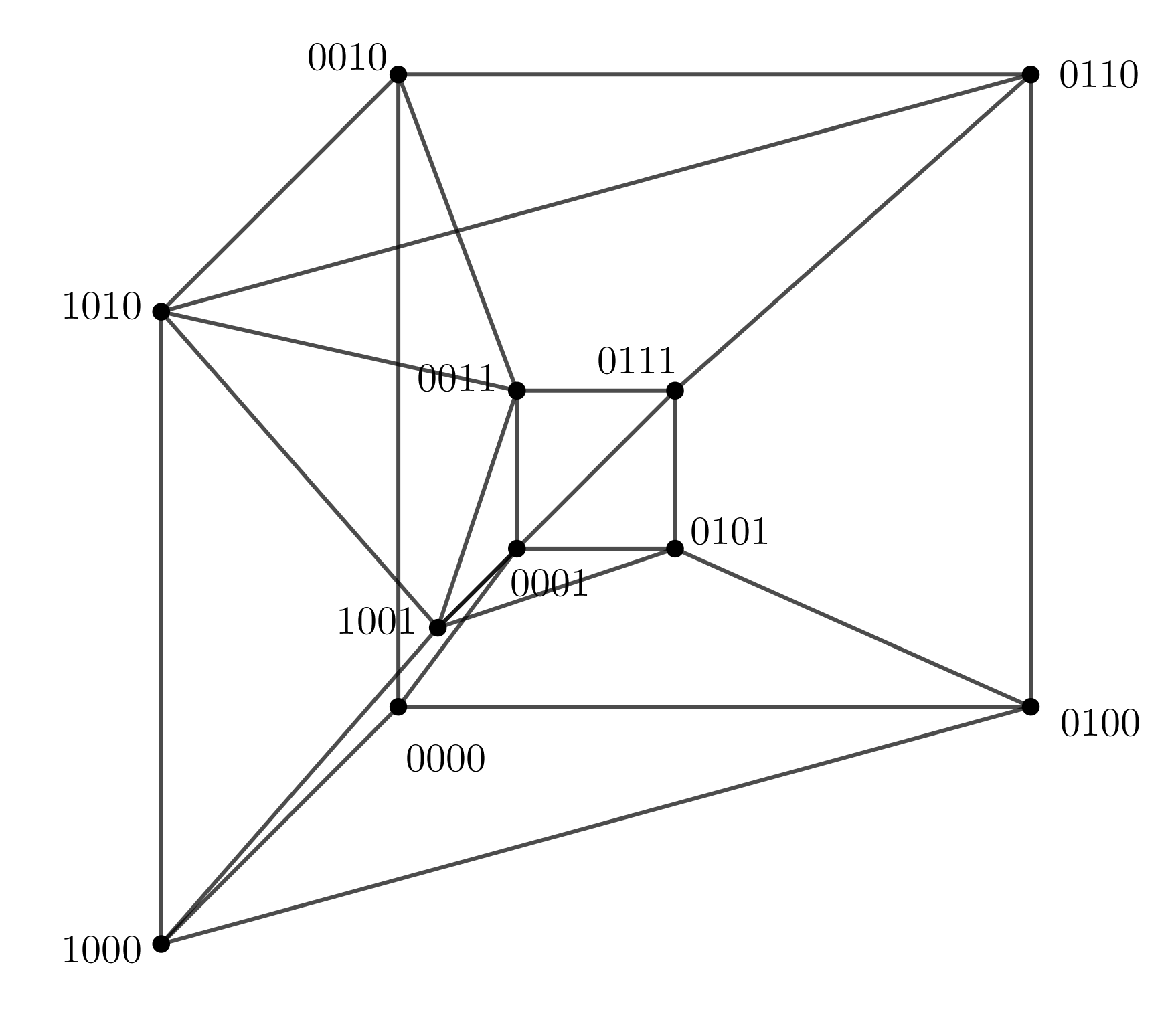}\\
$\PP'_3$&$\PP'_4$&$\PP'_5$\\
\includegraphics[scale=.3]{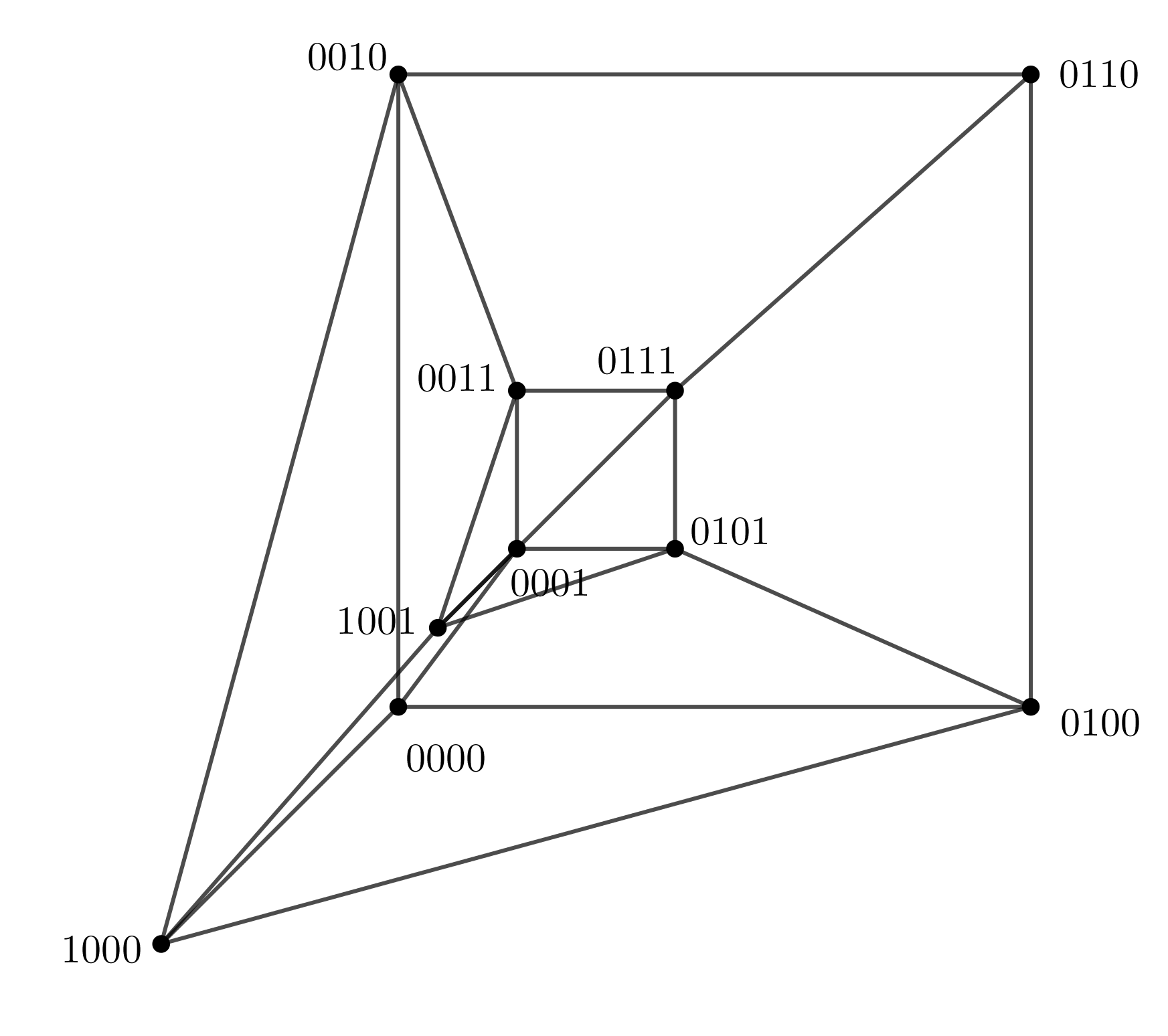}&
\includegraphics[scale=.3]{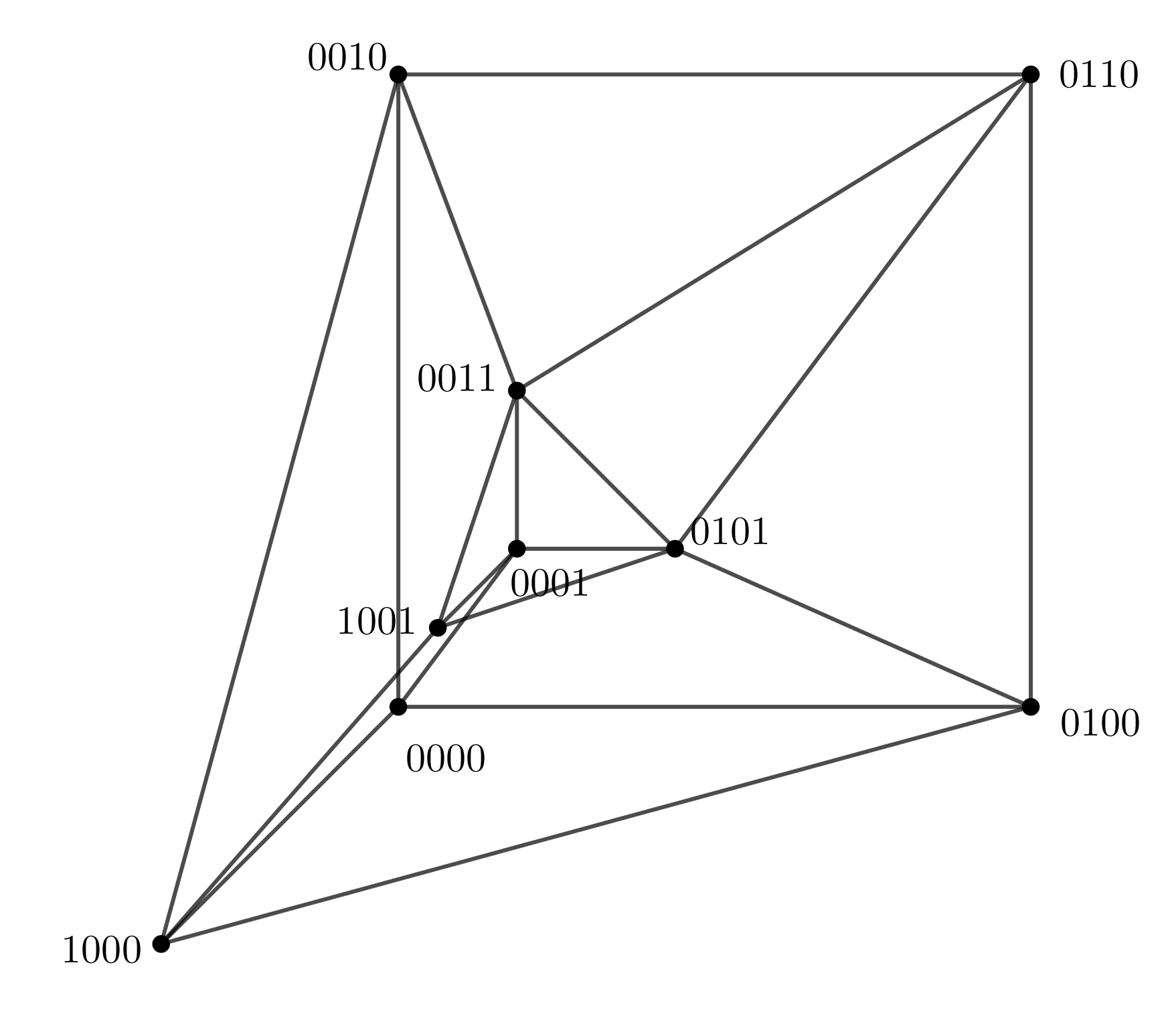}&
\includegraphics[scale=.3]{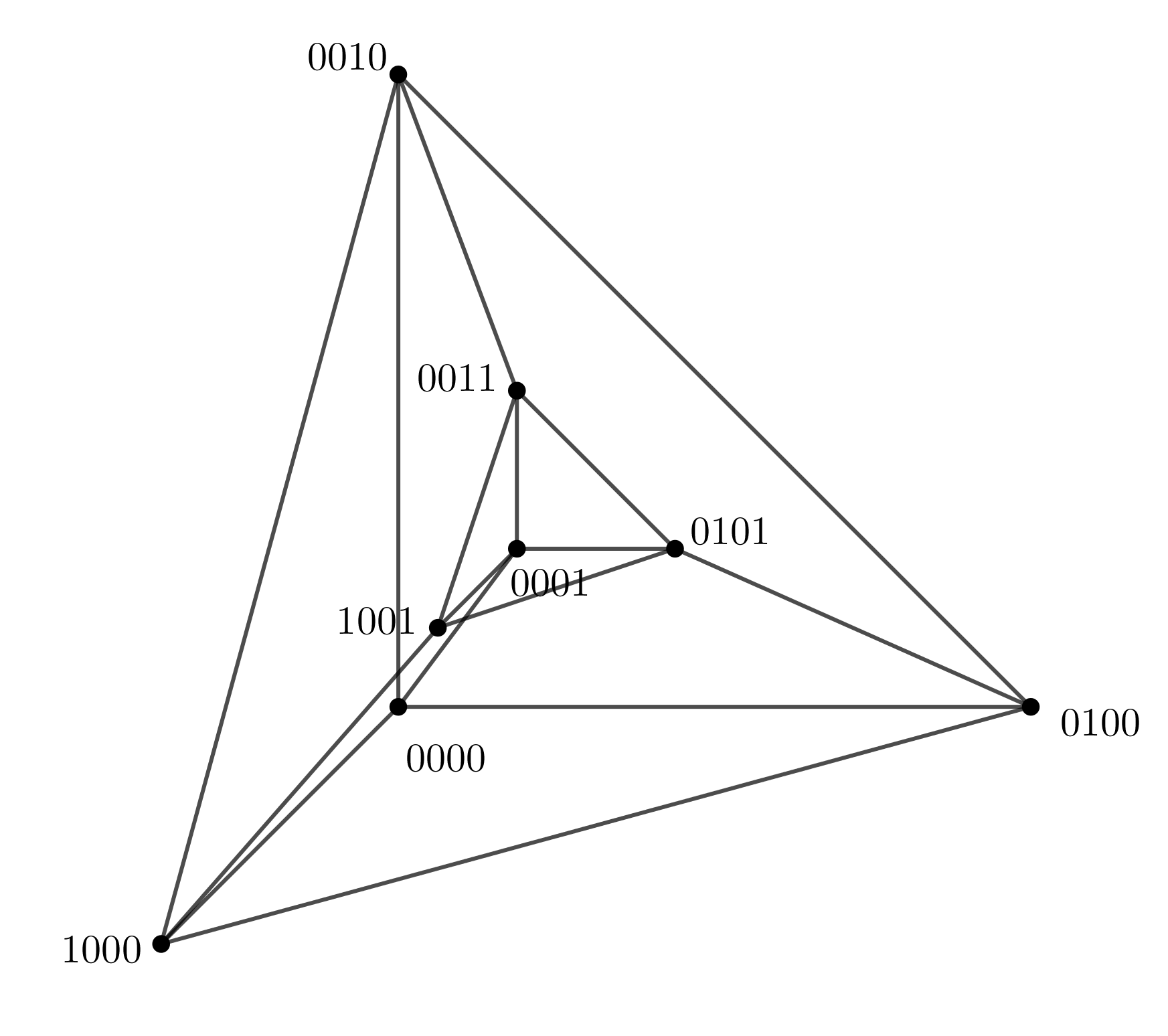}\\
$\PP'_6$&$\PP'_7$&$\PP'_8$\\
\includegraphics[scale=.3]{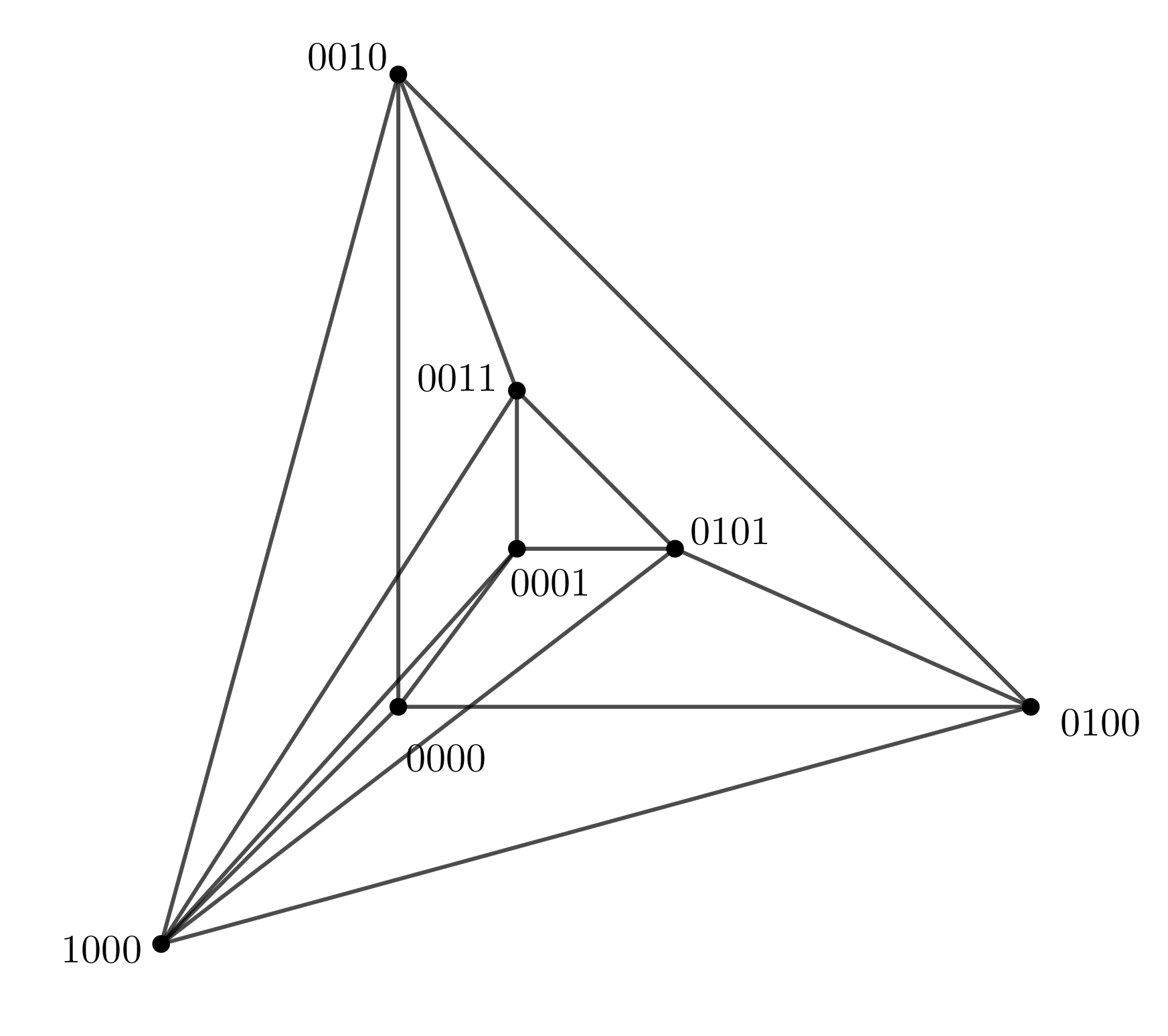}&
\includegraphics[scale=.3]{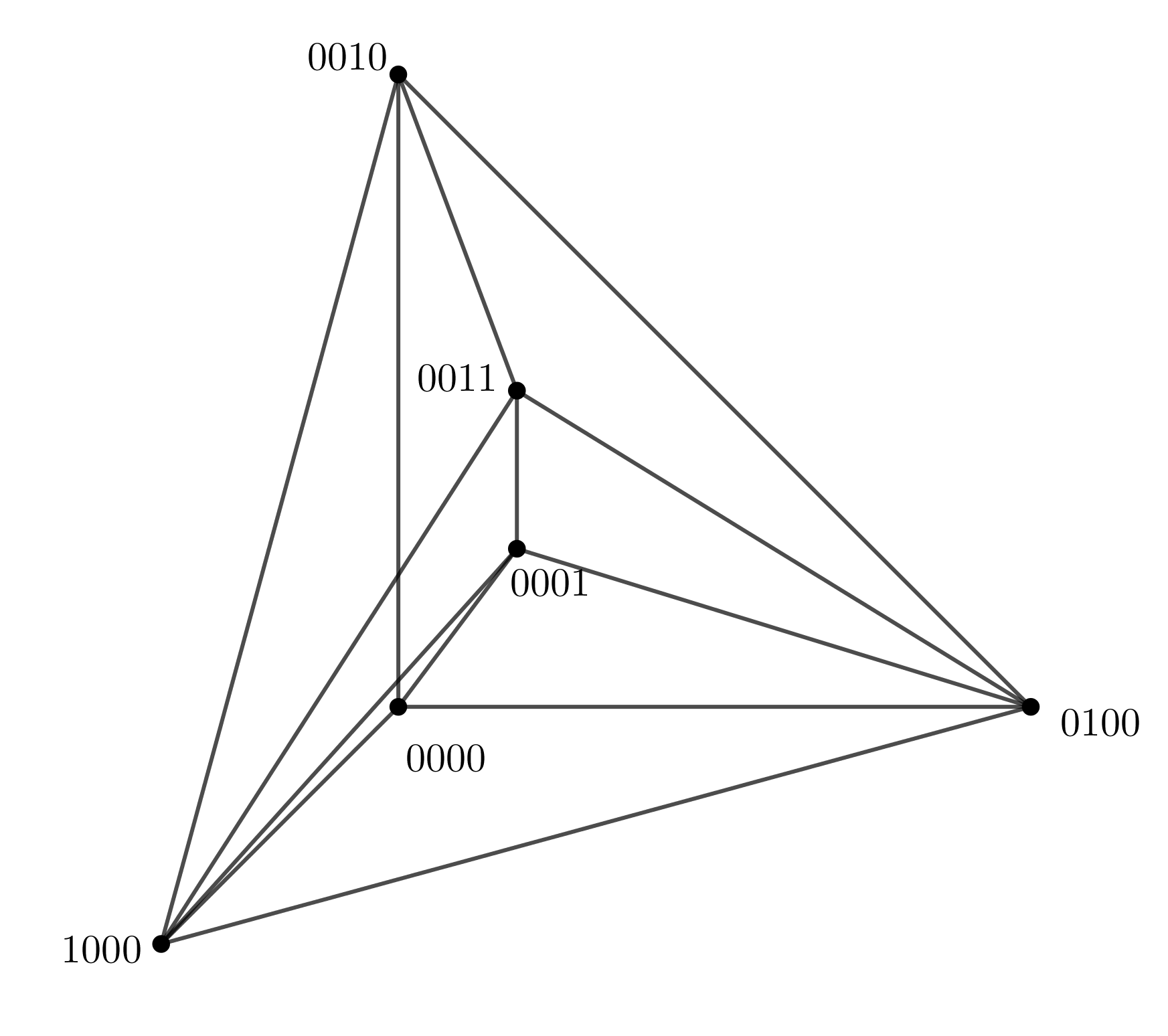}&
\includegraphics[scale=.3]{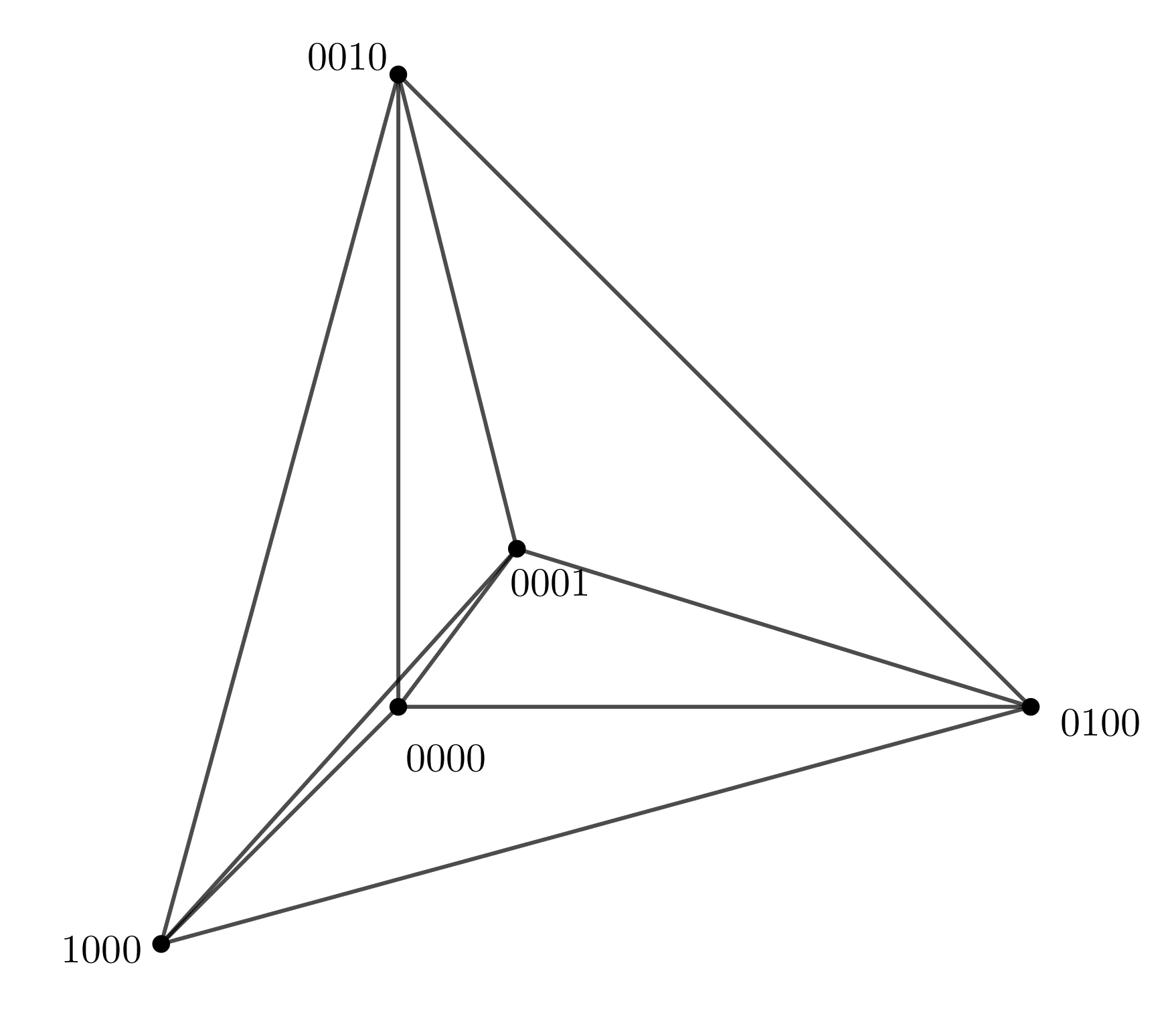}\\
$\PP'_9$&$\PP'_{10}$&$\PP'_{11}$
\end{tabular}
\caption{Sequence of elements in dimension 4}
\label{4elements}
\end{figure}

For $d=4$ the resulting elements $P'_i$ are displayed in Figure~\ref{4elements}. 
Here, the depictions are Schlegel diagrams---central projections of 4-polytopes onto a 3-dimensional hyperplane with respect to a point placed close to a facet (3-face) of the 4-polytope; see \cite{grunbaum}.
In particular,  $P'_0$ is the tesseract ($4$-cube), $P'_8$ is a tetrahedral prism (prism over the standard $3$-simplex), and $P'_{11}$ is the pentatope ($4$-simplex).
Note that $P'_0, P'_2, P'_4, P'_6$ are prisms over, respectively, the 3-dimensional elements $P_0,P_1,P_2,P_3$ of Figure~\ref{3elements}.

\begin{figure}[h!]
\centering
\begin{tabular}{ccc}
\includegraphics[scale=.3]{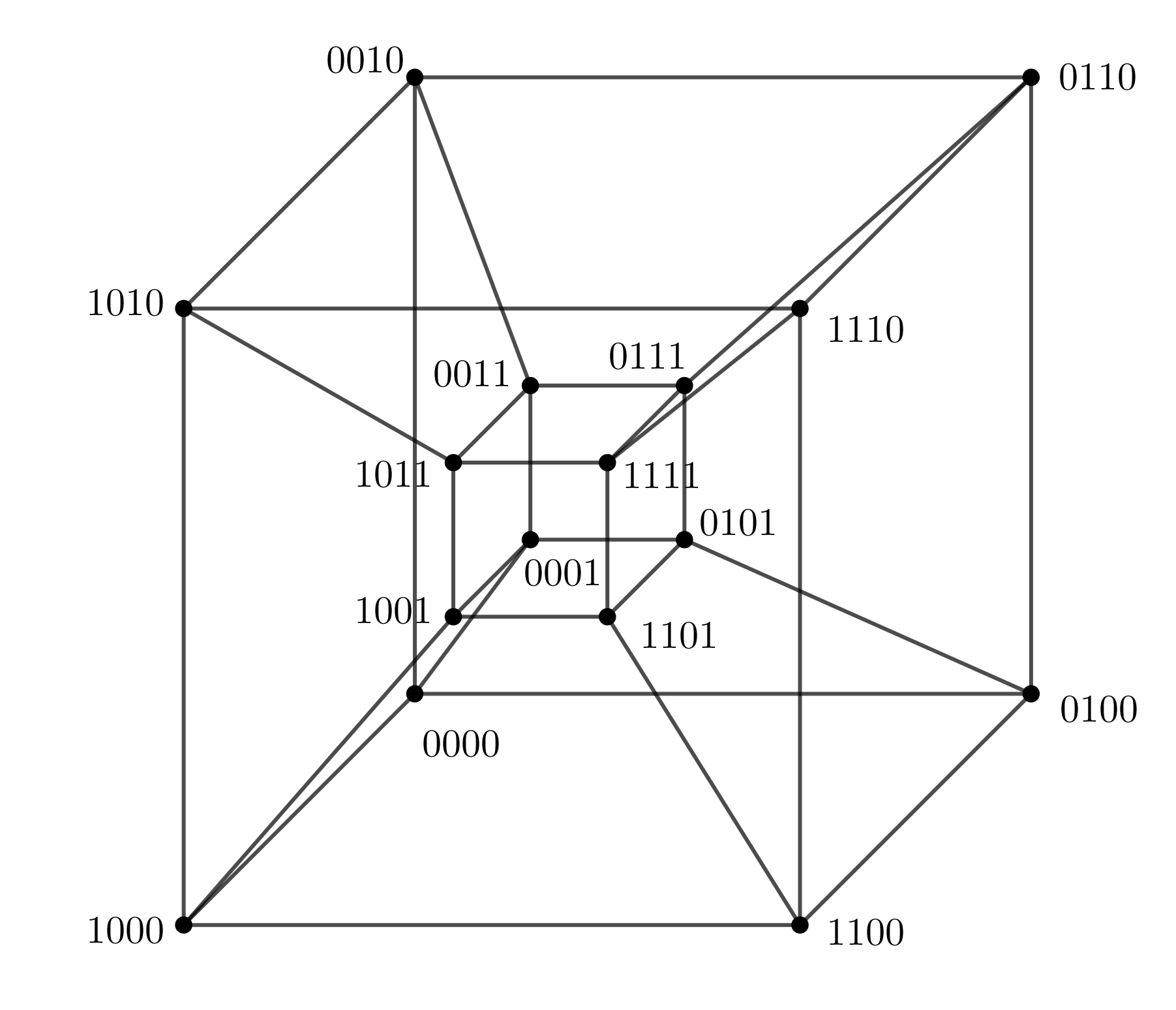}&
\includegraphics[scale=.3]{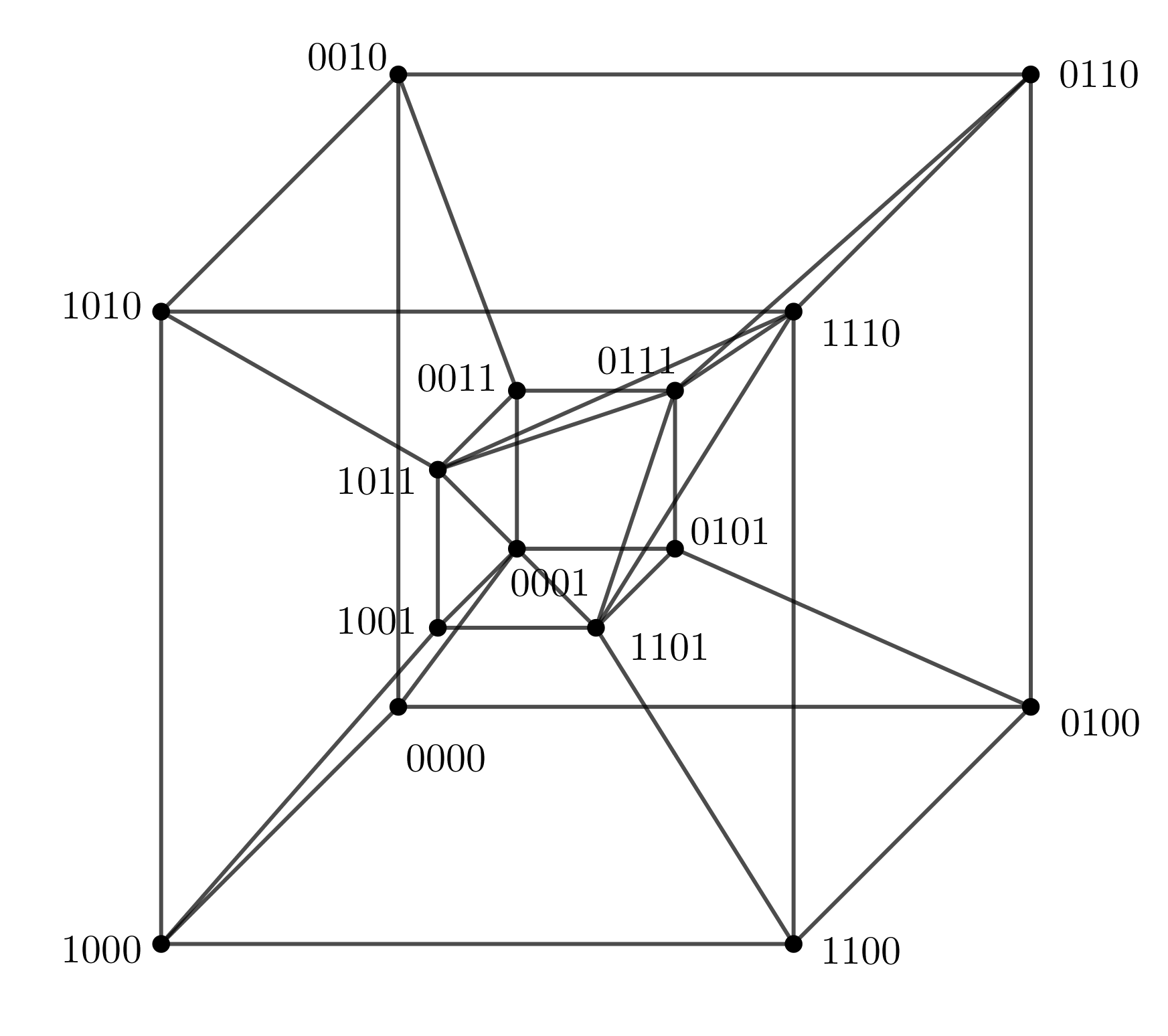}&
\includegraphics[scale=.3]{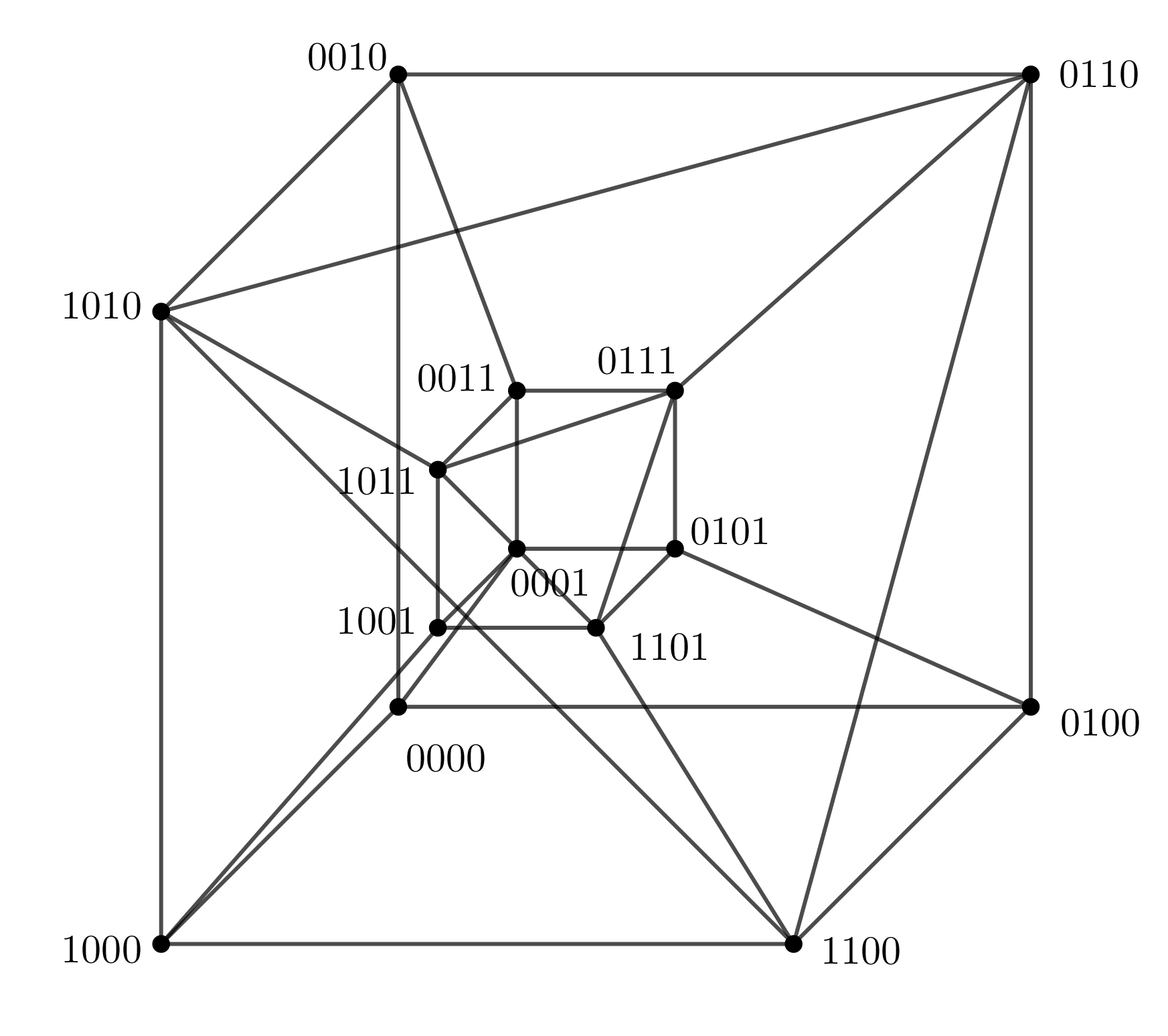}\\
$\PP'_0$&$\PP'_1$&$\PP'_2$\\
\includegraphics[scale=.3]{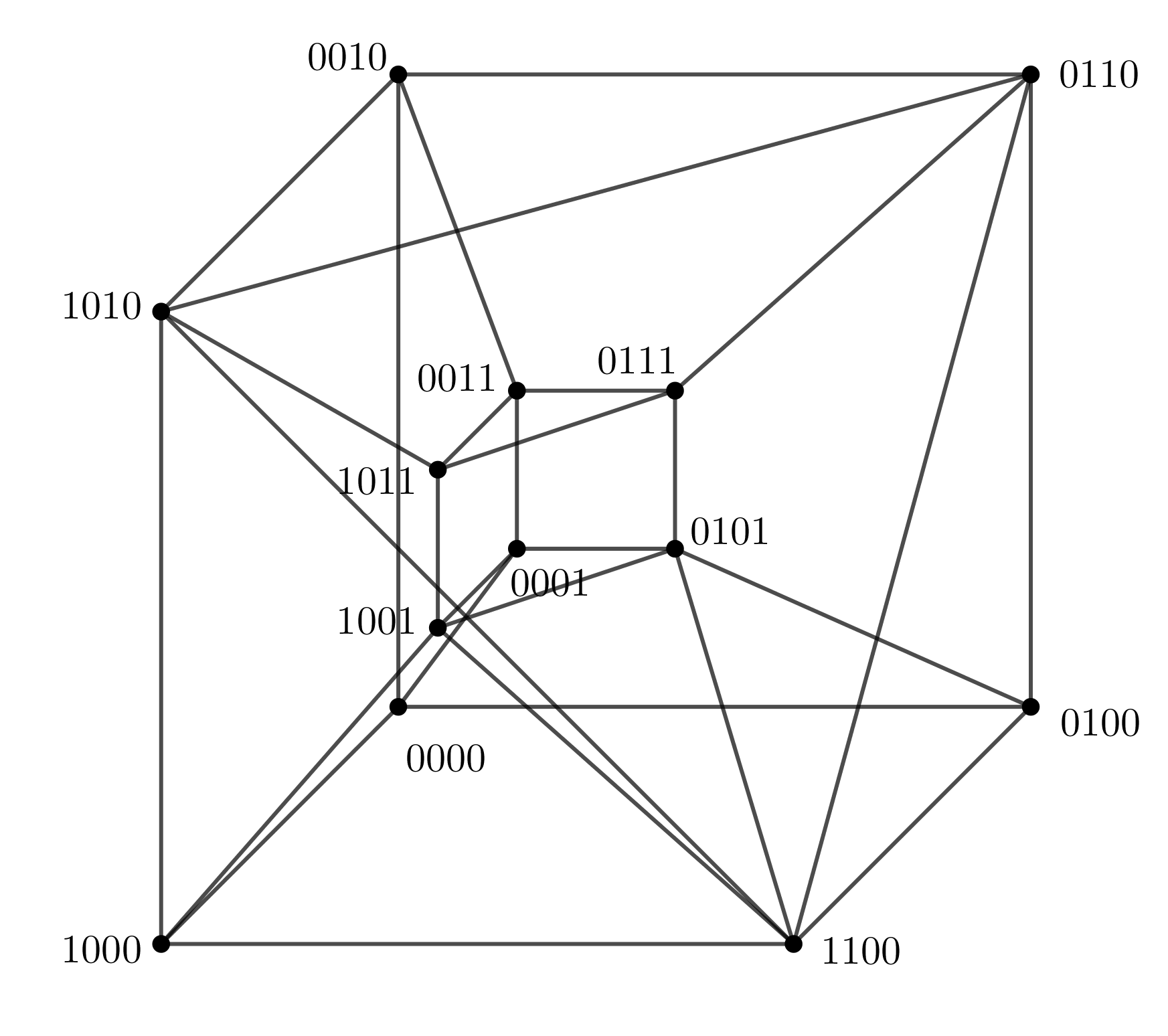}&
\includegraphics[scale=.3]{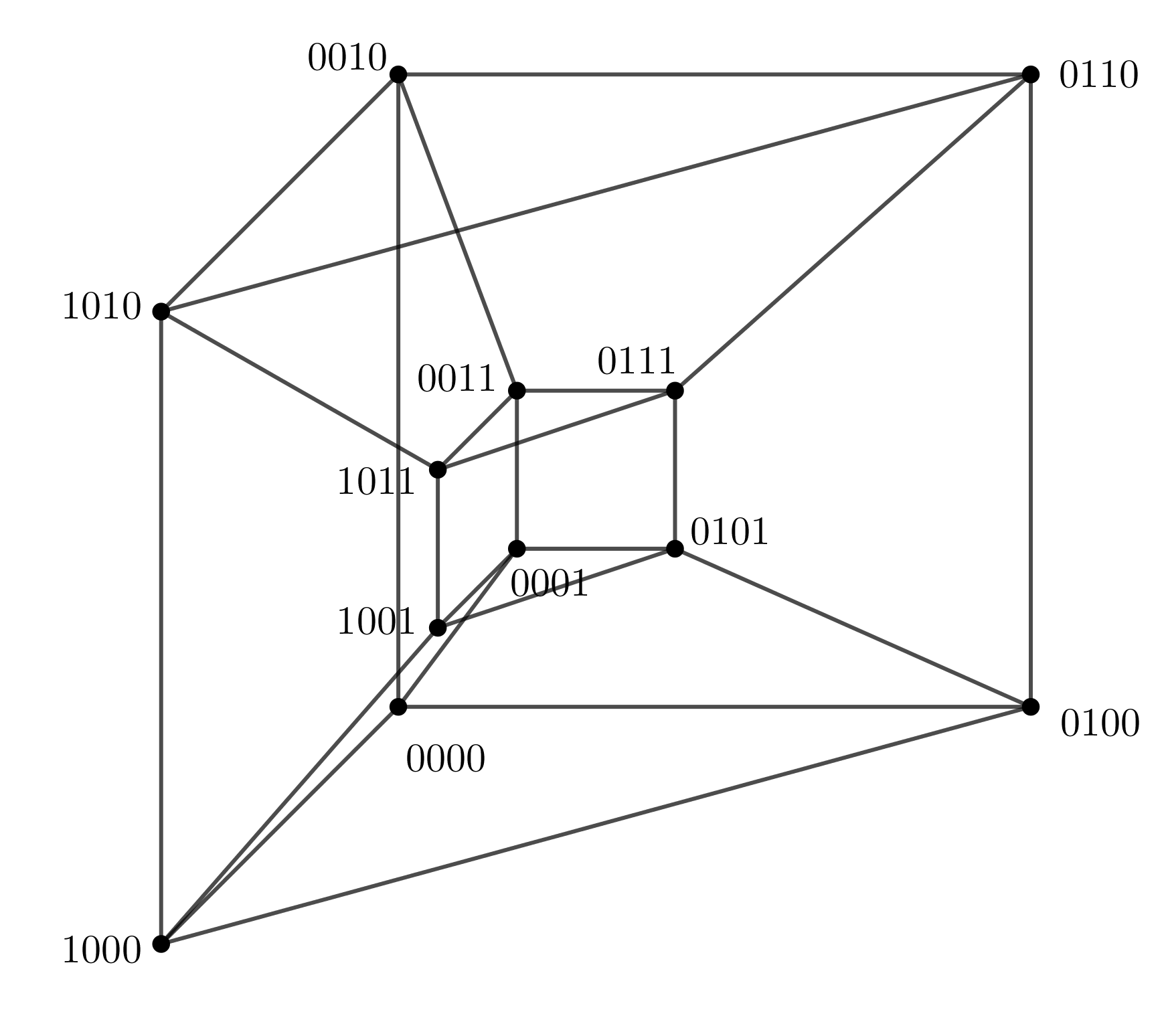}&
\includegraphics[scale=.3]{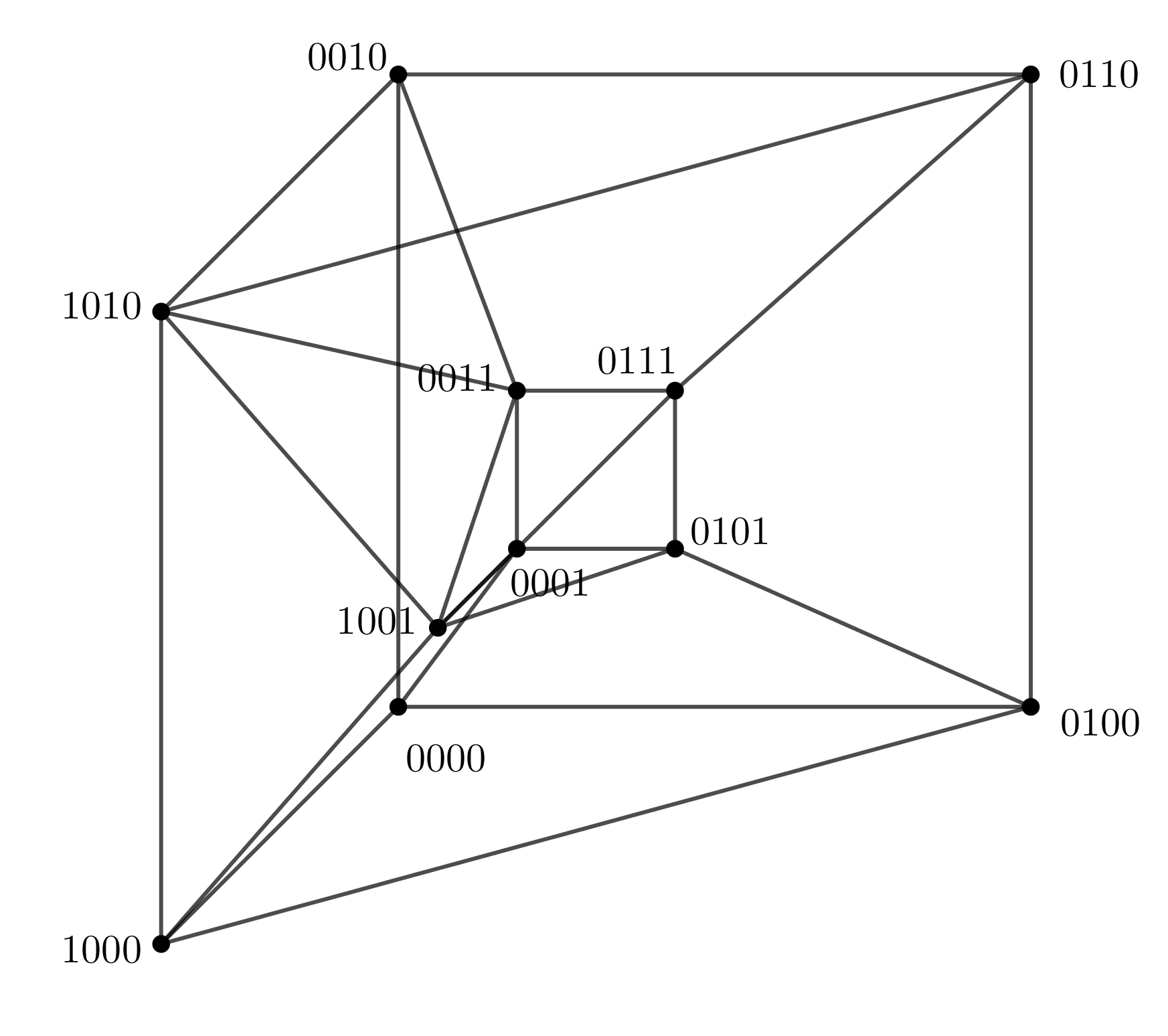}\\
$\PP'_3$&$\PP'_4$&$\PP'_5$\\
\includegraphics[scale=.3]{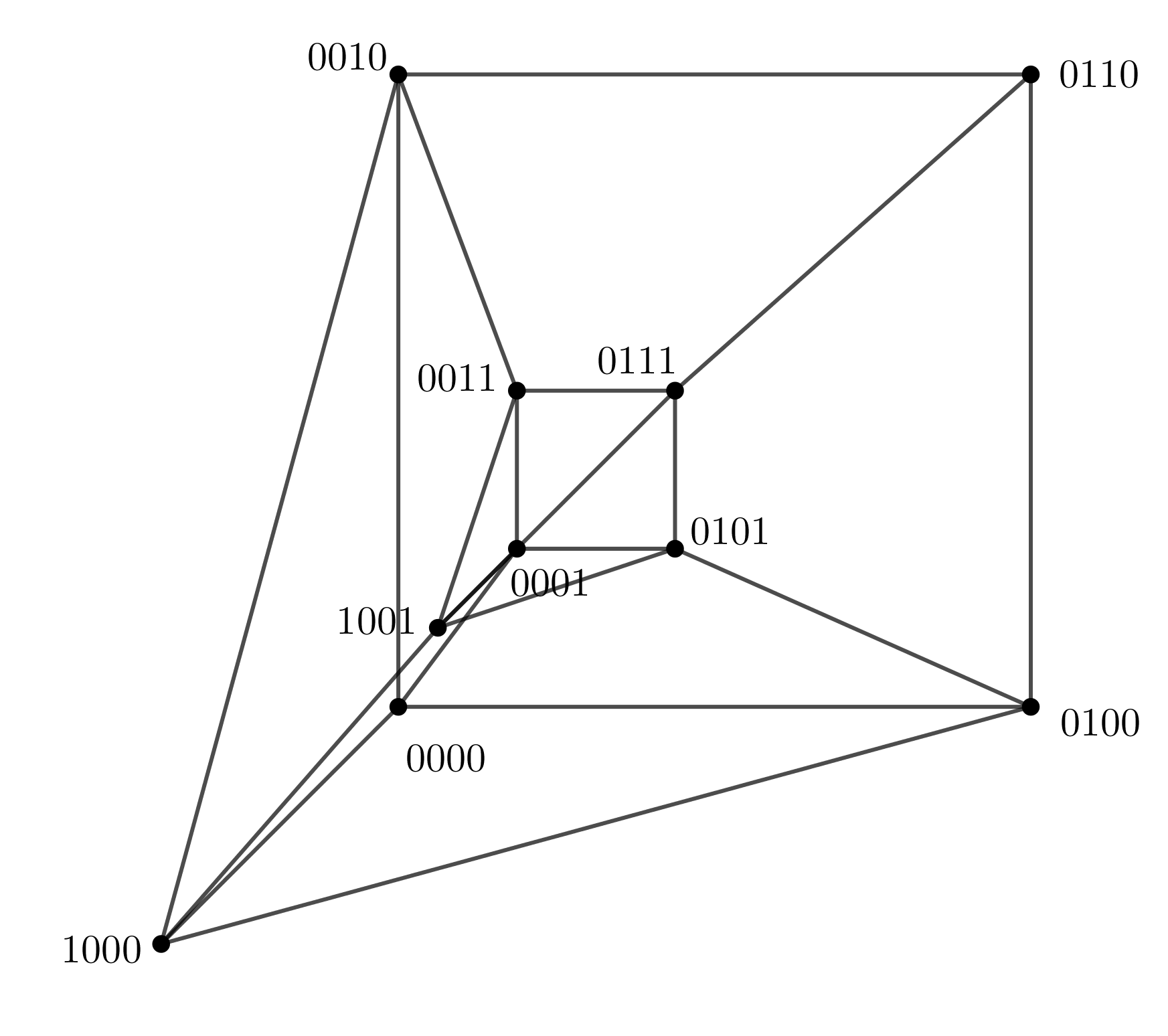}&
\includegraphics[scale=.3]{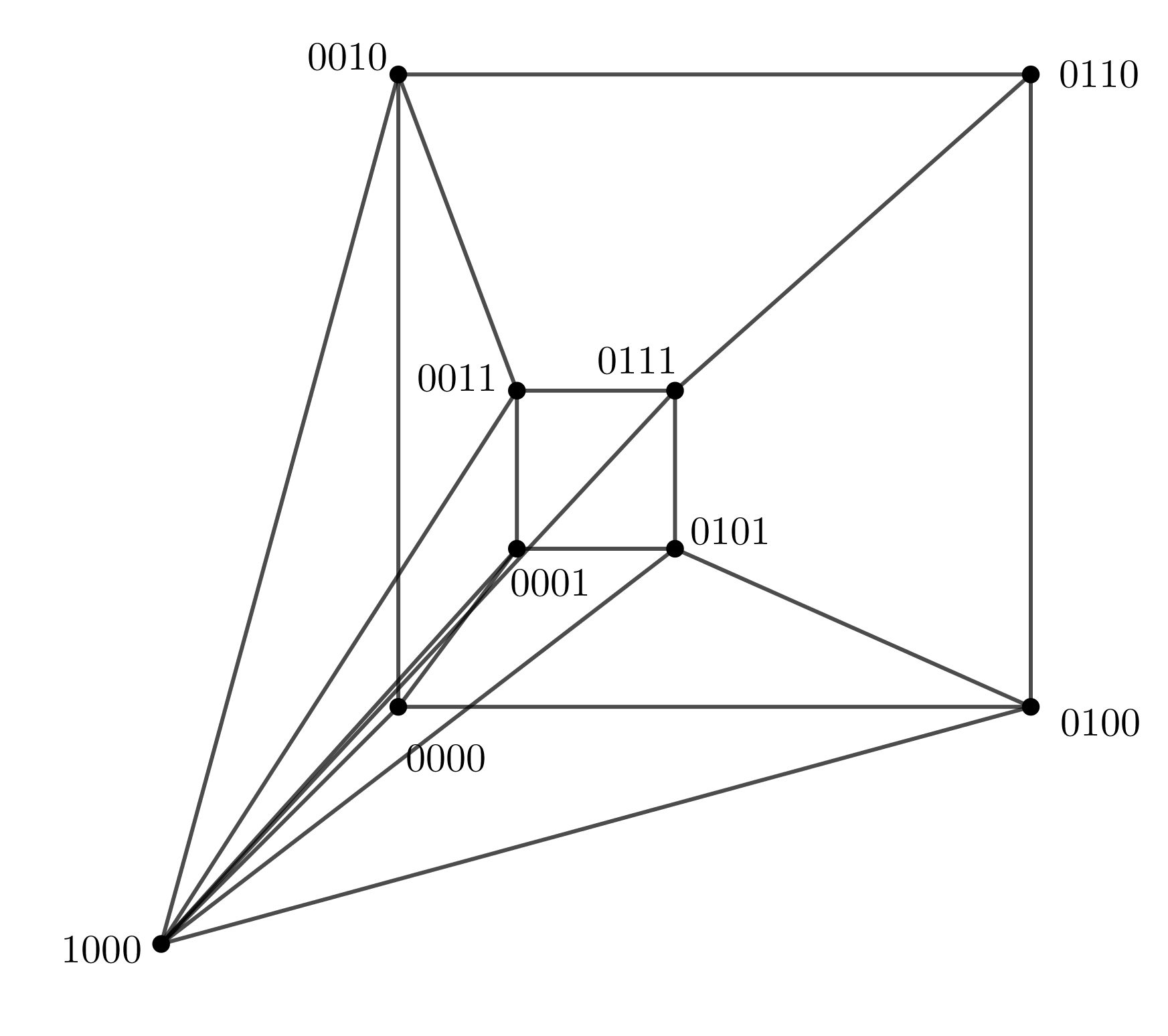}&
\includegraphics[scale=.3]{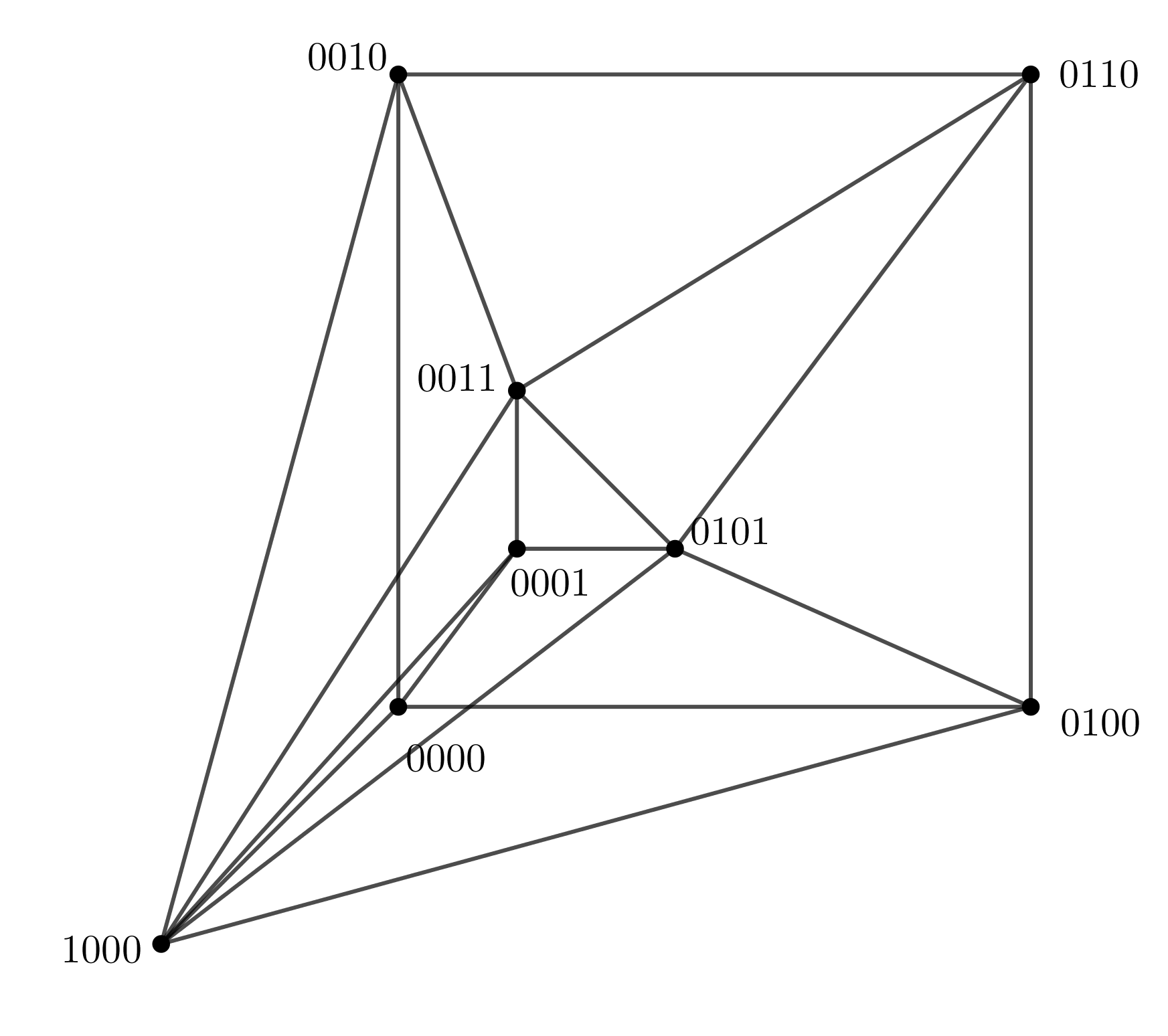}\\
$\PP'_6$&$\PP'_7$&$\PP'_8$\\
\includegraphics[scale=.3]{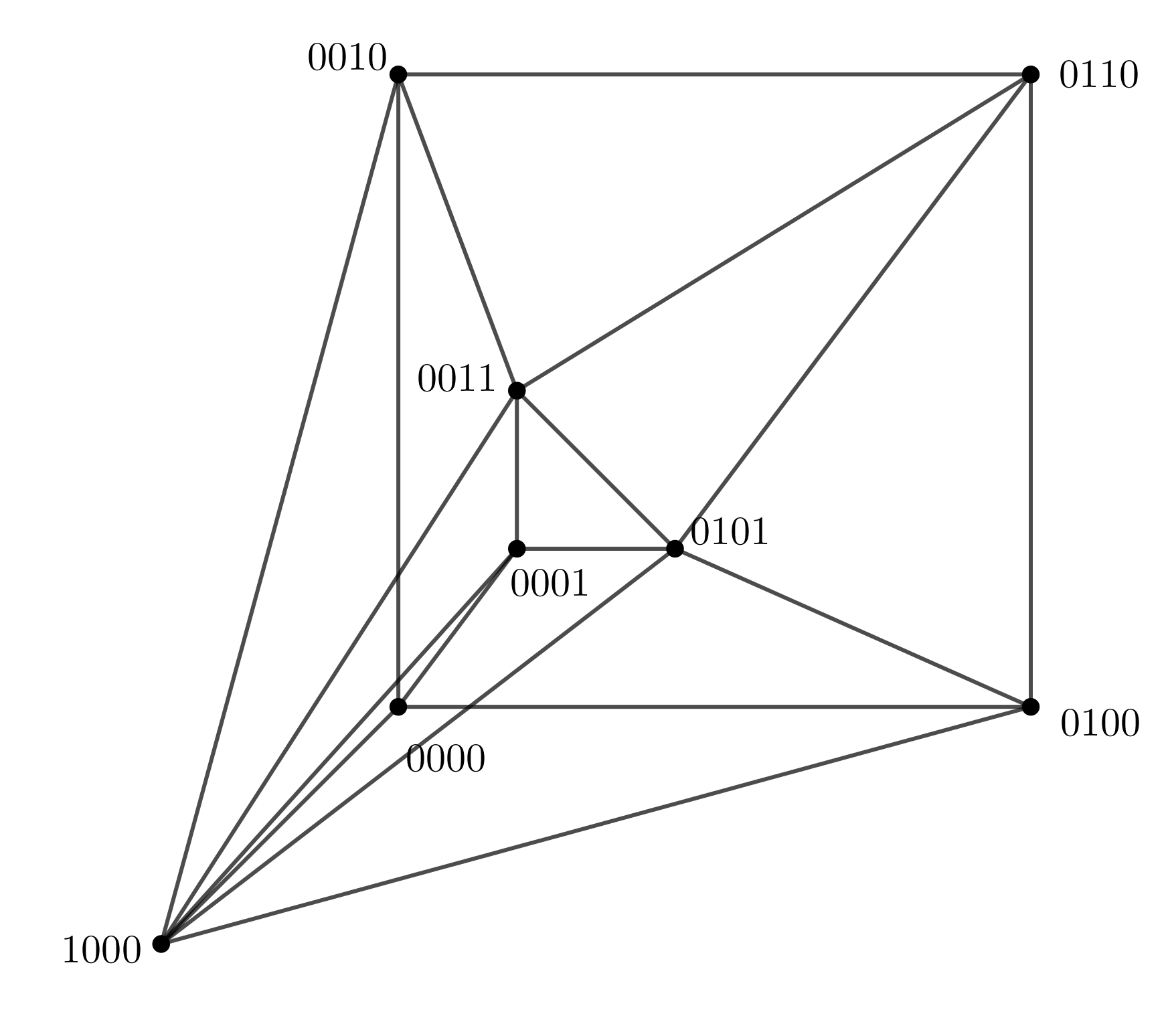}&
\includegraphics[scale=.3]{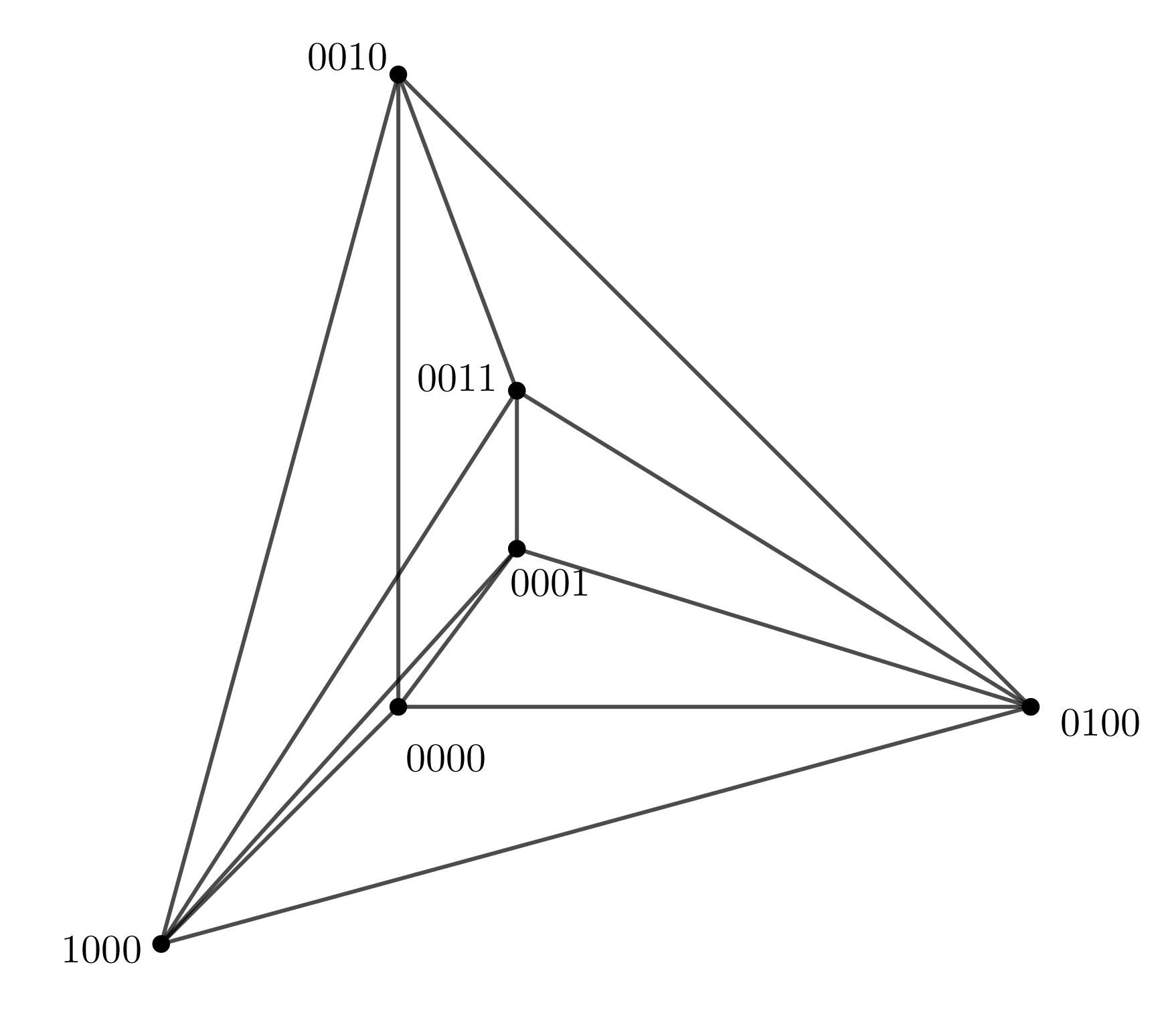}&
\includegraphics[scale=.3]{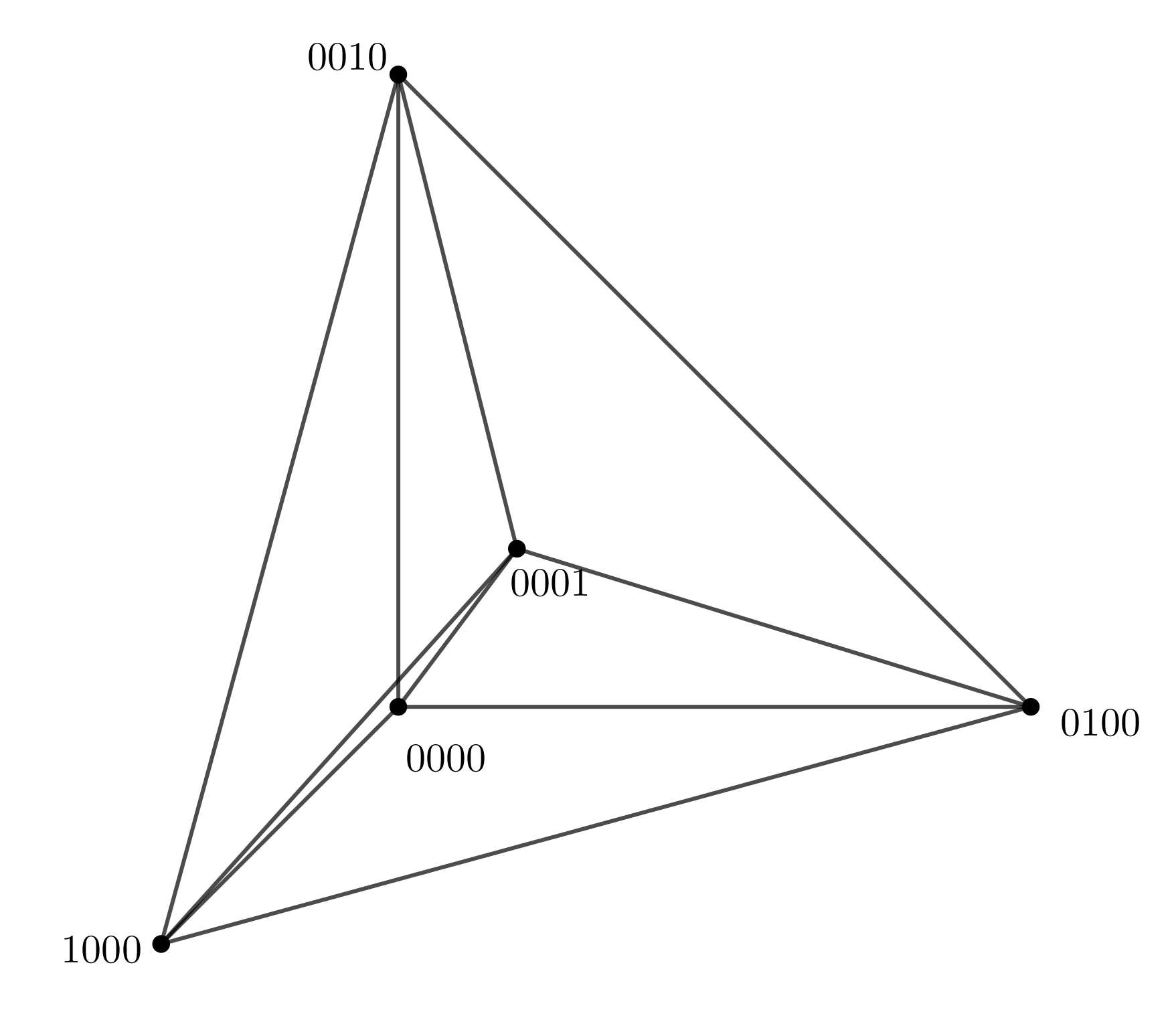}\\
$\PP'_9$&$\PP'_{10}$&$\PP'_{11}$
\end{tabular}
\caption{Second sequence of elements in dimension 4}
\label{4elementsb}
\end{figure}

The second vertex sequence of interest is the one obtained by regarding the vertices of the $d$-cube as numbers expressed in binary notation.  One may list them in decreasing order, but omit the vertices of the standard $d$-simplex.  The corresponding low dimensional sequences are shown in Table~\ref{vertexsequencebinary}.
\begin{table}[h!]
\centering
\[
\begin{array}{ccc}
d=2&d=3&d=4\\
\hline
11&111&1111\\
&110&1110\\
&101&1101\\
&011&1100\\
&&1011\\
&&1010\\
&&1001\\
&&0111\\
&&0110\\
&&0101\\
&&0011
\end{array}
\]
\caption{Decreasing binary sequences of vertices in low dimensions}
\label{vertexsequencebinary}
\end{table}
When $d=3$ this sequence matches the earlier one, but they are different when $d\geq4$.  In particular, this sequence does not contain the tetrahedral prism.
For $d=4$ the resulting elements $P'_i$ are displayed in Figure~\ref{4elementsb}.

\clearpage

These two types of sequences are far from exhaustive, since the number of combinatorially inequivalent 0/1 $d$-polytopes increases doubly-exponentially with 
$d$---see~\cite{ziegler}. For example, for $d =4$ there are 60,879 such 0/1 $d$-polytopes~\cite{zong2005known}. Nevertheless, we have highlighted the two sequences above because for $d =2$ and $d = 3$ they contain elements which are frequently used in conventional finite element methods (as we mentioned previously). Furthermore, for $d = 4$ they contain new element types. Finally, the first sequence also contains all of the standard 4D element types such as the tesseract, the tetrahedral prism, and the pentatope. In the next few sections, we will focus on developing numerical interpolation and integration procedures for these standard 4D elements.

\section{Polynomial Bases and Symmetry Groups} \label{sec;symmetry}

\subsection{Overview}

In the previous section, we constructed sequences of elements based on degenerations of the $4$-cube $\left[0,1\right]^4$, and we recovered the standard tesseract, tetrahedral prism, and pentatope elements. Our choice of the interval $\left[0,1\right]$ was consistent with standard mathematical conventions, however, for engineering purposes it is common practice to use the interval $\left[-1,1\right]$ instead. For certain operations (such as interpolation and integration) this interval is more convenient, as it is symmetrical and centered at the origin. Of course, formally speaking, the mathematical and engineering intervals are topologically equivalent, as they are related by a linear map. Therefore, the intervals are interchangeable, and we have a free choice of which one to use.

In what follows, we switch to the engineering convention, and construct a new set of reference elements for the tesseract, tetrahedral prism, and pentatope. In addition, we discuss the essential properties of these elements which are useful for implementation purposes. In particular,  integration bounds are introduced in order to define the domains of our reference elements, basis functions on each reference element are constructed for interpolation and quadrature purposes, and lastly we define symmetry groups for each reference element in order to facilitate the generation of fully symmetric quadrature rules. 

\subsection{Preliminaries}

Throughout this section, whenever we develop basis functions, we construct an orthonormal basis, i.e.,
\begin{equation*}
    \int_{\Omega} \psi_{ijkq} \left( \boldsymbol{x} \right) \psi_{rstv} \left( \boldsymbol{x} \right) d \boldsymbol{x} = \delta_{ir} \delta_{js} \delta_{kt} \delta_{qv},
\end{equation*}
where $\delta_{ir}$ is the Kronecker delta. The simplest orthonormal basis functions of degree $\porder$ have the form
\begin{equation*}
    \psi_{ijkq} \left( \boldsymbol{x} \right) = \zeta_{ijkq} \, \hat{P}_i^{(\alpha_1 ,\beta_1)} \left( a \right) \hat{P}_j^{(\alpha_2 ,\beta_2)} \left( b \right) \hat{P}_k^{(\alpha_3 ,\beta_3)} \left( c \right) \hat{P}_q^{(\alpha_4 ,\beta_4)} \left( d \right),
\end{equation*}
where $0 \leq i + j +k +q \leq \porder$, $a = a\left(x_1,x_4\right)$, $b = b\left(x_2,x_4\right)$, $c = c\left(x_3,x_4\right)$, and $d = d\left(x_4\right)$ are functions depending on the element type, $\zeta_{ijkq}$, $\alpha_1 ,\ldots, \alpha_4$ and $\beta_1, \ldots, \beta_4$ are constants, and $\hat{P}_n^{(\alpha ,\beta)}$ is the 1D orthonormal Jacobi polynomial defined as
\begin{equation*}
    \hat{P}_n^{(\alpha,\beta)} \left( x \right) = \frac{P_n^{(\alpha,\beta)} (x)}{\sqrt{\frac{2^{\alpha+\beta+1}}{2n + \alpha + \beta + 1}\frac{\left( n + \alpha \right)! \left( n + \beta \right)! }{ n! \left( n + \alpha + \beta \right)!}}}.
\end{equation*}
Here, $P_n^{(\alpha,\beta)}$ is the standard orthogonal (but not orthonormal) Jacobi polynomial.

\subsection{Tesseract}

\begin{figure}[h!]
\centering
\includegraphics[width=0.7\textwidth]{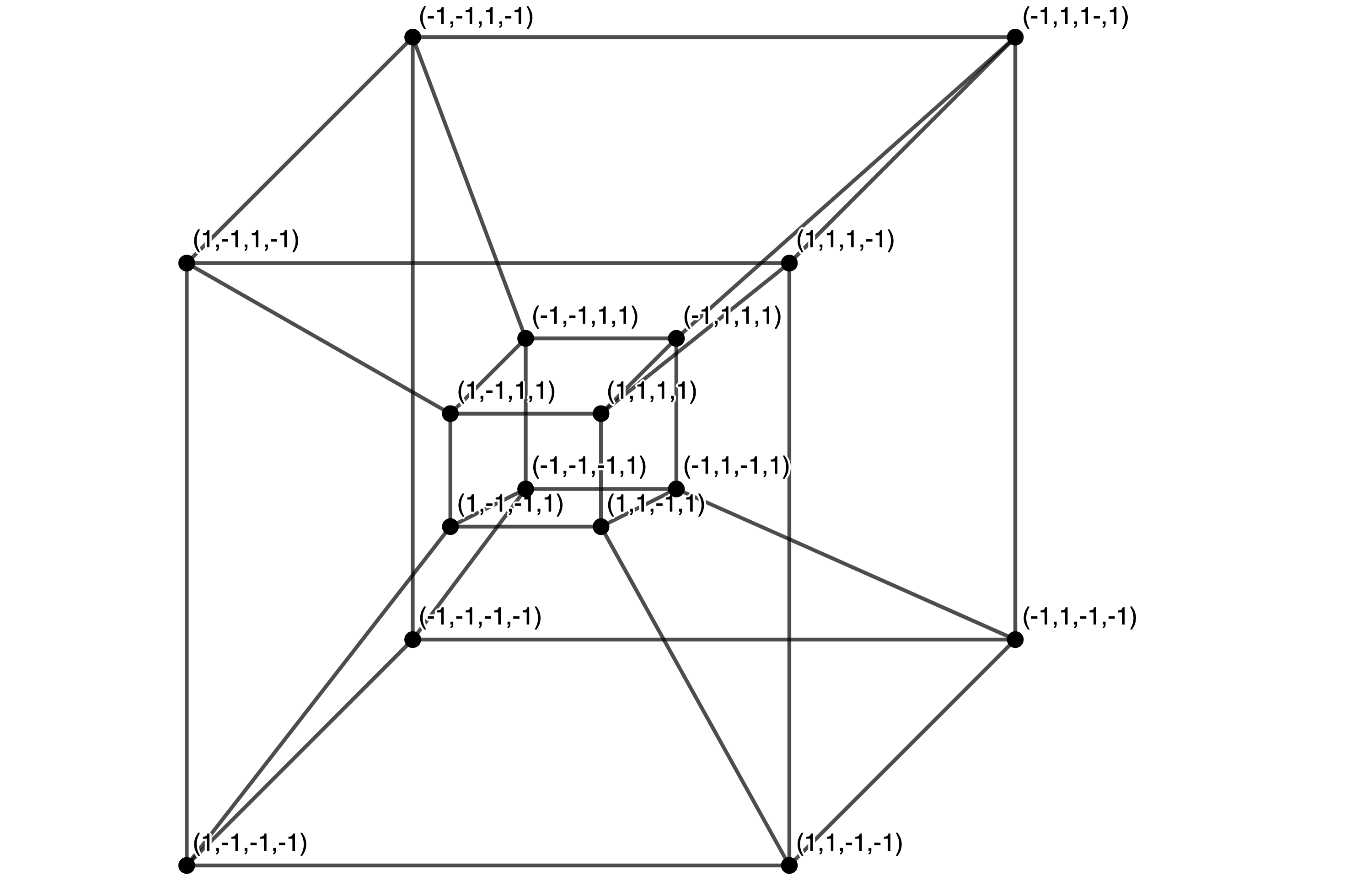}
\caption{Tesseract reference element.}
\label{eng_tesseract}
\end{figure}

The first reference element is the tesseract with edge length $2$ and volume $16$ shown in Figure~\ref{eng_tesseract}. It is straightforward to define bounds of integration as follows
\begin{equation*}
    \int_{-1}^1 \int_{-1}^1 \int_{-1}^1 \int_{-1}^1 dx_1 dx_2 dx_3 dx_4 = 16.
\end{equation*}
Trivially, the orthonormal basis inside of the reference tesseract is given by
\begin{equation*}
    \psi_{ijkq} \left( \boldsymbol{x} \right) =  \hat{P}_i \left( a \right) \hat{P}_j \left( b \right) \hat{P}_k \left( c \right) \hat{P}_q \left( d \right), 
\end{equation*}
where $0 \leq i, j, k, q \leq p$, $a = x_1$, $b = x_2$, $c = x_3$, and $d = x_4$. There are a total of $N_{\text{dof}} = \left(p+1\right)^4$ basis functions.

The tesseract has a total of $384$ symmetries. We can apply these symmetries to a point $(\alpha,\beta,\gamma,\delta)$ with the constraints $0 \leq \alpha,\beta,\gamma,\delta \leq 1$. This will yield a set $\chi(\alpha,\beta,\gamma,\delta)$ containing its counterparts noting that $\chi$ denotes the permutations. The cardinality of $\chi$ depends on if any of the symmetries give rise to identical points. By enumerating this process with the conditions above we obtain twelve symmetry orbits (groups)
\begin{equation*}
\begin{aligned}[c]
S_1 &= (0,0,0,0),\\
S_2 (\alpha) &= \chi(\alpha,0,0,0),\\
S_3 (\alpha) &= \chi(\alpha,\alpha,0,0), \\
S_4 (\alpha,\beta) &= \chi(\alpha,\beta,0,0),\\
S_5 (\alpha) &= \chi(\alpha,\alpha,\alpha,0), \\
S_6 (\alpha,\beta) &= \chi(\alpha,\alpha,\beta,0),\\
S_7 (\alpha,\beta,\gamma) &= \chi(\alpha,\beta,\gamma,0), \\
S_8 (\alpha) &= \chi(\alpha,\alpha,\alpha,\alpha), \\
S_9 (\alpha,\beta) &= \chi(\alpha,\alpha,\alpha,\beta), \\
S_{10} (\alpha,\beta) &= \chi(\alpha,\alpha,\beta,\beta), \\
S_{11} (\alpha,\beta,\gamma) &= \chi(\alpha,\alpha,\beta,\gamma), \\
S_{12} (\alpha,\beta,\gamma,\delta) &= \chi(\alpha,\beta,\gamma,\delta),
\end{aligned}
\qquad \qquad
\begin{aligned}[c]
|S_1| &= 1, \\
|S_2| &= 8, \\
|S_3| &= 24, \\
|S_4| &= 48, \\
|S_5| &= 32, \\
|S_6| &= 96, \\
|S_7| &= 192, \\
|S_8| &= 16, \\
|S_9| &= 64, \\
|S_{10}| &= 96, \\
|S_{11}| &= 192, \\
|S_{12}| &= 384.
\end{aligned}
\end{equation*}

\subsection{Tetrahedral Prism}
\begin{figure}[h!]
\centering
\includegraphics[width=0.7\textwidth]{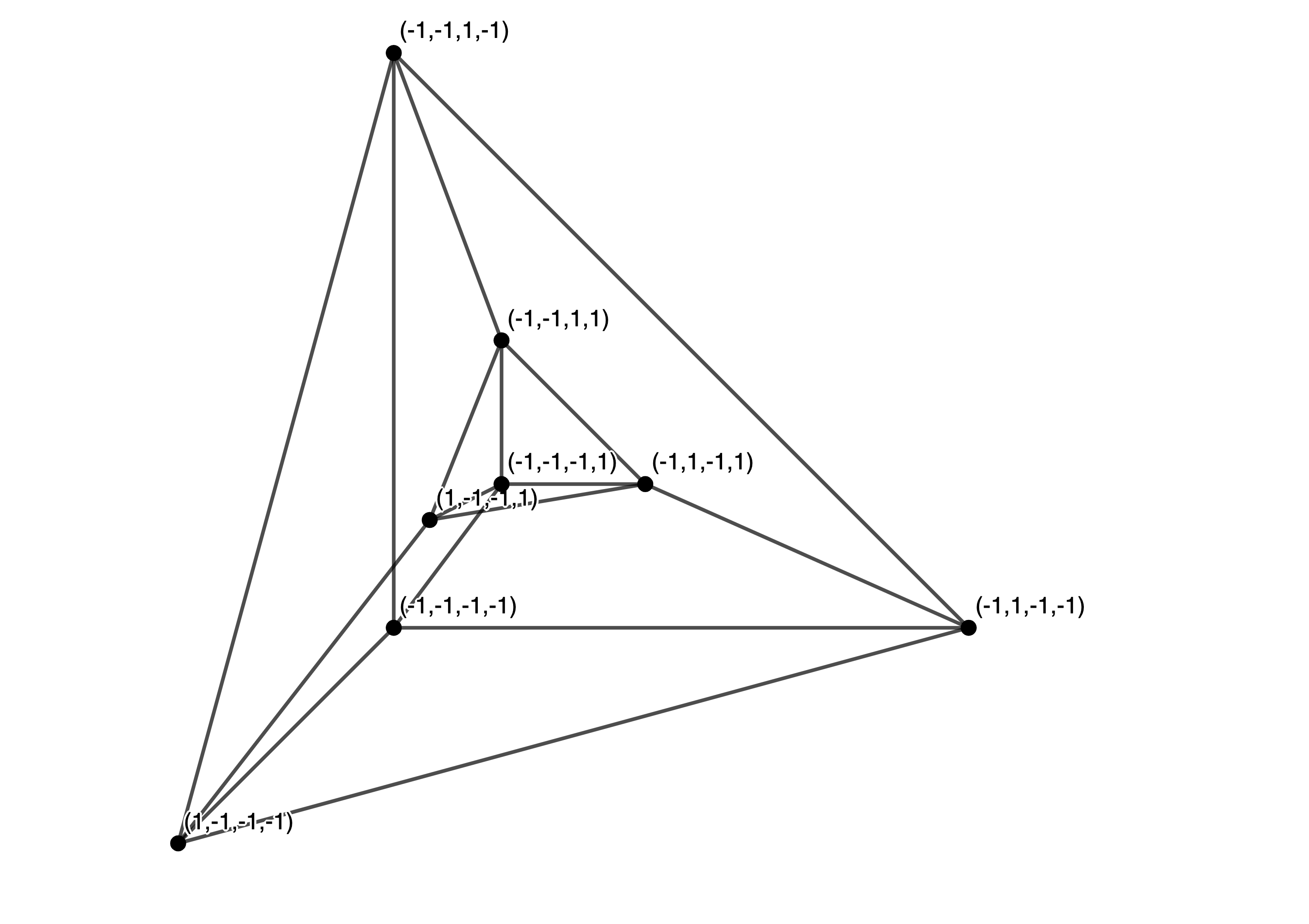}
\caption{Tetrahedral prism reference element.}
\label{eng_tetpri}
\end{figure}
Our second reference element is the tetrahedral prism which is created by an extrusion of a standard tetrahedron into four dimensional space; the element itself is shown in Figure~\ref{eng_tetpri}. (Note: the standard tetrahedron is the three dimensional extension of the right triangle.) The volume of the tetrahedral prism is $\frac{8}{3}$. The integration bounds for the tetrahedral prism are similar to the standard tetrahedron while accounting for the extrusion in the fourth dimension,
\begin{equation*}
    \int_{-1}^{1} \int_{-1}^{1} \int_{-1}^{-x_3} \int_{-1}^{-1-x_2 -x_3 } dx_1 dx_2 dx_3 dx_4 = \frac{8}{3}.
\end{equation*}
Furthermore, the orthonormal polynomial basis for the tetrahedral prism is found by extending the tetrahedron basis in~\cite{warburton2003constants,Hesthaven07,witherden2015identification},
\begin{equation*}
    \psi_{ijkq} \left( \boldsymbol{x} \right) = \sqrt{8} \hat{P}_i \left( a \right) \hat{P}_j^{(2i+1,0)} \left( b \right) \hat{P}_k^{(2i+2j+2,0)} \left( c \right) \hat{P}_q \left( d \right) \left( 1-b \right)^{i} \left( 1-c \right)^{i + j},
\end{equation*}
where $0 \leq i+j+k \leq p$, $0 \leq q \leq p$, $a = -2\frac{1+x_1}{x_2 + x_3}$ - 1, $b= 2\frac{1+x_2}{1-x_3} -1$, $c = x_3$, $d=x_4$, and $N_{\text{dof}} = \left(p+1\right)^2 \left(p+2\right) \left(p+3\right)/6$. 

For this element, it is convenient to work in barycentric coordinates due to the two tetrahedral facets. The barycentric coordinates  are specified as follows
\begin{equation*}
    \boldsymbol{\lambda} = \left( \lambda_1 , \lambda_2 , \lambda_3 , \lambda_4 , x_4 \right)^T  \quad \text{where} \ 0 \leq \lambda_i \leq 1, \quad \lambda_1 + \lambda_2 + \lambda_3 + \lambda_4  =1.
\end{equation*}
They are related to Cartesian coordinates by
\begin{equation*}
    \boldsymbol{x} = \begin{pmatrix}
    -1 & 1  & -1 & -1 & 0 \\
    -1 & -1 & 1  & -1 & 0 \\
    -1 & -1 & -1 & 1  & 0 \\
    0  &  0 &  0 & 0  & 1
    \end{pmatrix} \boldsymbol{\lambda} .
\end{equation*}
The symmetry orbits of the tetrahedral prism are similar to those of the tetrahedron, with an additional parameter that accounts for the extension of the prism in the fourth dimension. There are 48 total symmetries for the tetrahedral prism, and the symmetry orbits are written as follows
\begin{equation*}
\begin{aligned}[c]
S_1  &= (\tfrac{1}{4},\tfrac{1}{4},\tfrac{1}{4},\tfrac{1}{4},0),\\
S_2 (\delta) &= (\tfrac{1}{4},\tfrac{1}{4},\tfrac{1}{4},\tfrac{1}{4},\pm \delta),\\
S_3 (\alpha) &= \chi(\alpha,\alpha,\alpha,1-3\alpha,0),\\
S_4 (\alpha,\delta) &= \chi(\alpha,\alpha,\alpha,1-3\alpha,\pm \delta),\\
S_5 (\alpha) &= \chi(\alpha,\alpha,\tfrac{1}{2} - \alpha,\tfrac{1}{2} - \alpha,0), \\
S_6 (\alpha,\delta) &= \chi(\alpha,\alpha,\tfrac{1}{2} - \alpha,\tfrac{1}{2} - \alpha, \pm \delta), \\
S_7 (\alpha,\beta) &= \chi(\alpha,\alpha,\beta,1-2\alpha-\beta, 0), \\
S_8 (\alpha,\beta,\delta) &= \chi(\alpha,\alpha,\beta,1-2\alpha-\beta, \pm \delta), \\
S_9 (\alpha,\beta,\gamma) &= \chi(\alpha,\beta,\gamma, 1 - \alpha - \beta - \gamma,0), \\
S_{10} (\alpha,\beta,\gamma,\delta) &= \chi(\alpha,\beta,\gamma, 1 - \alpha - \beta - \gamma, \pm \delta),
\end{aligned}
\qquad \qquad
\begin{aligned}[c]
|S_1| &= 1, \\
|S_2| &= 2, \\
|S_3| &= 4, \\
|S_4| &= 8, \\
|S_5| &= 6, \\
|S_6| &= 12, \\
|S_7| &= 12, \\
|S_8| &= 24, \\
|S_9| &= 24, \\
|S_{10}| &= 48,
\end{aligned}
\end{equation*}
where $\alpha$, $\beta$, and $\gamma$ are constrained to ensure that $0 \leq \lambda_i \leq 1$, $\sum_i \lambda_i = 1$ and $0 \leq \delta \leq 1$.

\subsection{Pentatope}
\begin{figure}[h!]
\centering
\includegraphics[width=0.7\textwidth]{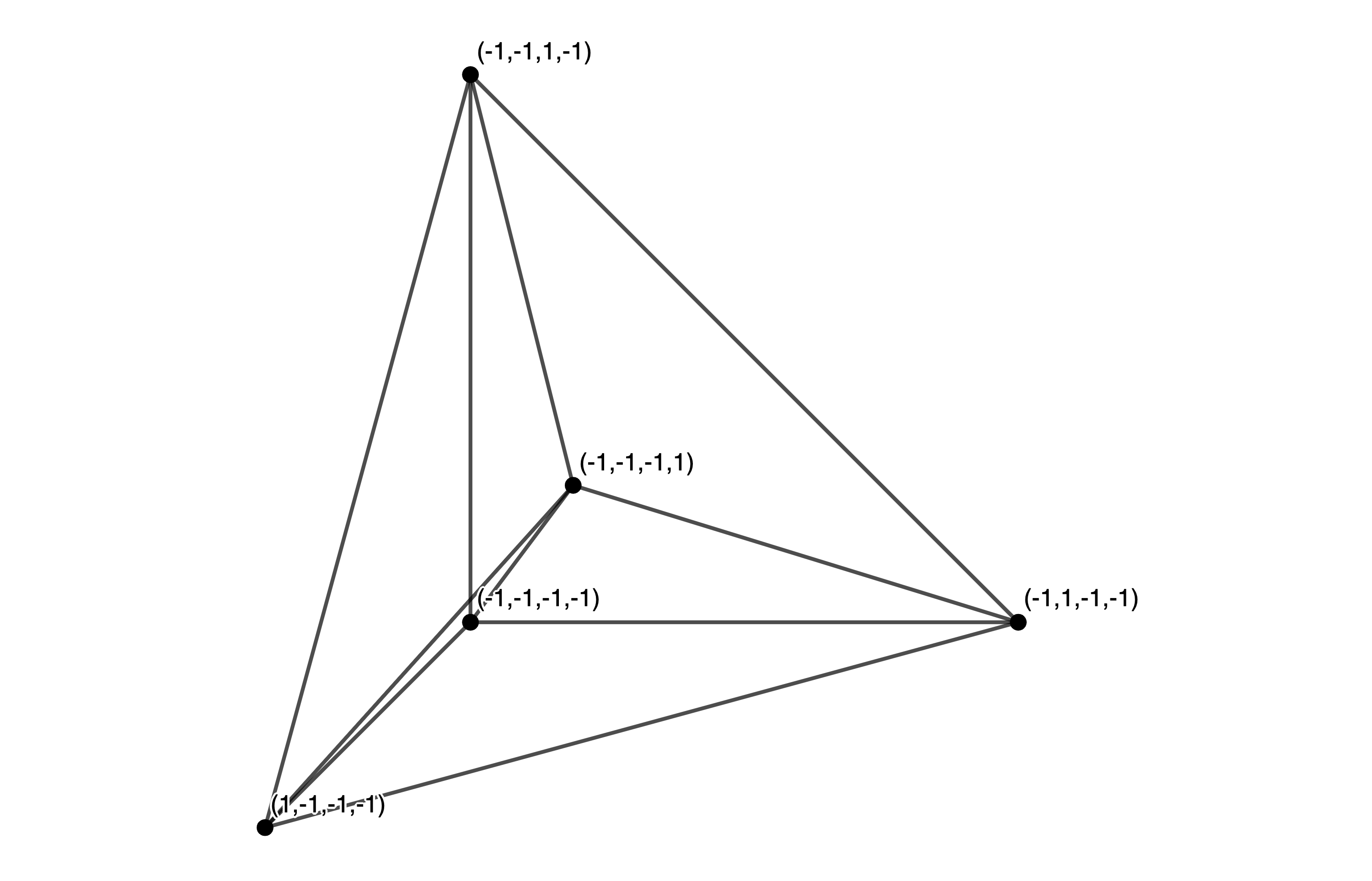}
\caption{Pentatope reference element.}
\label{eng_pentatope}
\end{figure}

The final reference element shown in Figure~\ref{eng_pentatope} is the pentatope or standard four-dimensional simplex which has a volume of $\frac{2}{3}$. The bounds of integration are
\begin{equation*}
    \int_{-1}^1 \int_{-1}^{-x_4} \int_{-1}^{-1-x_3 - x_4} \int_{-1}^{-2-x_2 - x_3 -x_4 } dx_1 dx_2 dx_3 dx_4 = \frac{2}{3}.
\end{equation*}
The orthonormal polynomial basis for the pentatope is obtained from the generalized simplex formula presented in~\cite{warburton2003constants}
\begin{align*}
    \nonumber \psi_{ijkq} \left( \boldsymbol{x} \right) &= 4 \hat{P}_i \left( a \right) \hat{P}_j^{(2i+1,0)} \left( b \right) \hat{P}_k^{(2i+2j+2,0)} \left( c \right) \hat{P}_q^{(2i+2j+2k+3,0)} \left( d \right) \\
    &\times \left( 1-b \right)^{i} \left( 1-c \right)^{i + j} \left( 1-d \right)^{i + j +k},
\end{align*}
where $0 \leq i+j+k+q \leq p$, $a = -2\frac{x_1 +1}{x_2 +x_3 +x_4 +1} - 1$, $b=-2\frac{1+x_2}{x_3 +x_4} - 1$, $c = 2\frac{1+x_3}{1-x_4} - 1$, $d=x_4$, and $N_{\text{dof}} = \left(p+1\right)\left(p+2\right)\left(p+3\right)\left(p+4\right)/24$.
Similar to the tetrahedron, it is convenient to work in barycentric coordinates. For the pentatope these are
\begin{equation*}
    \boldsymbol{\lambda} = \left( \lambda_1 , \lambda_2 , \lambda_3 , \lambda_4 , \lambda_5 \right)^T  \quad \text{where} \ 0 \leq \lambda_i \leq 1, \quad \lambda_1 + \lambda_2 + \lambda_3 + \lambda_4 + \lambda_5 =1,
\end{equation*}
and are related to Cartesian coordinates by
\begin{equation*}
    \boldsymbol{x} = \begin{pmatrix}
    -1 & 1 & -1 & -1 & -1 \\
    -1 & -1 & 1 & -1 & -1 \\
    -1 & -1 & -1 & 1 & -1 \\
    -1 & -1 & -1 & -1 & 1 
    \end{pmatrix} \boldsymbol{\lambda}.
\end{equation*}
Note that the columns of the matrix correspond to the vertices of the pentatope in Cartesian coordinates. The pentatope has a total of 120 symmetries as shown in \cite{maeztu1995consistent}, which take the following form
\begin{equation*}
\begin{aligned}[c]
S_1  &= (\tfrac{1}{5},\tfrac{1}{5},\tfrac{1}{5},\tfrac{1}{5},\tfrac{1}{5}),\\
S_2 (\alpha) &= \chi(\alpha,\alpha,\alpha,\alpha, 1-4\alpha),\\
S_3 (\alpha) &= \chi(\alpha,\alpha,\alpha,\tfrac{1}{2} - \tfrac{3}{2} \alpha,\tfrac{1}{2} - \tfrac{3}{2} \alpha),\\
S_4 (\alpha,\beta) &= \chi(\alpha,\alpha,\alpha,\beta,1 -3\alpha -\beta),\\
S_5 (\alpha,\beta) &= \chi(\alpha,\alpha,\beta,\beta,1-2\alpha-2\beta ), \\
S_6 (\alpha,\beta, \gamma) &= \chi(\alpha,\alpha,\beta,\gamma, 1-2\alpha-\beta-\gamma), \\
S_7 (\alpha,\beta,\gamma,\delta) &= \chi(\alpha,\beta,\gamma,\delta, 1 -\alpha -\beta -\gamma -\delta),
\end{aligned}
\qquad \qquad
\begin{aligned}[c]
|S_1| &= 1, \\
|S_2| &= 5, \\
|S_3| &= 10, \\
|S_4| &= 20, \\
|S_5| &= 30, \\
|S_6| &= 60, \\
|S_7| &= 120, 
\end{aligned}
\end{equation*}
where $\alpha$, $\beta$, $\gamma$, and $\delta$ are constrained to ensure that $0\leq \lambda_i \leq 1, \sum_i \lambda_i = 1$.

\section{Solution Strategy for Quadrature Rules} \label{sec;implementation}

\subsection{Our approach}

In this section, we describe a strategy for generating symmetric quadrature rules on space-time elements. Note that this strategy leverages the numerical infrastructure which was developed in the previous section.

The task of identifying a quadrature rule of strength $\qorder$ inside of a domain $\Omega$ can be framed as a non-linear optimization problem where we seek a set of $N_p$ abscissa $\{\bm{x}_i\}$ and associated weights $\{\omega_i\}$ such that
\begin{equation}\label{eq:mmatch}
    \int_{\Omega} g_{rstv}(\bm{x}) \, d\bm{x} = \sum_{i=1}^{N_p} \omega_i g_{rstv}(\bm{x}_i) \quad \text{for } \quad 0 \le r + s + t + v \le \qorder,
\end{equation}
where $g_{rstv}(\bm{x})=x_1^r x_2^s x_3^t x_4^v$ and $r,s,t,v \ge 0$.  When we evaluate the left hand side of this equation, the resulting system is a \emph{non-linear least squares problem}; albeit one where we seek a solution whose residual is zero.  Moreover, as formulated in Eq.~\eqref{eq:mmatch}, the system is ill-conditioned and fails to enforce the requirements that the resulting rules are symmetric, have positive weights, and have points which are contained within the domain of integration.

The first issue involving ill-conditioning can be addressed by taking
\[
 g_{rstv}(\bm{x}) = \psi_{rstv}(\bm{x}),
\]
where $\psi_{rstv}(\bm{x})$ corresponds to the orthonormal polynomial basis for $\Omega$ as defined in Section~\ref{sec;symmetry}.  Moreover, on account of the orthogonality of the basis, it follows that
\[
 \int_{\Omega} \psi_{rstv}(\bm{x})\, d\bm{x} = c\delta_{0r}\delta_{0s}\delta_{0t}\delta_{0v},
\]
where $c = 1/\psi_{0000}$.  This enables us to forego having to directly evaluate integrals inside of our domain in order to formulate our non-linear least squares problem.

The second issue involving symmetry can be resolved by decomposing the $N_p$ points into a set of symmetry groups (orbits).  For example, consider taking $\Omega$ to be the pentatope and $N_p$ to be $61$.  One possible symmetric arrangement of these points involves a combination of an $S_1$ orbit, an $S_3$ orbit, an $S_4$ orbit, and an $S_5$ orbit. Summing the number of points generated by each orbit we find
\[
 1 + 10 + 20 + 30 = 61,
\]
as required.  By symmetry, it clearly must be the case that all points belonging to an orbit must share the same quadrature weight.  Noting that an $S_1$ orbit has no parameters, $S_3$ orbits have one parameter, and $S_4$ and $S_5$ orbits both have two parameters, our non-linear least squares problem has a total of $9$ parameters; $5$ from the orbits and $4$ from the associated quadrature weights.  This is a substantial reduction compared to the $61\times 4 + 61 = 305$ parameters in our original---asymmetric---problem.  This greatly reduces the difficulty associated with obtaining an exact solution to the non-linear least squares problem.  The caveat here is that, in general, we do not know precisely which of the many distinct orbital decompositions is the correct one.  For example, in the case where $N_p = 61$ in a pentatope, there are $24$ unique decompositions, including the one we've already presented.  Others include a decomposition with one $S_1$ orbit, four $S_2$ orbits, one $S_3$ orbit, and an $S_5$ orbit for a total of $61$ points defined by $14$ parameters, and another decomposition with one $S_1$ orbit and 12 $S_2$ orbits for a total of $61$ points defined by $25$ parameters.  Given that there is no \emph{a priori} means of ascertaining which orbit is most likely to yield a rule, it is generally necessary to try \emph{all} of these decompositions.  Moreover, we remark here that the number of decompositions grows rapidly with $N_p$.  When $N_p = 600$ there are some $37,457$ decompositions inside of the pentatope.

Symmetry also enables us to reduce the number of basis functions which must be included in our non-linear least squares problem.  Specifically, the ability of a rule to exactly integrate one particular basis function oftentimes implies that it is also capable of exactly integrating other basis functions, too.  For example, inside of a tesseract, if a rule is capable of integrating $\psi_{rstv}(\bm{x})$ then, by symmetry, it must also be able to integrate $\psi_{srtv}(\bm{x})$ and all other permutations thereof.  This enables us to substantially reduce the size of the non-linear least squares problem.

The third issue involving positivity of the quadrature weights can be addressed in a variety of ways.  One possibility is to compute them using a \emph{non-negative least squares} routine.  However, the performance of such routines is typically not competitive with those of ordinary linear least squares solvers.  As such, we have found it more efficient to introduce a penalty parameter to the non-linear least squares problem.  Specifically, we add an additional equation into the non-linear system which vanishes if all the weights are positive
\begin{equation*}
 \varepsilon = \sum_{i=1}^{N_p} \frac{\omega_i - |\omega_i|}{2}.
\end{equation*}

The fourth issue involving confinement of quadrature points to  $\Omega$ can also be addressed in a variety of ways. One is to employ a non-linear least squares solver with support for linear equality constraints.  However, such solvers are less widespread than unconstrained non-linear least squares solvers.  An alternative approach (used here), which foregoes the need for explicit constraints, is to simply \emph{clamp} or cut off the relevant orbital parameters before passing them to the $S_i$ orbital functions which yield the Cartesian points.  This clamping discourages the solver from attempting to take orbital parameters outside of their respective ranges.

Another final improvement we can make to our baseline algorithm is to observe that the non-linear least squares problem is \emph{separable}.  If the orbital parameters---which define the abscissa---are known then our non-linear least squares problem in Eq.~\eqref{eq:mmatch} reduces down to a \emph{linear} least squares problem for the $\omega_i$.  Thus, it is possible to treat the quadrature weights as \emph{dependent} variables.  This enables us to further reduce the size of the non-linear least squares problem; albeit at the cost of having to solve a linear least squares problem at each iteration.

\subsection{Alternative approaches}

We remark here that in the case of the tesseract that it is possible to generate symmetric rules through a tensor-production construction of a suitable one-dimensional rule.  If a Gauss--Legendre rule with $l + 1$ points is employed, it follows that the resulting tesseract rule has $(l + 1)^4$ points and is capable of exactly integrating $g_{rstv}(\bm{x})$ for $0 \le r, s, t, v \le 2l + 1$. 
This property leads to the rules being extremely efficient for integrating functions with a  tensor-product structure.  Moreover, when paired with sum-factorization they can enable computationally efficient spectral element discretizations~\cite{orszag1979spectral}.

It is also possible to apply these tensor product rules to tetrahedral prisms and pentatopes by employing the well-known Duffy transformation~\cite{duffy1982quadrature}. In what follows, we briefly review the Duffy-based procedure for generating quadrature rules on the pentatope. One is encouraged to consult the appendix of~\cite{lehrenfeld2015nitsche} for more details about the origins of the Duffy transformation.

We begin by introducing $\bm{y} = \left(y_1,y_2,y_3 \right) \in \mathbb{R}^{3}$ and $\bm{x} = \left(\bm{y}, t \right) \in \mathbb{R}^{4}$. Next, we construct the following linear map
\begin{align*}
    \left(\bm{y},t\right) \rightarrow \left(\frac{2\bm{y} + \left(1+t\right)\bm{1}}{1-t},t \right) = \left(\bm{y}',t\right),
\end{align*}
where $\bm{1} = \left(1,1,1\right)$. Thereafter, this map can be used to reformulate a generic integral on the reference pentatope as follows
\begin{align}
    \nonumber &\int_{-1}^{1} \int_{-1}^{-t} \int_{-1}^{-1-t-y_3} \int_{-1}^{-2-t-y_3-y_2} f\left(\bm{y},t\right) dy_1\, dy_2 \, dy_3 \, dt \\[1.5ex]
    \nonumber &= \frac{1}{8} \int_{-1}^{1} \left(1-t\right)^{3} \int_{-1}^{1} \int_{-1}^{-y_{3}'} \int_{-1}^{-1-y_{3}'-y_{2}'} f\left(\frac{\left(1-t\right)\bm{y}'-\left(1+t\right)\bm{1}}{2},t\right) dy_{1}' \, dy_{2}' \, dy_{3}' \, dt \\[1.5ex]
    \nonumber & = \frac{1}{8} \int_{-1}^{1} \left(1-t\right)^{3} \int_{-1}^{1} \int_{-1}^{-y_{3}'} \int_{-1}^{-1-y_{3}'-y_{2}'} g\left(\bm{y}',t\right) d\bm{y}' dt \\[1.5ex]
    &= \frac{1}{8} \int_{-1}^{1} \left(1-t\right)^{3} G\left(t\right) dt, \label{duffy_one}
\end{align}
where on the last two lines we have introduced the following quantities
\begin{align}
    \nonumber g\left(\bm{y}',t\right) &= f\left(\frac{\left(1-t\right)\bm{y}'-\left(1+t\right)\bm{1}}{2},t\right),
    \\[1.5ex]
    G\left(t\right) &= \int_{-1}^{1} \int_{-1}^{-y_{3}'} \int_{-1}^{-1-y_{3}'-y_{2}'} g\left(\bm{y}',t\right) d\bm{y}'. \label{duffy_two}
\end{align}
The key observation is that the integral in Eq.~\eqref{duffy_two} is a three-dimensional integral over the standard tetrahedron. Therefore, one may utilize a generic three-dimensional integration formula to evaluate Eq.~\eqref{duffy_two}, in conjunction with a 1D Gauss-Legendre rule for evaluating Eq.~\eqref{duffy_one}. In this fashion, the Duffy transformation replaces a four-dimensional integration procedure with a combination of three- and one-dimensional procedures. Furthermore, the transformation can be repeated recursively until only 1D Gauss-Legendre rules are required. In this case, the Duffy transformation effectively maps a Gauss-Legendre tensor-product rule from the standard tesseract on to the standard pentatope. 

Somewhat unsurprisingly, none of the rules which arise from the Duffy transformation are optimal. In addition, although such rules are functional and easy to generalize, they are asymmetric and suffer from an undesirable clustering of points.


\section{Summary of Quadrature Rules} 
\label{sec;summary}

Using the optimisation procedure outlined in Section~\ref{sec;implementation}, we obtained a series of symmetric quadrature rules on the tesseract, tetrahedral prism, and pentatope. The implementation was carried out by modifying the polyquad code, an open-source, parallel, highly robust quadrature-generating code written in C++ (see~\cite{witherden2015identification} for details). The resulting quadrature rules, which span orders 2 to 16, are summarised in Table~\ref{tab:rules}. In addition, Table~\ref{tab:example} contains an example of a new rule with strength $\mathcal{P} = 9$ and number of points $N_p = 151$ on the pentatope. All of the rules---\textbf{which are provided as electronic supplementary material}---have positive weights with none of the points being outside of the domain. Looking at Table~\ref{tab:rules}, there are several points of note.  Firstly, due to symmetry all of the tesseract rules are capable of exactly integrating \emph{any} odd function---including odd-degree polynomials.  A consequence of this is that all of the even rules pick up an extra order of accuracy.  Secondly, we also note that the rules in the tetrahedral prism are limited to $\qorder = 14$, whereas those inside of the tesseract and the pentatope go up to $\qorder = 16$.  This is due to the fact that the number of unique topological arrangements of points is an order of magnitude higher in the tetrahedral prism than in the other elements.  As such, it is necessary to consider a greater number of non-linear least squares problems, which in turn, translates into an increased computational~cost.

\begin{table}[h!]
    \centering
    \begin{tabular}{rrrrrrr} \toprule
    & \multicolumn{3}{c}{$N_p$} \\ \cmidrule{2-4}
    $\qorder$ & Tesseract & Pentatope & Tetrahedral prism\\ \midrule
    2 & 16 & 5 & \underline{6} \\
    3 & 16 & 15 & \underline{12}\\
    4 & 24 & 20 & \underline{20}\\
    5 & 24 & 30 & \underline{27}\\
    6 & 57 & 56 & \underline{61}\\
    7 & 57 & \underline{70} & \underline{72}\\
    8 & 160 & \underline{105} & \underline{114}\\
    9 & 160 & \underline{151} & \underline{159}\\
    10 & \underline{272} & \underline{210} & \underline{259}\\
    11 & \underline{272} & \underline{281} & \underline{322}\\
    12 & \underline{512} & \underline{445} & \underline{468}\\
    13 & \underline{512} & \underline{555} & \underline{608}\\
    14 & \underline{728} & \underline{725} & \underline{921}\\
    15 & \underline{728} & \underline{905}\\
    16 & \underline{1384} & \underline{1055} \\
    \bottomrule
    \end{tabular}
    \caption{Number of points $N_p$ required for a rule of strength $\qorder$ inside of each element type. Note that each newly discovered quadrature rule has been underlined. In addition, the explicit parameters (weights and abscissa) for each rule are provided in the electronic supplemental material, in quadruple precision.}
    \label{tab:rules}
\end{table}

\begin{table}[h!]
    \centering
    \begin{tabular}{r|l|l} \toprule
    Orbit & Abscissa & Weight \\ 
    \hline
    $S_1$ & 0.20000000000000000000000000000000 & 0.026283450664919790544554931007647 \\
    \hline
    $S_2$ & 0.033396077194956640848502386475095 & 0.00095166545512843622669550284918655 \\
    \hline
    $S_2$ & 0.24389669924344346067467117680514 & 0.010657347396678176509183437341993 \\
    \hline
    $S_3$ & 0.0050791921061062008978905902787788 & 0.00031411345650965132958719133177732 \\
    \hline
    $S_3$ & 0.2967495659603128481315951580135 & 0.011334029741446587013061345750027 \\
    \hline
    $S_4$ & 0.039703540947493234138821823418935 & 0.0034567258020126958766999595741076 \\
    & 0.18556651270749874050254409289523 & \\
    \hline
    $S_4$ & 0.098399070565063719667340340181978 & 0.0089501016561121758766577213589523 \\
    & 0.25309277814594493439811072010579 & \\
    \hline
    $S_4$ & 0.15177815296591115799210754308869 & 0.0042400257059868015877992812892785 \\
    & 0.011737624799814566830560054028264 & \\
    \hline
    $S_5$ & 0.0088561456625601727465721736167059 & 0.0015526370637512027432010914510579 \\
    & 0.21639650331182284793855392295174 & \\
    \hline
    $S_5$ & 0.085781201932866919896943973150337 & 0.0028780181522793959965686558635525 \\
    & 0.41353118973888920995638429357091 & \\
    \bottomrule
    \end{tabular}
    \caption{The $\mathcal{P} = 9$ quadrature rule on the pentatope with $N_p = 151$ points. The remaining rules are provided in the electronic supplemental material, in quadruple precision.}
    \label{tab:example}
\end{table}

\pagebreak
\clearpage

\section{Numerical Experiments}
\label{sec;experiments}

In order to validate the quadrature rules developed in Section~\ref{sec;summary}, we will now perform a set of numerical experiments on the rules of strength 6-12. First, we will validate that the quadrature rules of degree $\qorder$ exactly integrate polynomials of total order $\porder \leq \qorder$, to within machine precision.
Thereafter, we will evaluate the ability of the quadrature rules to integrate transcendental functions reliably, with well-behaved approximation capabilities.
All experiments in the following section are performed using quad-precision numerical computations to avoid precision issues.

\subsection{Integration of polynomial functions on solitary elements}

Let us consider the integration of polynomials of total order $\porder$, which can be given by
\begin{equation*}
    f_\mathrm{poly}(\bm{x}; \porder)= \sum_{(r, s, t, v)= 0}^{r + s + t + v \leq \porder} C_{rstv} g_{rstv}(\bm{x}) ,
\end{equation*}
where $C_{rstv}$ is a constant factor. In this work, $C_{rstv}$ is generated procedurally by choosing for each $(r, s, t, v)$ a value from a standard normal distribution (see the detailed discussion in~\cite{williams2020family}).

\begin{figure}[h!]
\centering
\includegraphics[width=0.95 \textwidth]{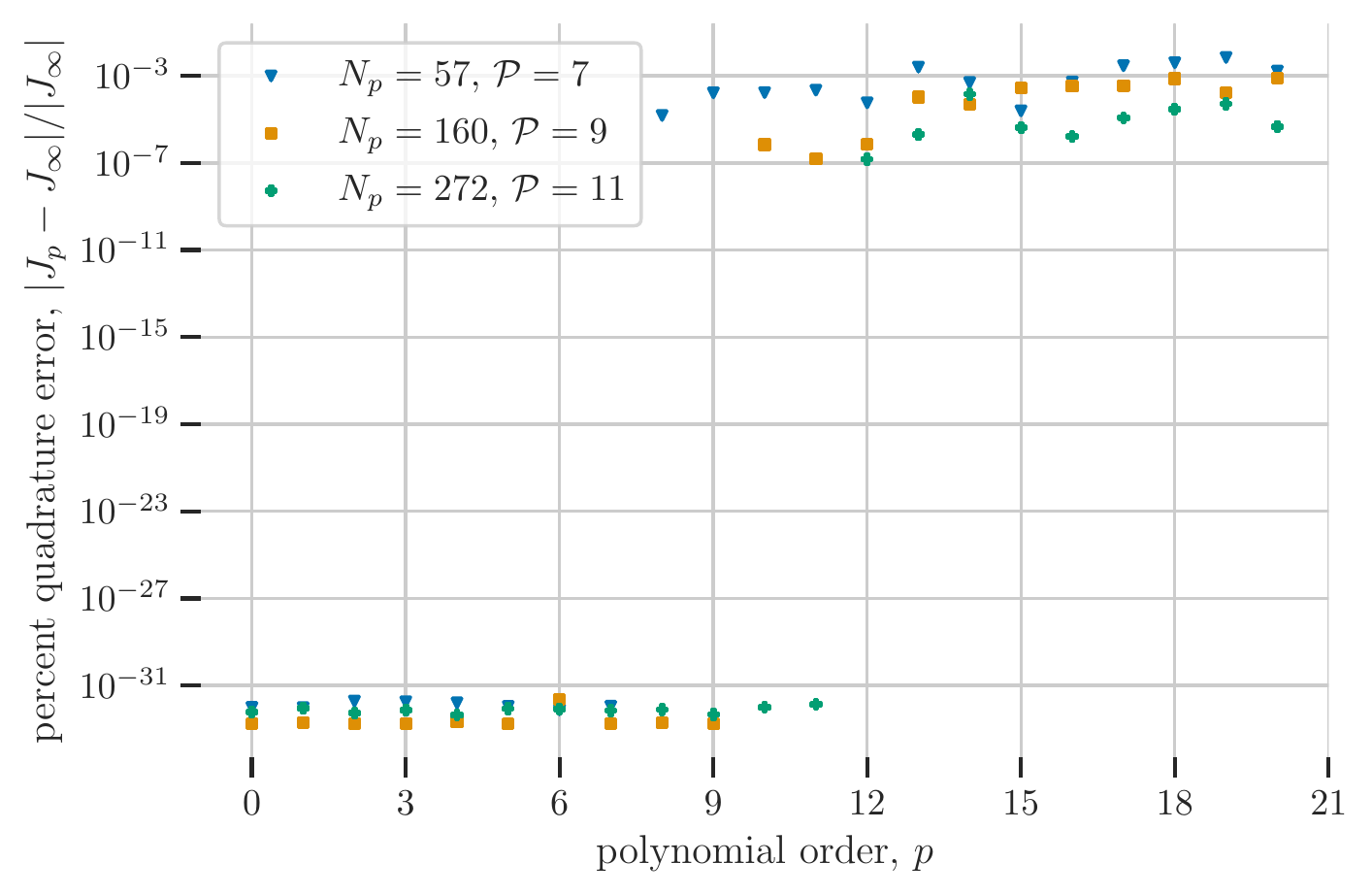}
\caption{Percent error in polynomial integration on a tesseract element due to quadrature.}
\label{quaderror_elem_tesseract}
\end{figure}

\begin{figure}[h!]
\centering
\includegraphics[width=0.95 \textwidth]{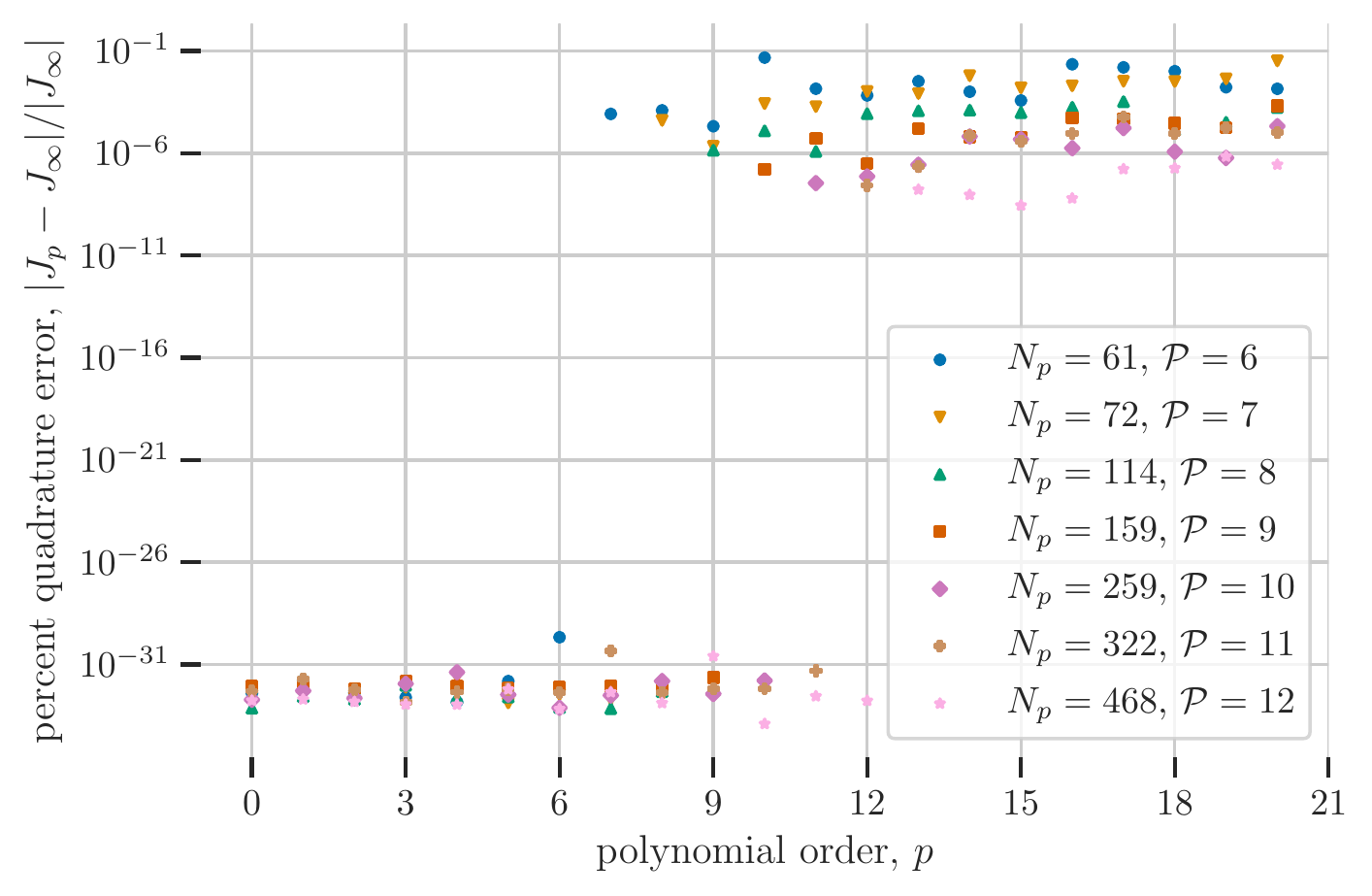}
\caption{Percent error in polynomial integration on a tetrahedral prism element due to quadrature.}
\label{quaderror_elem_prism}
\end{figure}

\begin{figure}[h!]
\centering
\includegraphics[width=0.95 \textwidth]{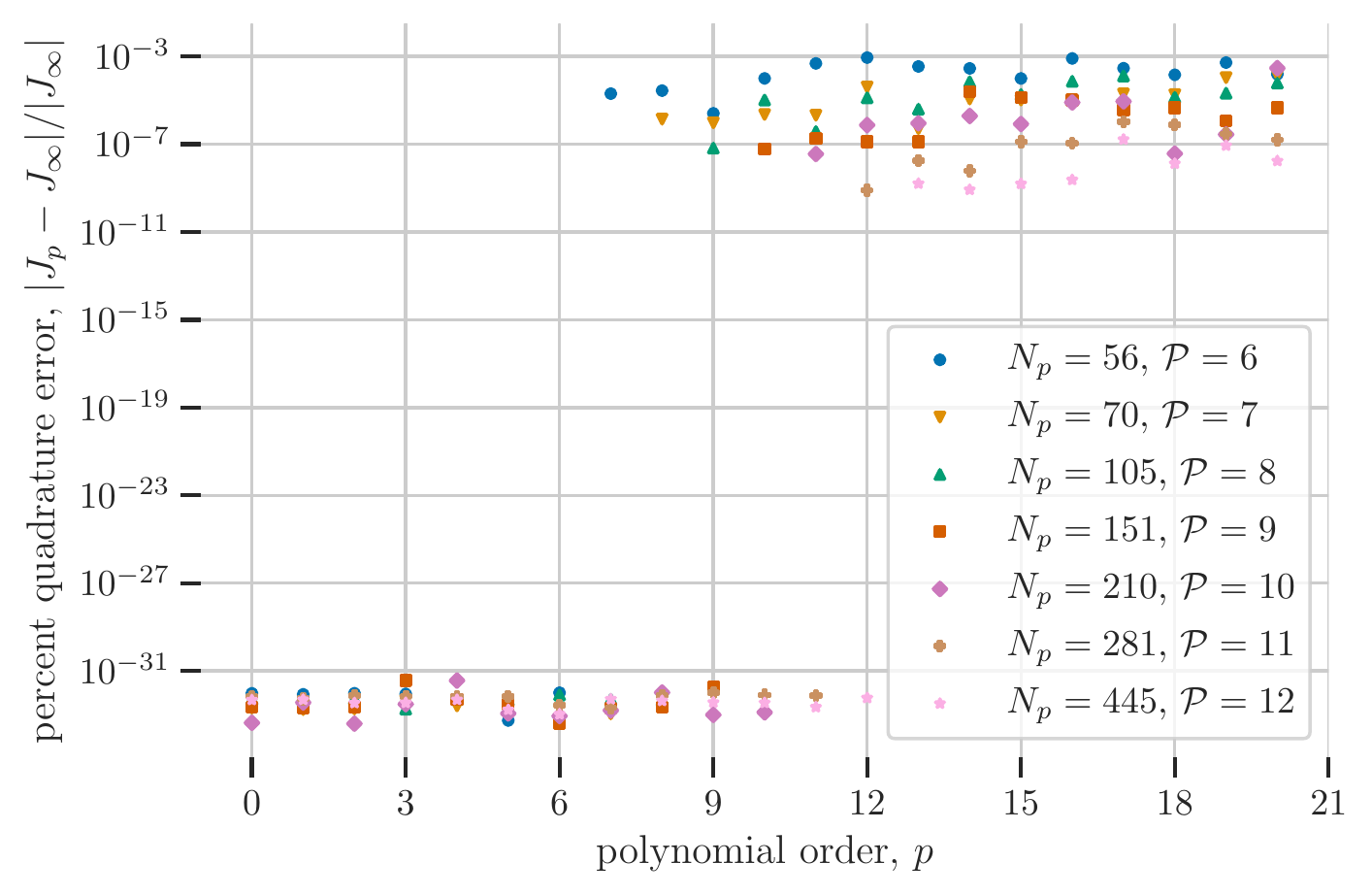}
\caption{Percent error in polynomial integration on a pentatope element due to quadrature.}
\label{quaderror_elem_pentatope}
\end{figure}

For each of the tesseract, tetrahedral prism, and pentatope elements, we have taken the quadrature rules of order $\qorder$ and used them to estimate the integrals of polynomials $f_\mathrm{poly}$ of order $\porder$ as follows
\begin{equation*}
    J = \sum_{i= 1}^{N_p} \omega_i f_\mathrm{poly} (\bm{x}_i) .
\end{equation*}
For each monomial term, $g_{rstv}$, we can write the exact integral $I_{g_{rstv}}$ on an element as follows
\begin{equation*}
    I_{g_{rstv}}= \int_{\kappa} g_{rstv} \left(\bm{x}\right) \, d \bm{x}=
    \begin{cases}
        \frac{1}{(r + 1)(s + 1)(t + 1)(v + 1)} & \text{$\kappa$ is a tesseract} \\
        \frac{s! t! v!}{(r + 1)(s + t + v + 3)!} & \text{$\kappa$ is a tetrahedral prism} \\
        \frac{r! s! t! v!}{(r + s + t + v + 4)!} & \text{$\kappa$ is a pentatope}
    \end{cases}
\end{equation*}
Thus, the exact integral of a polynomial that is being approximated by $J$ can be given analytically by
\begin{equation*}
    J_\infty(\bm{x}; \porder)= \sum_{(r, s, t, v)= 0}^{r + s + t + v \leq \porder} C_{rstv} I_{g_{rstv}} .
\end{equation*}
In Figures~\ref{quaderror_elem_tesseract}, \ref{quaderror_elem_prism}, and~\ref{quaderror_elem_pentatope}, the percent error produced by the quadrature rules on a single element, for each of the element types is shown as a function of the polynomial order $\porder$.
In these plots, it can be seen that the quadrature rules guarantee exactness to machine precision for quadruple precision floating point numbers whenever $\qorder \geq \porder$, as we expect. 

\subsection{Integration of transcendental functions on multi-element grids}

Next, we seek to quantify the performance of the quadrature rules in situations that are immediately relevant to finite element applications. We will consider a domain $\Omega= [0, 1]^4$, which we will tessellate with three uniform triangulations of $\Omega$: $\mathcal{T}_\mathrm{tesseract}$, $\mathcal{T}_\mathrm{prism}$, and $\mathcal{T}_\mathrm{pentatope}$. For $\mathcal{T}_\mathrm{tesseract}$, each dimension of the 4-cube domain is subdivided $m \in \{ 1, 2, 3, \ldots, 12 \}$ times, into $N_\mathrm{elem}= m^4$ tesseract elements. In the case of the prismatic and pentatope grids, each tesseract element is further subdivided into either tetrahedral prisms or pentatopes by using the Kuhn-Freudenthal triangulation on the first three or all four dimensions of the tesseract element, respectively. This results in $N_\mathrm{elem}= 3! m^4$ tetrahedral prism or $N_\mathrm{elem}= 4!m^4$ pentatope elements respectively. In each of these three cases, the integral of a function $f \in L_1 \left(\Omega\right)$ is then computed by
\begin{equation*}
    J_\infty= \int_\Omega f(\bm{x}) \mathrm{~d} \bm{x} \approx J \equiv \sum_{\kappa \in \mathcal{T}} \sum_{i= 1}^{N_p} \omega_i^{(\kappa)} f \left(\bm{x}_i^{(\kappa)} \right) ,
\end{equation*}
where $\bm{x}_i^{(\kappa)}$ and $\omega_i^{(\kappa)}$ are transformed  from the appropriate reference element to each physical element.

In order to understand the performance of the quadrature rules, we will evaluate their use on three transcendental functions, namely an exponential function:
\begin{equation*}
    f_1= \exp(x^2 + 2 y^3 + 3 z^4 + 4 t^5) ,
\end{equation*}
a sinusoidal function:
\begin{equation*}
    f_2= \sin(x^2 + 2 y^3 + 3 z^4 + 4 t^5) ,
\end{equation*}
and an even function (also sinusoidal):
\begin{equation*}
    f_3= \sin(x^2 + y^2 + z^2 + t^2) .
\end{equation*}

\begin{figure}[h!]
\centering
\includegraphics[width=0.95 \textwidth]{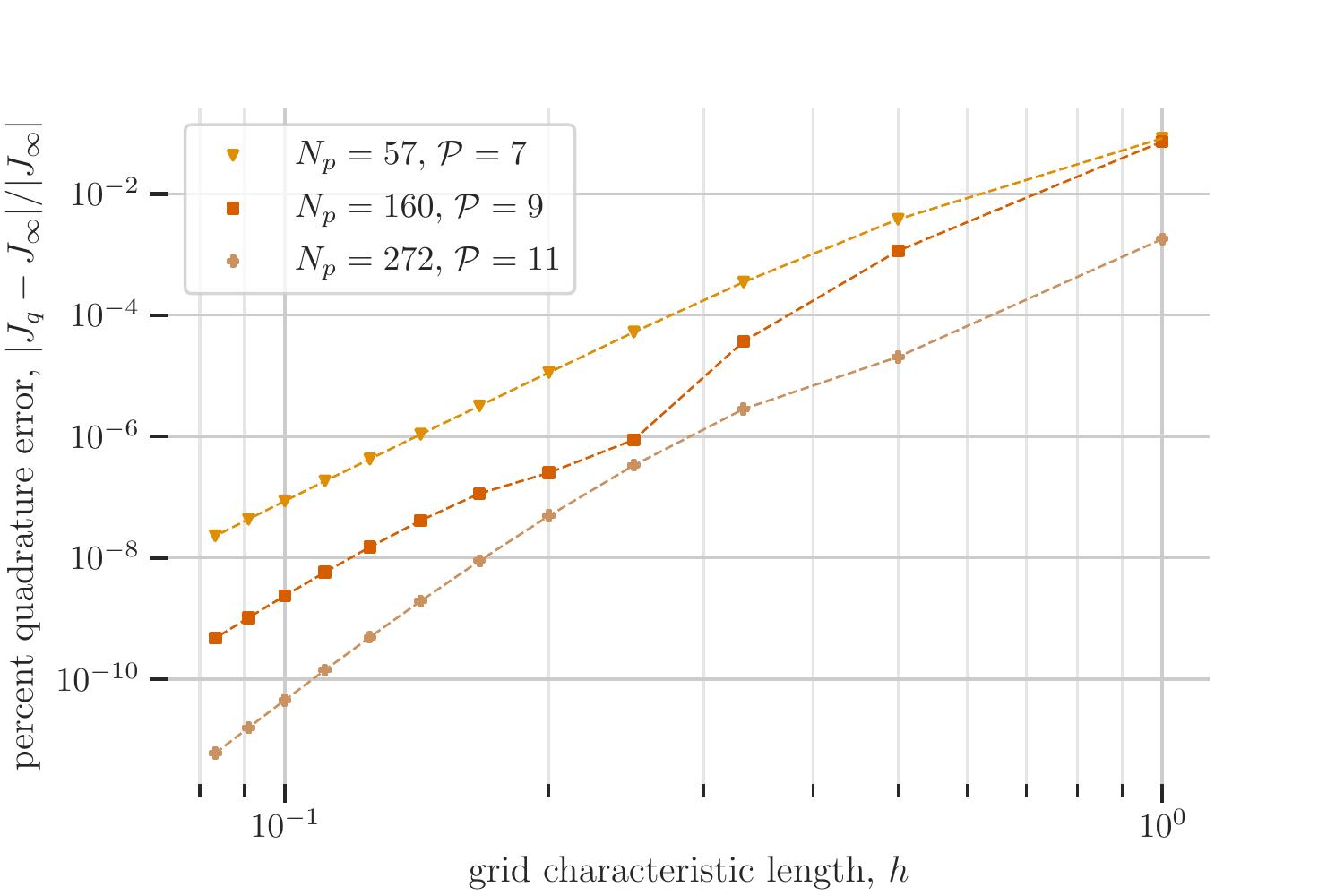}
\caption{Percent error in integration of $f_1$ on a tesseract grid due to quadrature.}
\label{quaderror_grid_tesseract_exp}
\end{figure}

\begin{figure}[h!]
\centering
\includegraphics[width=0.95 \textwidth]{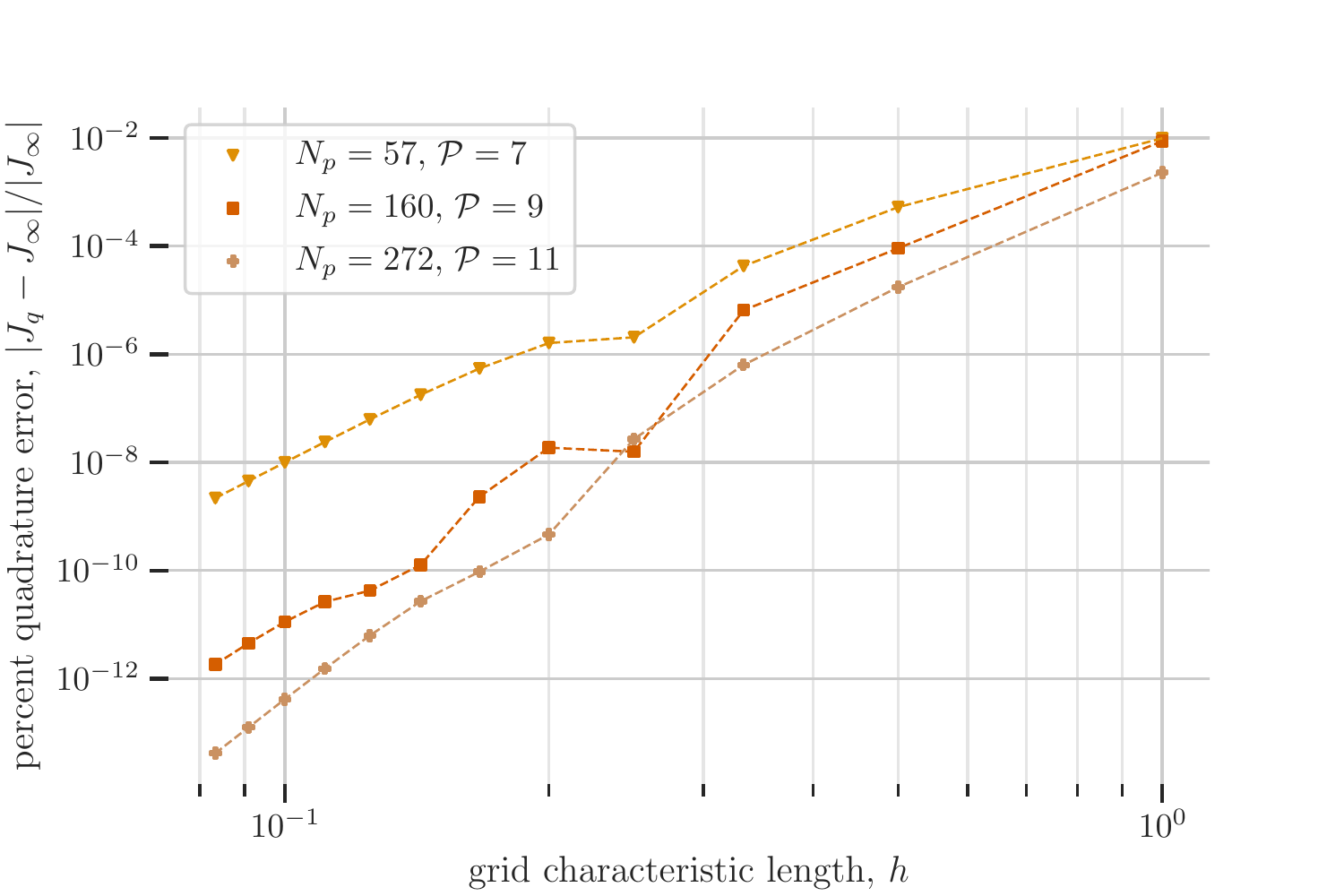}
\caption{Percent error in integration of $f_2$ on a tesseract grid due to quadrature.}
\label{quaderror_grid_tesseract_sin}
\end{figure}

\begin{figure}[h!]
\centering
\includegraphics[width=0.95 \textwidth]{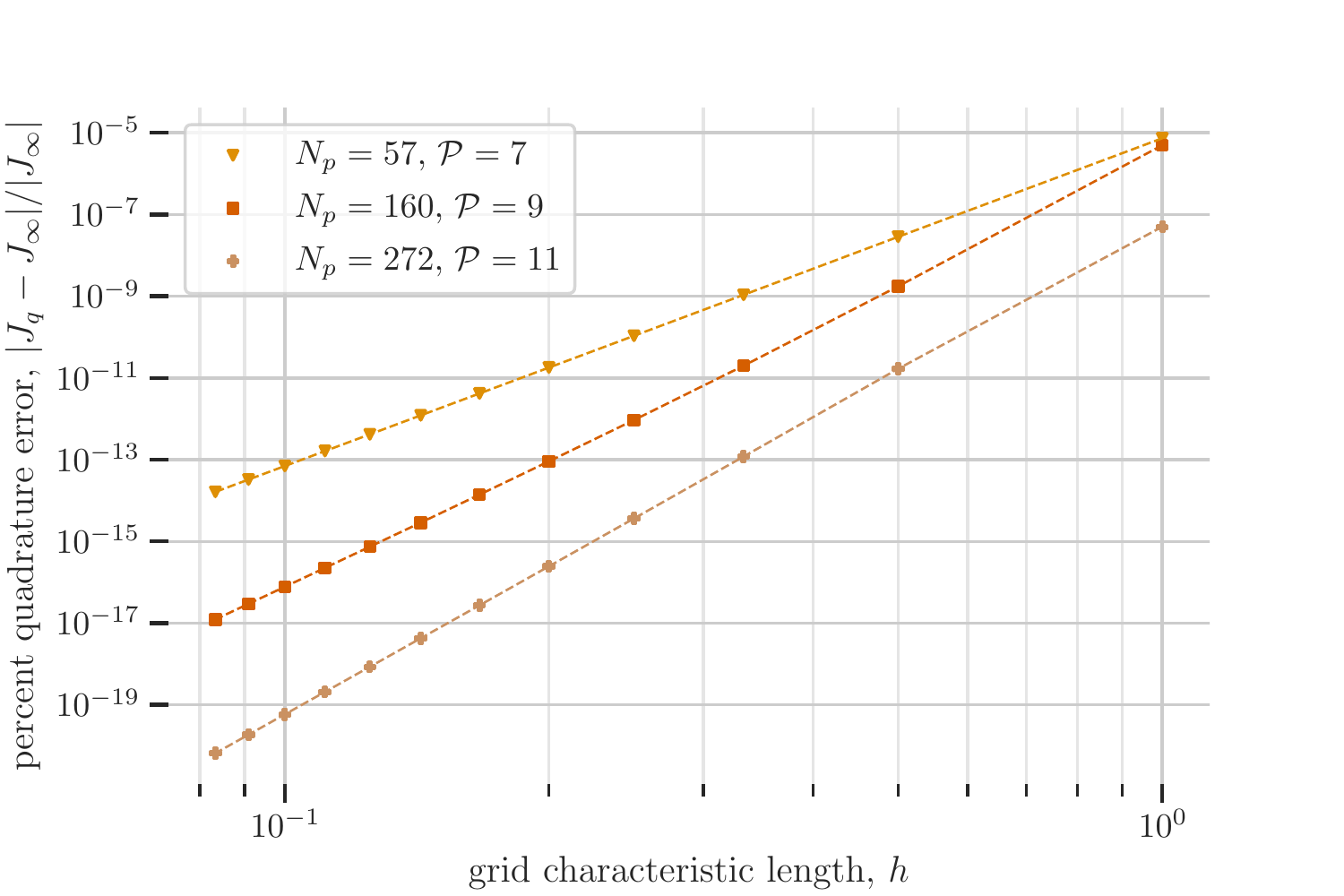}
\caption{Percent error in integration of $f_3$ on a tesseract grid due to quadrature.}
\label{quaderror_grid_tesseract_odd}
\end{figure}

In Figures~\ref{quaderror_grid_tesseract_exp}, \ref{quaderror_grid_tesseract_sin}, and~\ref{quaderror_grid_tesseract_odd}, the percent error induced by the quadrature approximation for each function is plotted with respect to the grid characteristic length
\begin{equation*}
    h= \sqrt[4]{\frac{|\Omega|}{N_\mathrm{elem}}} ,
\end{equation*}
with $|\Omega|$ the volume of the domain and $N_\mathrm{elem}$ the number of elements.

\begin{figure}[h!]
\centering
\includegraphics[width=0.95 \textwidth]{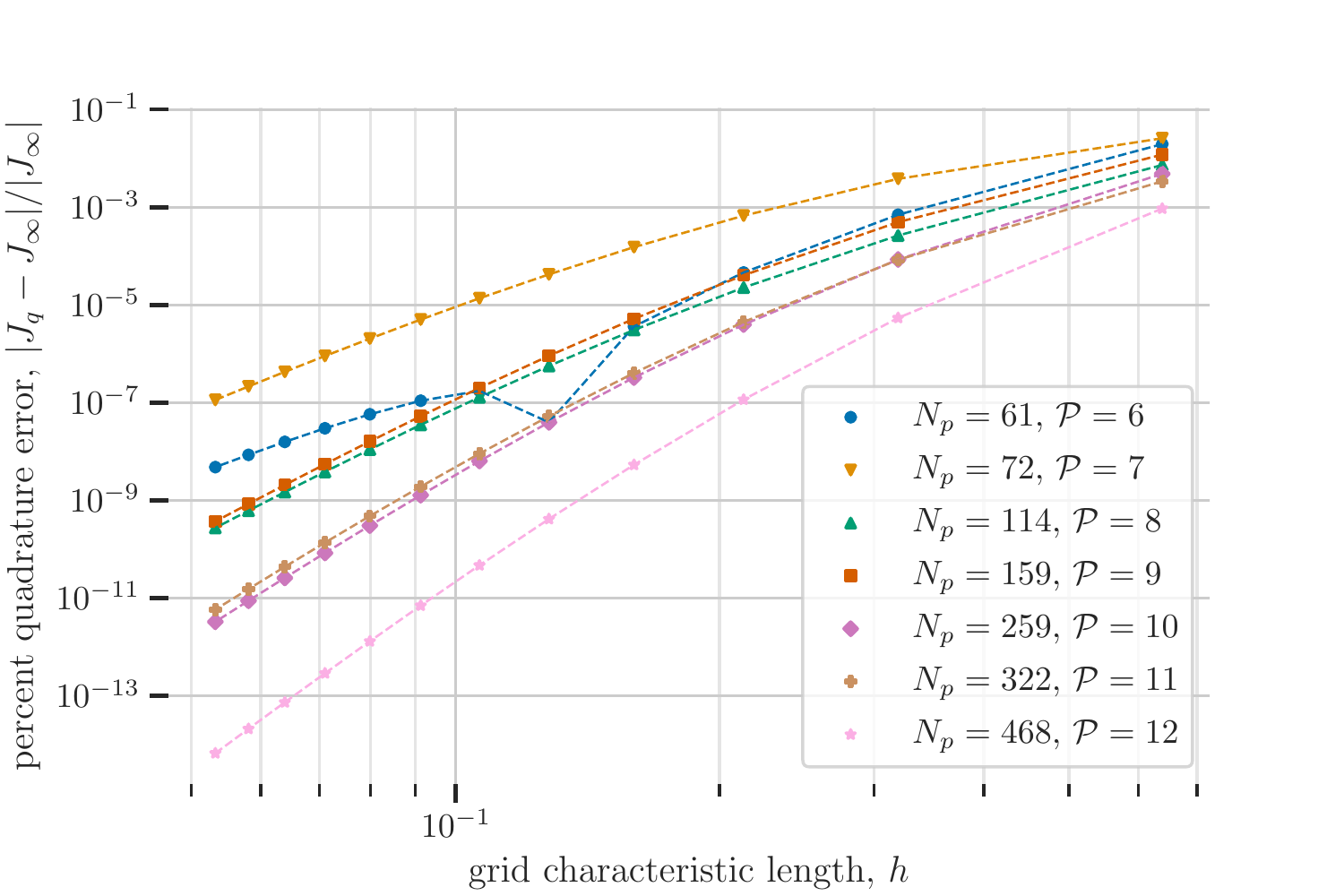}
\caption{Percent error in integration of $f_1$ on a prismatic grid due to quadrature.}
\label{quaderror_grid_prism_exp}
\end{figure}

\begin{figure}[h!]
\centering
\includegraphics[width=0.95 \textwidth]{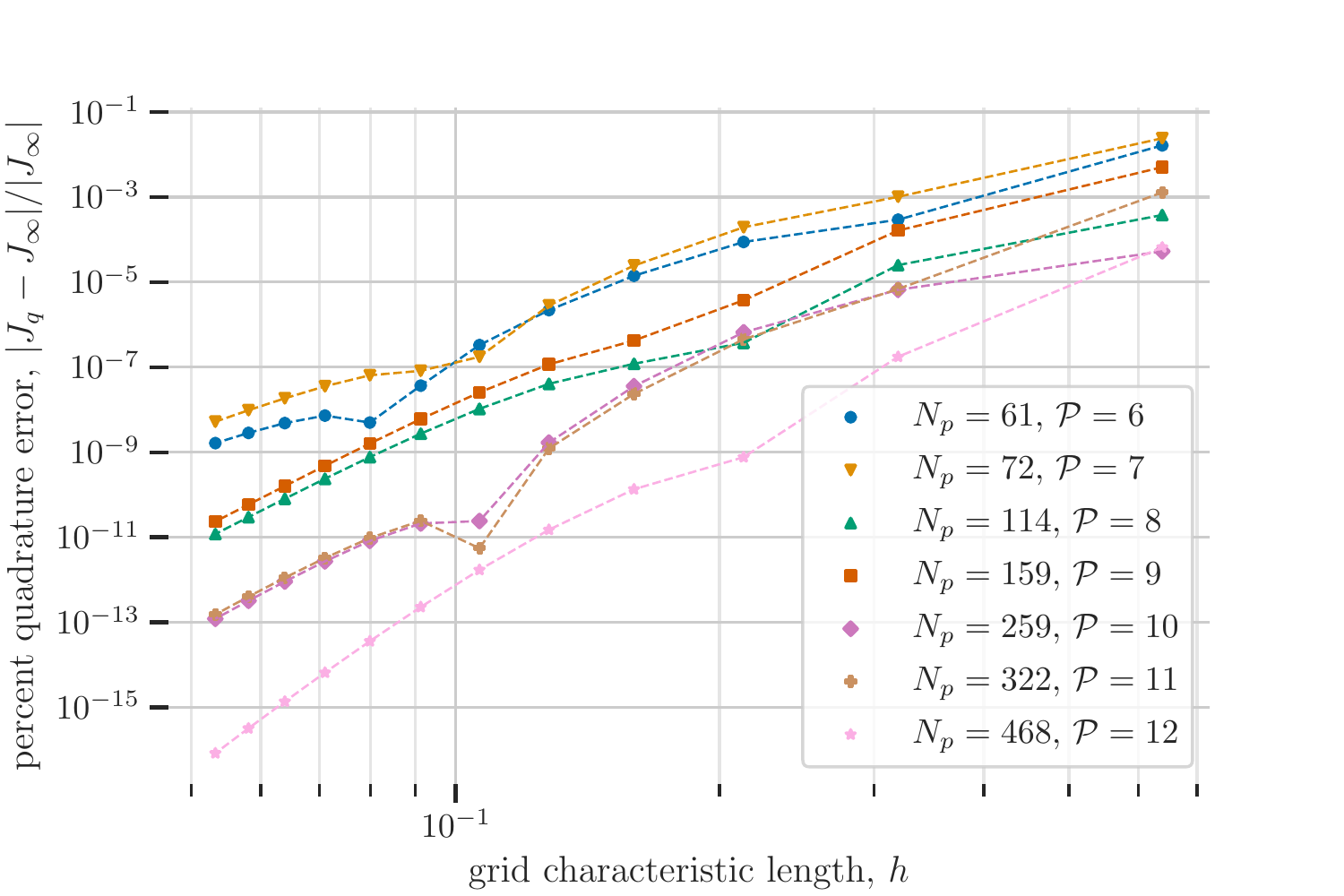}
\caption{Percent error in integration of $f_2$ on a prismatic grid due to quadrature.}
\label{quaderror_grid_prism_sin}
\end{figure}

\begin{figure}[h!]
\centering
\includegraphics[width=0.95 \textwidth]{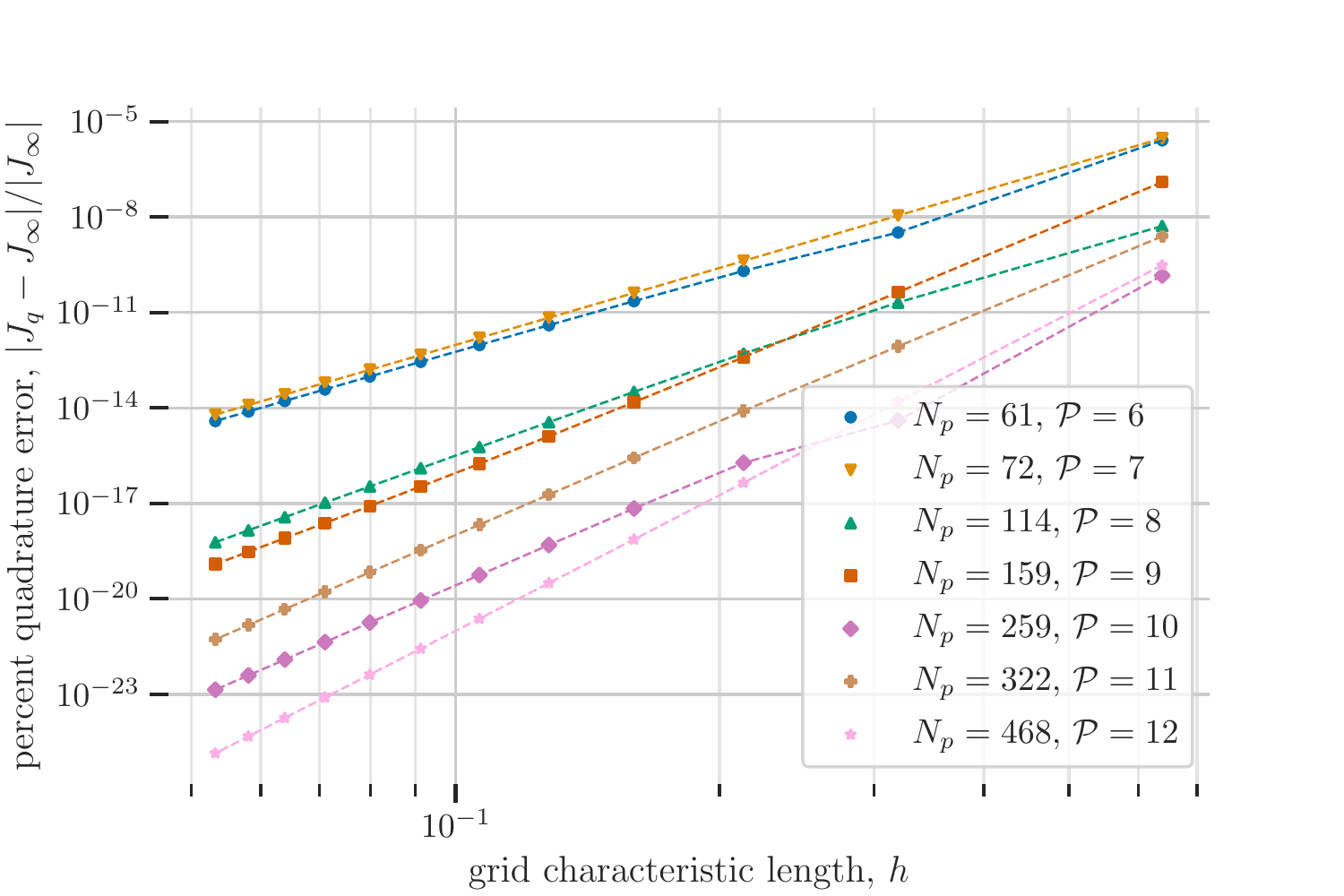}
\caption{Percent error in integration of $f_3$ on a prismatic grid due to quadrature.}
\label{quaderror_grid_prism_odd}
\end{figure}

\begin{figure}[h!]
\centering
\includegraphics[width=0.95 \textwidth]{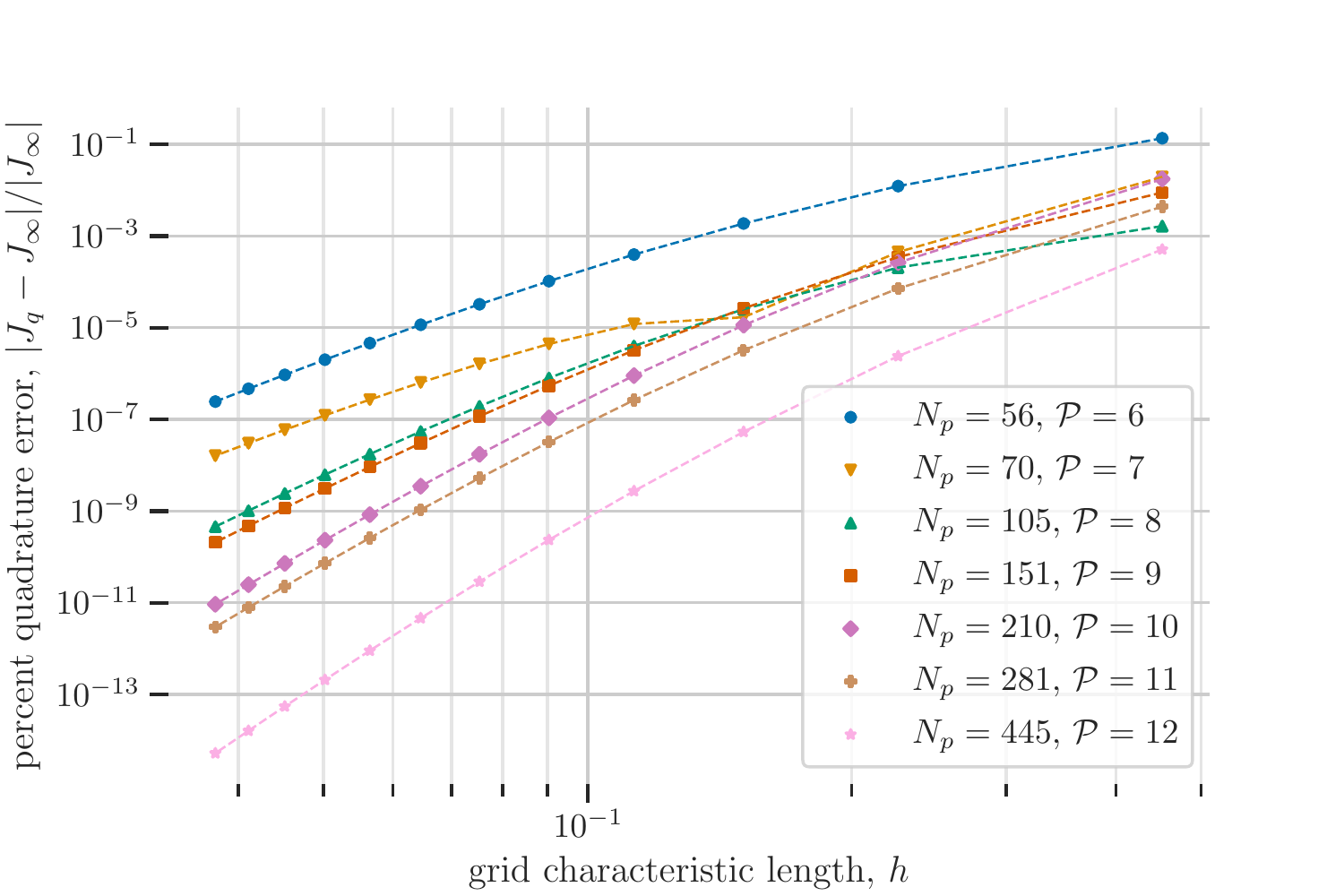}
\caption{Percent error in integration of $f_1$ on a pentatope grid due to quadrature.}
\label{quaderror_grid_pentatope_exp}
\end{figure}

\begin{figure}[h!]
\centering
\includegraphics[width=0.95 \textwidth]{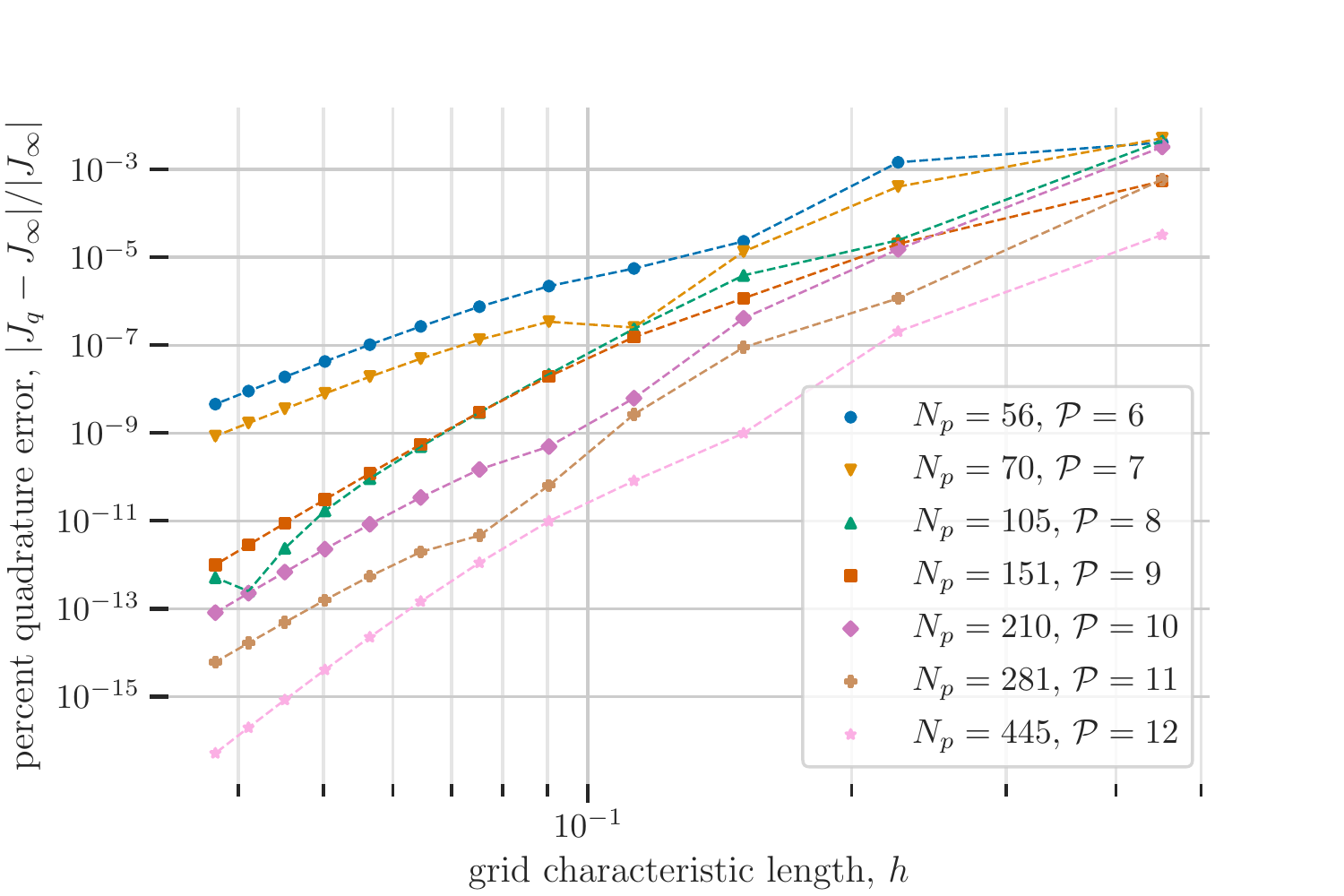}
\caption{Percent error in integration of $f_2$ on a pentatope grid due to quadrature.}
\label{quaderror_grid_pentatope_sin}
\end{figure}

\begin{figure}[h!]
\centering
\includegraphics[width=0.95 \textwidth]{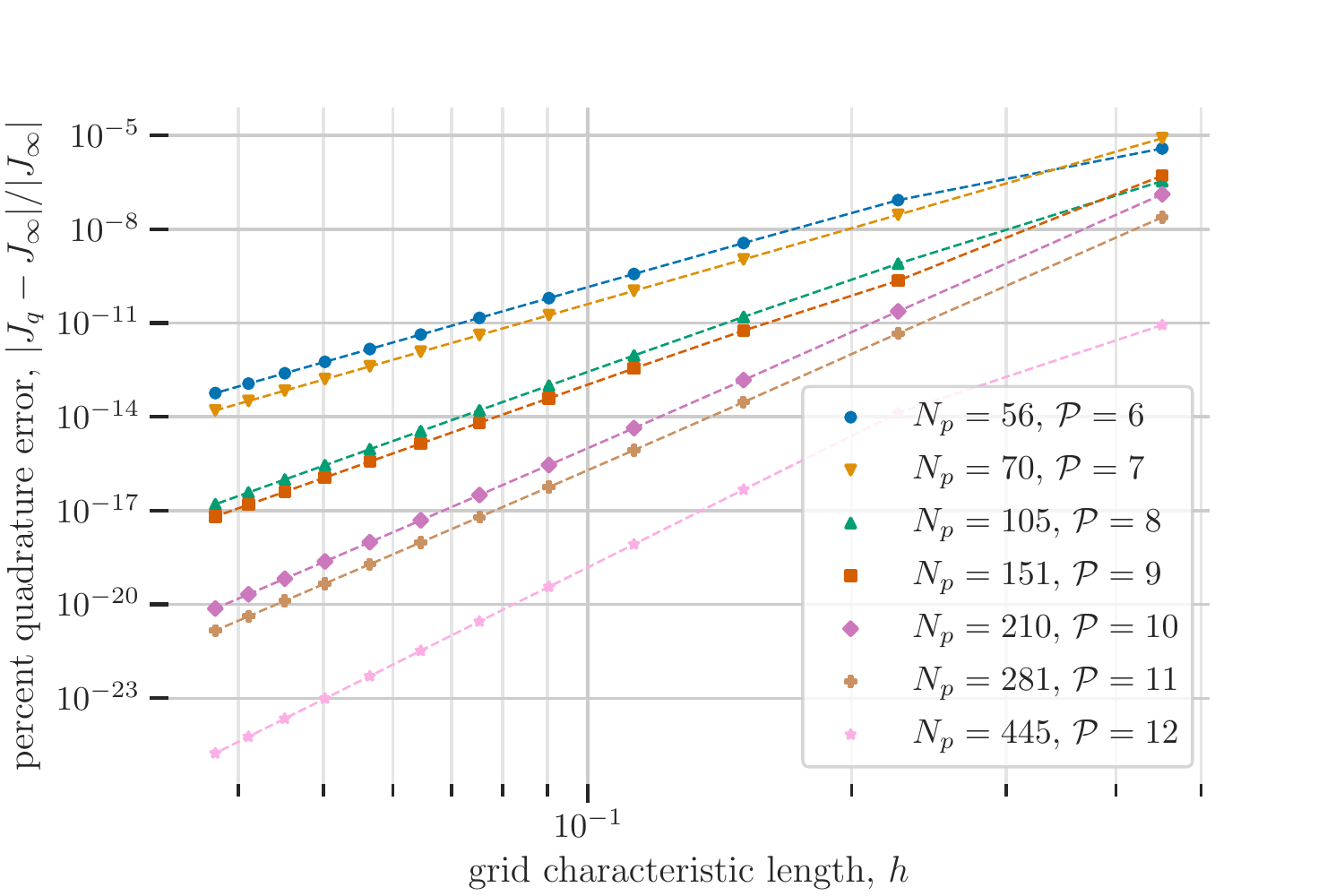}
\caption{Percent error in integration of $f_3$ on a pentatope grid due to quadrature.}
\label{quaderror_grid_pentatope_odd}
\end{figure}

Similarly, in Figures~\ref{quaderror_grid_prism_exp}, \ref{quaderror_grid_prism_sin}, and~\ref{quaderror_grid_prism_odd} we can see the results for the tetrahedral prism elements, and in Figures~\ref{quaderror_grid_pentatope_exp}, \ref{quaderror_grid_pentatope_sin}, and~\ref{quaderror_grid_pentatope_odd}, the results for the pentatope elements. The conclusion from these collected results is that quadrature error seems to be well behaved for each of the included transcendental functions, and approximations of integrals of similar functions by these quadrature rules can be expected to be accurate. In each case, quadrature rules of design order $\qorder$ produce convergence rates that are at least $h^{\qorder}$ and frequently better ($h^{\qorder + 1}$ or $h^{\qorder + 2}$), which makes them well suited for finite element methods. 

We note that, in a few cases lower degree rules outperform higher degree rules. This is likely due to a coincidental and beneficial alignment between the symmetry groups of a lower degree rule and the topology of one of our transcendental functions. Evidently, this trend will not hold for arbitrary transcendental functions, and the stronger quadrature rule will tend to outperform the weaker rule in most cases. 


\clearpage

\section{Conclusion}
\label{sec;conclusion}

In this article, we have explored the numerical foundations of finite element methods for real-world 3$d+t$ problems. We began our discussion by introducing a rigorous approach for generating four-dimensional elements using sequences of 0/1-polytopes. Thereafter, we demonstrated that at least one such sequence can be naturally encoded in terms of four digit binary numbers. Next, we narrowed our focus to consider the tesseract, tetrahedral prism, and pentatope elements. For each element type, we presented a brief review of the basic numerical properties, including explicit definitions of the reference domain, limits of integration, orthonormal polynomial basis, and symmetry groups. In addition, for each element, we constructed a new set of fully symmetric quadrature rules. Finally, we demonstrated the effectiveness of our quadrature rules by performing a detailed set of numerical experiments. 

We anticipate that our work (above) will significantly encourage the development of \emph{general} space-time methods for practical applications, and will not merely be limited to finite element methods. In particular, the newly developed element sequences will be useful in the construction of hybrid four-dimensional meshes for space-time visualization and data manipulation purposes. In addition, the newly developed quadrature rules will be immediately useful for both space-time finite element \emph{and} finite volume methods.

\section*{Declaration of Competing Interests}

The authors declare that they have no known competing financial interests or personal relationships that could have appeared to influence the work reported in this paper.

\section*{Acknowledgements}

The authors would like to thank Dr.~Chuluunbaatar Ochbadrakh for helping error-check the quadrature rules.

\section*{Funding}

This research did not receive any specific grant from funding agencies in the public, commercial, or not-for-profit sectors.

\pagebreak
\clearpage


{\footnotesize\bibliography{technical-refs}}

\end{document}